\documentclass[11pt,a4paper,reqno]{amsart}
\usepackage[T1]{fontenc}
\usepackage[utf8]{inputenc}
\usepackage{lmodern,microtype,amsmath,amssymb,amsthm,mathtools,mathrsfs}
\usepackage[margin=27mm,headheight=16pt]{geometry}
\usepackage{enumitem,booktabs,longtable,tabularx,array,xcolor,colortbl,hyperref,tikz-cd,needspace}
\definecolor{ink}{HTML}{193E55}
\hypersetup{colorlinks=true,linkcolor=ink,citecolor=ink,urlcolor=ink,pdftitle={Fermat curves and explicit class field theory of real quadratic fields},pdfauthor={Bora Yalkinoglu},pdfkeywords={Fermat curves, Jacobi sums, ray class fields, Kummer theory, Galois descent}}
\setlist{itemsep=3pt,topsep=4pt,leftmargin=1.8em}
\numberwithin{equation}{section}
\newtheorem{theorem}{Theorem}[section]
\newtheorem{proposition}[theorem]{Proposition}
\newtheorem{lemma}[theorem]{Lemma}
\newtheorem{corollary}[theorem]{Corollary}
\theoremstyle{definition}\newtheorem{definition}[theorem]{Definition}
\newtheorem{example}[theorem]{Example}
\newtheorem{convention}[theorem]{Convention}
\theoremstyle{remark}\newtheorem{remark}[theorem]{Remark}
\newcommand{\Q}{\mathbb Q}\newcommand{\Z}{\mathbf Z}
\newcommand{\F}{\mathbf F}\newcommand{\C}{\mathbb C}
\newcommand{\OO}{\mathcal O}\newcommand{\Gal}{\operatorname{Gal}}
\newcommand{\Hom}{\operatorname{Hom}}\newcommand{\Spec}{\operatorname{Spec}}
\newcommand{\Cl}{\operatorname{Cl}}\newcommand{\Tr}{\operatorname{Tr}}
\newcommand{\Nm}{\operatorname N}
\newcommand{\Frob}{\operatorname{Frob}}\newcommand{\divi}{\operatorname{div}}
\newcommand{\id}{\operatorname{id}}
\newcommand{\Jac}{\operatorname{Jac}}
\newcommand{\End}{\operatorname{End}}\newcommand{\Div}{\operatorname{Div}}
\newcommand{\lcm}{\operatorname{lcm}}
\newcommand{\im}{\operatorname{im}}

\newcommand{\Sel}{\operatorname{Sel}}

\newcommand{\Art}{\operatorname{Art}}
\newcommand{\rank}{\operatorname{rank}}

\newcommand{\geoFrac}[1]{\left\{#1\right\}}
\title[Explicit class field theory of real quadratic fields]{Fermat curves and explicit class field theory\protect\\for real quadratic number fields}
\author{Bora Yalkinoglu}
\address{Institut de Recherche Math\'ematique Avanc\'ee,
Universit\'e de Strasbourg, CNRS, UMR 7501,
7 rue Ren\'e Descartes, 67000 Strasbourg, France}
\email{yalkinoglu@math.unistra.fr}
\date{October 30, 2026}
\begin{document}
\begin{abstract}
\noindent We use Fermat curves with real multiplication and classical CM theory to construct an explicit class field theory for arbitrary real quadratic number fields. The underlying strategy is well-established in computational class field theory, our main contribution is a geometric model to realize this strategy explicitly. This approach complements the recent breakthrough works solving conjectures by Brumer, Stark and Shintani and gives an alternative constructive solution to Hilbert's 12th problem for real quadratic number fields.
\end{abstract}
\maketitle
\begingroup\small\setlength{\parskip}{0pt}\linespread{1}\selectfont
\setcounter{tocdepth}{1}\tableofcontents\endgroup
\clearpage

\section{Introduction}
\label{root170:introduction}

\noindent Let $K$ be a real quadratic number field of discriminant
$\Delta$, and let $K_{\mathfrak f}^+$ be a strict ray class field of
$K$. The goal of this paper is to construct explicit generators
for $K_{\mathfrak f}^+$ and thus for the maximal abelian extension
$K^{\rm ab}=\bigcup_{\mathfrak f}K_{\mathfrak f}^+$ of $K$. \\ \\
The underlying mechanism of exploiting Kummer theory is well-established in
computational class field theory \cite{CohenAdvanced,FiekerArtin,FHS}.
Our main innovation is to produce explicit Kummer bases from
Frobenius eigenvalues of Fermat Jacobians with real multiplication by $K$. \\ \\
\noindent For this we fix a character $\chi : \mathrm{G}^+_\mathfrak f \to \Z/n\Z$ of the strict ray class group, which defines a subfield $H_\chi \subset K^+_\mathfrak f$. It will be enough to consider characters $\chi$ landing in $\Z/p^a\Z$, for $p$ a prime number, and our construction works most elegantly if $p$ is ordinary, i.e., $p \nmid 2\Delta$.\\ \\ %(Proposition~\ref{v221:computed-characters} and Lemma~\ref{root170:assembly}). \\ \\      
In the ordinary case, for anti-invariant $\chi$, we will construct a finite module of Kummer classes $W_\Sigma$ and will show that the Kummer element $[\beta_\chi]$, which will produce a generator $\theta_\chi \in H_\chi$, is contained in $W_\Sigma$ (Theorem~\ref{v203:ordinary-support}). In the general case, we will construct a finite module of Kummer classes $W_\Sigma^{(1)}$ and will show that the Kummer element $[\beta_\chi]$ is contained in $W_\Sigma^{(1)}$ (with $\chi$ replaced by the signature-corrected character $\psi$ at $p=2$; Theorem~\ref{exc170:allcoverage}). \\ \\
\noindent The main contribution of this paper is to produce explicit Kummer bases for the Kummer modules $W_\Sigma$ and $W_\Sigma^{(1)}$ from Frobenius eigenvalues of Fermat Jacobians with real multiplication by $K$, which will give rise to an explicit description of the Kummer classes $[\beta_\chi]$. \\ \\
In the ordinary construction, Jacobi sums explicitly evaluate the Frobenius determinants. In the general construction, this remains true at good primes; at bad primes, 
the Shimura--Taniyama formula from CM theory is used to construct explicit replacement elements yielding the same Kummer module (Theorems~\ref{geo170:ordinary-complete} and~\ref{geo170:all-directions}). \\ \\
\noindent An alternative construction using Fermat curves without requiring real multiplication by $K$ is also given, which allows to use Jacobi sums at all primes, and produces the same Kummer module $W_\Sigma^{(1)}$ and the same Kummer class $[\beta_\chi]$ (after restoring the signature factor at $p=2$; Theorem~\ref{v228:diagonal-span}). 

\noindent Schematically our construction is given by
\[
\begin{tikzcd}[column sep=2em,row sep=3.2em]
 \underbrace{(K,\mathfrak f,\chi)}_{\text{arithmetic input}}
       \arrow[d,"\text{Fermat geometry}"'] & & & \\
 \underbrace{(A_i)_i}_{\substack{\text{Frobenius eigenvalues}\\
                              \text{or CM replacements}}}
       \arrow[r,"\text{saturation}"]
 &\underbrace{([z_j])_j}_{\substack{\text{Kummer basis}\\
                                  \text{}}}
       \arrow[r,"\text{reciprocity}"]
 &\underbrace{[\beta_\chi]}_{\substack{\text{Kummer class of $\chi$}\\
                                    }}
       \arrow[r,"\text{descent}"]
 &\underbrace{\theta_\chi}_{\text{generator of }H_\chi/K}
\end{tikzcd}
\]
\noindent and leads to our main result, Theorem \ref{root170:equality}, which is a constructive solution to Hilbert's twelfth problem for arbitrary real quadratic fields.
\Needspace{7\baselineskip}
\begin{theorem}\label{root170:main}
\noindent Our construction provides explicit elements
$\theta_1,\ldots,\theta_s$ such that
\[
 K_{\mathfrak f}^+=K(\theta_1,\ldots,\theta_s),
\]
\noindent together with their defining equations and Galois action.
\end{theorem}
\noindent Let's explain our approach in some more detail now.

\Needspace{7\baselineskip}
\subsection{The Kummer modules $W_\Sigma$ and $W_\Sigma^{(1)}$}
\noindent Fix $n=p^a$ and suppose $p\nmid2\Delta$. Let $E=K(\zeta_n)$,
let $j$ be complex conjugation on $E$, and let $\tau$ be the automorphism of $E$ that extends the nontrivial automorphism $\tau_K$ of $K$ while fixing $\zeta_n$. With the support $\Sigma$ chosen as in
Theorem~\ref{v203:ordinary-support}, the ordinary Kummer module is defined by
\begin{equation}
 W_\Sigma\simeq U_\Sigma/U_\Sigma^n,\qquad
 U_\Sigma=\{x\in\OO_E[\Sigma^{-1}]^\times:
                         x\tau(x)=xj(x)=1\}.
\end{equation}
\noindent For an (anti-invariant) character $\chi$, Theorem~\ref{v203:ordinary-support} provides a unique Kummer class $[\beta_\chi]\in W_\Sigma$, such that
\begin{equation}\label{v230:intro-target}
 \nu_{\beta_\chi}=\chi|_{G_E},\qquad E(\sqrt[n]{\beta_\chi})=EH_\chi.
\end{equation}
\noindent Here $\nu$ is the Kummer correspondence of Proposition~\ref{v203:kummer}.

\noindent In the general case, allowing $p\mid2\Delta$, we use
the field and support of Theorem~\ref{exc170:allcoverage} and the general Kummer module is defined by
\begin{equation}
 W_\Sigma^{(1)}\simeq U_\Sigma^{(1)}/(U_\Sigma^{(1)})^n,
 \qquad U_\Sigma^{(1)}=
       \{x\in\OO_E[\Sigma^{-1}]^\times:xj(x)=1\}.
\end{equation}
\noindent For odd $p$ it contains $[\beta_\chi]$ for every $\chi$.
At $p=2$, we first obtain $[\beta_\psi]\in W_\Sigma^{(1)}$ for the signature-corrected character $\psi$, then recover $[\beta_\chi]$; see Sections~\ref{exc170:coefficientfields}
and~\ref{v222:general-selection-section}. \\ \\
Once these modules are known, Artin reciprocity selection and Galois descent
are the standard mechanisms to extract the corresponding generator $\theta_\chi \in H_\chi$, see Sections~\ref{root170:selection}, \ref{root170:ordinary-descent}, \ref{v222:general-selection-section} and~\ref{root170:general-algebra}, and \cite[Sections~4.2--4.3]{FHS}. 

\subsection{The Fermat geometry}
\noindent We can now explain our main contribution, an explicit geometric construction of a Kummer basis for the Kummer modules $W_\Sigma$ and $W_\Sigma^{(1)}$. The construction is based on Fermat curves with real multiplication by $K$ and their Frobenius eigenvalues. \\
Again, we first explain the ordinary case, and then the general case. \\ \\
\noindent For $p\nmid2\Delta$, consider the smooth projective curve $C_{n,\Delta}$ defined by
\begin{equation} 
                       C_{n,\Delta}:Y^n+Z^\Delta=1.
\end{equation}
\noindent It is a quotient of an ordinary Fermat curve, its commuting deck actions give cyclotomic CM structures.
The inclusion $K\subset\Q(\zeta_\Delta)$ further gives rise to a RM structure. The two embeddings of $K$ specify equal-dimensional eigenspaces $M_v^+,M_v^-$
of the rational Tate module of $\Jac(C_{n,\Delta})$ (see section~\ref{geo170:rm}). \\ \\For $F=\Q(\zeta_{n\Delta})$, $\mathfrak q$ a prime of $F$ above $\Sigma$, and
$v\in\Z_n^\times$, put
\begin{equation}
 A_{\mathfrak q}(v)
   =\frac{\det(\Frob_{\mathfrak q}\mid M_v^+)}
          {\det(\Frob_{\mathfrak q}\mid M_v^-)}
   =\prod_{c\in\Z_\Delta^\times}
                J_{\mathfrak q,c}(v)^{\chi_\Delta(c)}\in U_\Sigma.
\end{equation}
\noindent Here $\chi_\Delta$ is the quadratic Dirichlet character of
$K$, and $J_{\mathfrak q,c}(v)$ is the Jacobi sum of
\eqref{geo170:jacobi}. CM describes the individual eigenvalues;
RM specifies the quotient with the two required norm conditions
(Theorem~\ref{geo170:detquotient}).

\begin{theorem}\label{v213:intro-ordinary}
\noindent Write $\Sigma_F$ for the primes of $F$ above $\Sigma$, and put
\[
 \Gamma_{\rm ord}=\langle\mu(E)\cap U_\Sigma,
      A_{\mathfrak q}(v):\mathfrak q\in\Sigma_F,
                                  \ v\in\Z_n^\times\rangle.
\]
\noindent Then $[U_\Sigma:\Gamma_{\rm ord}]<\infty$ and 
\[
 W_\Sigma\simeq
       \Gamma_{\rm ord}^{\rm sat}/(\Gamma_{\rm ord}^{\rm sat})^n,
\]
\noindent where
\[
 \Gamma_{\rm ord}^{\rm sat}
   =\{x\in U_\Sigma:x^{p^b}\in\Gamma_{\rm ord}
                                      \text{ for some }b\ge0\}.
\]
\end{theorem}

\noindent The first assertion is Theorem~\ref{geo170:ordinary-complete},
whose proof reduces to nonvanishing of Dirichlet $L$-values;
the second is Corollary~\ref{root170:saturation-theorem}.
Saturation extracts the missing $p$th roots inside $E$ and gives
basis representatives explicitly related to the geometric inputs:
\begin{equation}\label{v215:intro-jacobi-coordinates}
             z_j^{p^{b_j}}=\xi_j\prod_i A_i^{m_{ji}},
                 \qquad \xi_j\in\mu(E),\quad m_{ji}\in\Z.
\end{equation}
\noindent The saturation procedure is unnecessary when the index is prime to $p$
(Corollary~\ref{v224:no-saturation}).
The underlying Fermat decomposition and Frobenius-Jacobi formulas
are classical \cite{Weil1949,Weil1952,KoblitzRohrlich,ColemanFermat,OY}.

\noindent For arbitrary $p$, we enlarge the first exponent and consider the curve $C_{M,\Delta}$ which is defined by
\[
 C_{M,\Delta}:Y^M+Z^\Delta=1,\qquad
 \Delta\mid M,\quad n\mid M,\quad E\subseteq F=\Q(\zeta_M).
\]
\noindent The map $(Y,Z)\mapsto(Y^{M/n},Z)$ connects this curve to
$C_{n,\Delta}$, preserving the same real multiplication. \\ 
For $v\in\Z_M\setminus\{0\}$, the eigenspaces $M_v^+,M_v^-$ can
now have different dimensions;
instead we compare $M_v^+$ with its conjugate $M_{-v}^+$.
For each rational prime $\ell$ below $\Sigma$, choose a prime
$\mathfrak q_\ell$ of $F$ above $\ell$ and a good-reduction
extension with a prime $\mathfrak P_\ell$ above $\mathfrak q_\ell$.
Then Theorem~\ref{v223:uniform-determinant} gives
\begin{equation}\label{v222:intro-RM-values}
 A_{\mathfrak q_\ell}(v)=N_{F/E}\!\left(
   \frac{\det(\Frob_{\mathfrak P_\ell}\mid M_v^+)}
        {\det(\Frob_{\mathfrak P_\ell}\mid M_{-v}^+)}
                                  \right)\in U_\Sigma^{(1)}.
\end{equation}
\noindent The determinants are identified with their algebraic values
in $F$. The norm puts their quotient in $E$; it is unnecessary in
the ordinary construction because its determinants already lie in $E$
(Proposition~\ref{v215:ordinary-general-comparison}).

\noindent Here two kinds of primes must be distinguished.
The character-order prime $p$ is \emph{exceptional} when
$p\mid2\Delta$. Good or bad reduction concerns the prime $\ell$
at which Frobenius is evaluated. In the ordinary construction,
$\Sigma$ avoids $p\Delta$, so every evaluation prime lies in the
good-reduction range $\ell\nmid n\Delta$. In the general construction,
the curve has good reduction at $\ell\nmid M$, but the support may
also contain primes dividing $M$, even when $p$ is ordinary.
Conversely, exceptional $p$ does not mean that every Frobenius
evaluation is at bad reduction.

\noindent At $\ell\nmid M$, the determinants are again products of
Jacobi sums. At primes dividing $M$, the Shimura--Taniyama formula
gives their ideals, from which principal-ideal arithmetic constructs
explicit replacements $A'_{\mathfrak q_\ell}(v)$
(Proposition~\ref{v223:principalized-input-proposition}).
These need not equal the Frobenius values, but give the same
Kummer module after saturation.

\begin{theorem}\label{v222:intro-uniform}
\noindent For the field and support of Theorem~\ref{exc170:allcoverage},
and the choices above, put
\[
 \Gamma^{(1)}=\langle\mu(E),A_{\mathfrak q_\ell}(v):
       1\le v<M,\ \ell\text{ lies below }\Sigma\rangle.
\]
\noindent Then $[U_\Sigma^{(1)}:\Gamma^{(1)}]<\infty$ and
\[
 W_\Sigma^{(1)}\simeq
       \Gamma^{(1),\rm sat}/(\Gamma^{(1),\rm sat})^n,
\]
\noindent where saturation is taken inside $U_\Sigma^{(1)}$ at $p$.
The same assertions hold for the subgroup obtained by replacing
the bad-prime inputs with $A'_{\mathfrak q_\ell}(v)$.
\end{theorem}
\noindent These are Theorem~\ref{geo170:all-directions} and
Proposition~\ref{exc170:saturation}. The proof combines classical
CM theory \cite{SerreTate1968,TateCM1969} with nonvanishing of
Dirichlet $L$-values at one. \\ \\
Alternatively, the smooth projective curves $C_c$, for $c > 2$ and $c \, \vert \, M$, defined by
\begin{equation}
      C_c:z^c=x(1-x),
\end{equation}
 which are also quotients of ordinary Fermat curves, give rise to the same Kummer
module $W_\Sigma^{(1)}$ using Jacobi sums at all support primes $\ell$, choosing $c$ with $\ell\nmid c$ for each contributing character component. This construction does not require RM by $K$
(Theorem~\ref{v228:diagonal-span}). \\ \\
\noindent To finish, Artin reciprocity selects the target class $[\beta_\chi]$ from
the Kummer basis, and equivariant Galois descent gives the generator $\theta_\chi \in H_\chi$, which now has a clear geometric interpretation in terms of Fermat geometry.
This is worked out in
Theorems~\ref{root170:reciprocity-proved},
\ref{root170:ordinary-descent-theorem},
\ref{v216:general-radicand} and~\ref{root170:general-R}. \\ \\
We also have readers in mind who have a background in mathematical physics, as the recent breakthrough work \cite{RadchenkoWheeler2026} uses mainly ideas from mathematical physics, therefore we try to give many details which are not strictly necessary for the number-theoretic reader. We hope that this will make our paper accessible to a wider audience.

\subsection*{Motivation and background}
\noindent
Our solution to Hilbert's 12th problem for real quadratic fields is complementary to the recent breakthrough works of Dasgupta-Kakde and Radchenko-Wheeler, which solve conjectures of Brumer, Stark and Shintani \cite{DasguptaKakde2024,RadchenkoWheeler2026}.
Our construction is explicitly geometric, and is based on the classical geometry of Fermat curves and the theory of Complex Multiplication \cite{ShimuraTaniyama,ShimuraCM} pepped up with the RM structure coming from $K$. Our generators are not defined as special values of Artin $L$-functions or other special functions.
\\ \\ \noindent
Nevertheless, our construction is heavily inspired by earlier work of the author which relates Shintani-Stark invariants and cyclic quantum
dilogarithms \cite{YalkinogluShintani,YalkinogluCyclic2025}. The latter special functions live on Fermat curves and playing around with these structures, including extensive calculations with Sage and Pari, led to the present paper. Related special-value constructions include
\cite{DasguptaKakde2024,RadchenkoWheeler2026,Yamamoto2010,
KoppShintaniFaddeev,BergeronCharolloisGarcia}, Fermat towers and $p$-adic beta-gamma formulas appear in \cite{AndersonHyperadelic,BaldassarriBetaGamma}. It is an extremely interesting problem to understand the relation between these different approaches to Hilbert's 12th problem for real quadratic fields. \\ \\
The main mystery of whether there is a geometric replacement for elliptic curves with CM for a real quadratic number field remains completely open!

\subsection*{Outline of the paper}
\noindent Section~\ref{root170:target} recalls the character and Kummer
correspondences (Proposition~\ref{v203:kummer}).
Section~\ref{root170:ordinary-space} defines the Kummer modules
(Propositions~\ref{v203:Selmer-exact} and~\ref{root170:U-W});
Section~\ref{v206:coverage} proves target membership
(Theorems~\ref{v203:ordinary-support} and~\ref{exc170:allcoverage}).

\noindent In the ordinary case, Section~\ref{geo170:curves} constructs the
Frobenius determinants (Theorem~\ref{geo170:detquotient});
Section~\ref{geo170:ordinary-spanning} proves their divisor span
(Theorem~\ref{geo170:ordinary-complete});
Section~\ref{root170:saturation} gives the Kummer basis
(Corollary~\ref{root170:saturation-theorem});
Sections~\ref{root170:selection} and~\ref{root170:ordinary-descent}
perform reciprocity selection and descent
(Theorems~\ref{root170:reciprocity-proved}
and~\ref{root170:ordinary-descent-theorem}).

\noindent In the general case, Section~\ref{geo170:universal} defines the
inputs (Theorem~\ref{v223:uniform-determinant} and
Proposition~\ref{v223:principalized-input-proposition});
Section~\ref{v222:general-spanning} proves their span and the
diagonal alternative (Theorems~\ref{geo170:all-directions}
and~\ref{v228:diagonal-span});
Section~\ref{exc170:classnumber} gives the basis
(Proposition~\ref{exc170:saturation} and
Corollary~\ref{v226:general-basis});
Sections~\ref{v222:general-selection-section}
and~\ref{root170:general-algebra} perform selection and descent
(Theorems~\ref{v216:general-radicand} and~\ref{root170:general-R}).

\noindent Section~\ref{mx171:chapter} assembles the fields
(Theorem~\ref{root170:equality}), gives their Artin action
(Proposition~\ref{v211:artin-action}), and treats mixed moduli
and the strict Hilbert class field
(Propositions~\ref{v216:mixed-modulus} and~\ref{hi171:hilbertcorrectness}).
Section~\ref{root170:examples} gives explicit examples.
Appendix~\ref{v207:ray-character-appendix} records the ray-character
calculation (Proposition~\ref{v204:arithmetic-input-proposition});
Appendix~\ref{geo203:CM-appendix} relates the CM ideal to the reflex
norm. Appendix~\ref{mx201:towers} explains a useful classical
feature of real quadratic class fields: at fixed support, only
cyclotomic growth continues indefinitely. More precisely,
Theorem~\ref{mx201:stabilization} gives
\[
 \bigcup_{a\ge1}K_{\mathfrak f^a}^+
    =K_{\mathfrak f^{a_0}}^+\,
      K\bigl(\Q_\infty^{(p)}:p\in T\bigr)
       \quad\text{for some }a_0,
\]
\noindent where $T$ consists of the rational primes whose full support in
$K$ divides $\mathfrak f$, and $\Q_\infty^{(p)}$ is the
cyclotomic $\widehat\Z_p$-extension of $\Q$.

\subsection*{Notation and conventions}
\label{v207:global-conventions}
\noindent Fix an algebraic closure $\overline\Q\subset\C$ of $\Q$;
all number fields lie in $\overline\Q$. For a number field $F$,
write $G_F=\Gal(\overline\Q/F)$ and let $\mu(F)$ be its roots
of unity. Field products denote composita; products of finite
algebras denote direct products.

\noindent Frobenius is arithmetic: if $L/F$ is a finite Galois extension,
$\mathfrak q$ is unramified in $L$, and $\mathfrak P\mid\mathfrak q$, then
\[
 \Frob_{\mathfrak P/\mathfrak q}\in\Gal(L/F),\qquad
 \Frob_{\mathfrak P/\mathfrak q}(x)\equiv
 x^{\Nm(\mathfrak q)}\pmod{\mathfrak P}
 \quad(x\in\OO_L).
\]
\noindent Its conjugacy class is denoted $\Frob_{\mathfrak q}$; in an
abelian extension it is independent of $\mathfrak P$ and equals
$\Art_{L/F}(\mathfrak q)$ \cite[Section~2]{FHS}.
On a Tate module at good reduction, this is the action induced by
the $\Nm(\mathfrak q)$-power Frobenius endomorphism of the reduction,
as in~\eqref{geo170:frob-convention}.

\noindent Write $K=\Q(\sqrt d)$ with $d>0$ squarefree, and let $\tau_K$
be its nontrivial automorphism. Fix $\sigma_1(\sqrt d)>0$ and
$\sigma_2=\sigma_1\circ\tau_K$; $x\gg0$ means positivity at
both embeddings. Choose
$\OO_K^{\times,+}=\langle\epsilon\rangle$ with
$\sigma_1(\epsilon)>1$. All ray class fields are strict:
both real places are included in the modulus.

\noindent We use $\Z_n=\Z/n\Z$, reserving $\widehat\Z_p$ for the
$p$-adic integers. Characters take values in $\Z_n$; their
multiplicative form is $\rho_\chi(g)=\zeta_n^{\chi(g)}$.
Kummer classes are written multiplicatively, with scalar action
$a\cdot[x]=[x]^a$. The notation $W_\Sigma(E,n)$ or
$W_\Sigma^{(1)}(E,n)$ displays the field and exponent when needed;
otherwise they are fixed. The subscript $\Sigma$ records the prime
support. These are $\Z_n$-modules, and vector
spaces when $n$ is prime. For an abelian group $A$ written additively,
\[
 A[n]=\{a\in A:na=0\},\qquad
 A[p^\infty]=\bigcup_{r\ge1}A[p^r].
\]

\noindent Prime ideals use fraktur letters. Valuations satisfy
$v_{\mathfrak p}(\pi_{\mathfrak p})=1$ and
$v_{\mathfrak p}(0)=+\infty$.
For a nonzero integral ideal $\mathfrak a\subset\OO_F$,
$\Nm(\mathfrak a)=|\OO_F/\mathfrak a|$ is its absolute norm,
extended multiplicatively to fractional ideals. Subscripts
indicate relative norms. For $B\in\Z_{>0}$, $v_p(B)$ is its
$p$-adic valuation.

\noindent A lattice is a finitely generated free abelian group.
Smith normal form computes, for an integer relation matrix
$D\in\operatorname{Mat}_{r\times s}(\Z)$, matrices
$P\in\operatorname{GL}_r(\Z)$ and
$Q\in\operatorname{GL}_s(\Z)$ such that
\[
 PDQ=\operatorname{diag}(d_1,\ldots,d_t,0),\qquad
 d_i>0,\quad d_i\mid d_{i+1}.
\]
\noindent The diagonal matrix may be rectangular. It gives explicit coordinates
\[
 \Z^r/D\Z^s\xrightarrow{\sim}
    \bigoplus_{i=1}^t\Z_{d_i}\oplus\Z^{r-t},
                    \qquad[x]\longmapsto[Px].
\]
\noindent The quotient is finite exactly when $t=r$, with order
$\prod_i d_i$. We use this standard algorithm for ray-group
coordinates and lattice indices \cite[Section~2.4.4]{CohenANT}.

\subsection*{Use of AI tools}
\noindent When we started this project at the beginning of 2026, we used SageMath and Pari extensively for explicit calculations and 
exploited AI tools to help with those calculations. Later, we used AI-tools to proof-read the manuscript, help with references and check the examples in section~\ref{root170:examples}. 
All mathematical ideas are due to the author, who takes responsibility for the content of this paper and the correctness of the results. 
\section{Ray class characters and their Kummer representatives}
\label{root170:target}

\noindent We fix conventions for standard class field and Kummer theory.
For ray groups and reciprocity see \cite[Chapter VI]{NeukirchANT}
and \cite[Part III]{NeukirchBonn}; for Kummer theory and its
constructive use see \cite[Theorem 5.30]{MilneFT},
\cite[Chapters 4--5]{CohenAdvanced}, and \cite[Sections 3--4]{FHS}.
The statements below are classical or elementary consequences.

\subsection{The ray class group and its characters}
\label{v204:arithmetic-input}
\noindent For a nonzero integral ideal
$\mathfrak f$, let $I_K(\mathfrak f)$ be the multiplicative group
of nonzero fractional ideals prime to $\mathfrak f$. The strict
ray class group is
\begin{equation}\label{root170:raygroup}
 G_{\mathfrak f}^+=I_K(\mathfrak f)/
       \{(x):x\in K^\times,\ x\gg0,\ x\equiv1\pmod{\mathfrak f}\}.
\end{equation}
\noindent The congruence $x\equiv1\pmod{\mathfrak f}$ means
$v_{\mathfrak p}(x-1)\ge v_{\mathfrak p}(\mathfrak f)$ for
$\mathfrak p\mid\mathfrak f$.

\noindent Recall the well-known exact sequence
\cite[Chapter VI, Section 1]{NeukirchANT}:
\begin{equation}\label{v214:ray-exact}
 1\longrightarrow
 \frac{(\OO_K/\mathfrak f)^\times}{\langle\bar\epsilon\rangle}
 \longrightarrow G_{\mathfrak f}^+
 \longrightarrow\Cl^+(K)\longrightarrow1.
\end{equation}
\noindent Here $\Cl^+(K) = G_{(1)}^+$ is the narrow ideal class group, and
$\bar\epsilon$ is the reduction of the fixed generator of totally
positive units modulo $\mathfrak f$.
The first map sends $\bar b$ to $[(b)]$ for a totally positive
lift, which exists by approximation; the second forgets the
congruence condition. Thus
\begin{equation}\label{v214:ray-quotient}
 \Cl^+(K)=1\quad\Longrightarrow\quad
 G_{\mathfrak f}^+\simeq
       (\OO_K/\mathfrak f)^\times/\langle\bar\epsilon\rangle.
\end{equation}
\noindent In general this quotient is only the kernel in
\eqref{v214:ray-exact}; the sequence need not split. No class-number
assumption is made here.

\noindent Artin reciprocity and the existence theorem identify this group with
the Galois group of the strict ray class field:
\[
 \Art_{\mathfrak f}:G_{\mathfrak f}^+
    \xrightarrow{\sim}\Gal(K_{\mathfrak f}^+/K),
       \qquad[\mathfrak p]\longmapsto\Frob_{\mathfrak p},
\]
\noindent where Frobenius acts on residue fields by $x\mapsto x^{\Nm(\mathfrak p)}$
\cite[Part III, Theorems 7.8, 9.2 and 9.4]{NeukirchBonn}. For a character
$\chi:G_{\mathfrak f}^+\to \Z_n$, put
\[
 H_\chi=(K_{\mathfrak f}^+)^{\ker\chi},
 \qquad [H_\chi:K]=|\im\chi|.
\]
\noindent We use the same symbol for its inflation to $G_K$:
\[
 \chi(g)=\chi\bigl(\Art_{\mathfrak f}^{-1}
                        (g|_{K_{\mathfrak f}^+})\bigr),\qquad g\in G_K.
\]
\noindent Thus $\chi|_{G_E}$ below always means restriction of this
continuous Galois character.

\subsection{Recovering the full ray class field}
\noindent Finite abelian characters separate points. Consequently
\begin{equation}\label{v214:character-compositum}
 K_{\mathfrak f}^+=\prod_\chi H_\chi,
 \qquad
 \Gal\bigl(K_{\mathfrak f}^+/\prod_\chi H_\chi\bigr)
             =\bigcap_\chi\ker\chi=\{1\},
\end{equation}
\noindent where $\chi$ ranges over $\Hom(G_{\mathfrak f}^+,\Z_e)$ with
$e$ the least common multiple of the orders of its elements
(the group exponent), and the product denotes the compositum
inside $\overline\Q$. It suffices to take just the coordinate
characters of a decomposition
$G_{\mathfrak f}^+\simeq\bigoplus_j \Z_{d_j}$ computed by Smith normal form.
Their prime-power components also separate points
\cite[Chapter 3]{MilneFT}\cite[Section 4.1]{FHS}; Section~\ref{mx171:chapter}
computes these characters and assembles their field generators.
Characters valued in one fixed $\Z_n$ suffice precisely when $e\mid n$.

\noindent Appendix~\ref{v207:ray-character-appendix} computes a presentation
$G_{\mathfrak f}^+\simeq\Z^r/A\Z^s$ from ideals, units and residues.
The columns of the integer matrix $A$ are the exponents in the
relations among the chosen ideal classes. Its Smith normal form
gives the cyclic factors and their coordinate characters
\cite[Section~2.4.4]{CohenANT}; no equation for $H_\chi$ is an input.

\subsection{Kummer theory after adjoining roots of unity}
\noindent For $n>1$, the field $E=K(\zeta_n)$ always suffices for Kummer theory. We allow
larger coefficient fields because the dyadic norm-one construction
of Section~\ref{exc170:coefficientfields} may require them.
Let $E/K$ be a finite Galois extension with $\mu_n\subset E$,
and choose a primitive $n$th root of unity $\zeta_n$.
Set $G_E=\Gal(\overline\Q/E)$. For $\beta\in E^\times$ and
$\alpha\in\overline\Q$ with $\alpha^n=\beta$, define
\begin{equation}\label{root170:kummer-character}
 \nu_\beta:G_E\longrightarrow \Z_n,
                  \qquad g(\alpha)=\zeta_n^{\nu_\beta(g)}\alpha.
\end{equation}
\noindent Since $g(\alpha)/\alpha\in\mu_n$, this determines
$\nu_\beta(g)\in\Z_n$ uniquely. The group $G_E$ fixes $\mu_n$,
so $\nu_\beta$ is a homomorphism, with open kernel
$G_{E(\alpha)}$. Replacing $\alpha$ by $\zeta_n^k\alpha$, or
$\beta$ by $\beta c^n$ with $c\in E^\times$, leaves the character
unchanged. The following standard Kummer isomorphism therefore
has the indicated domain \cite[Theorem 5.30]{MilneFT}.

\begin{proposition}
% Former title: Kummer correspondence
\label{v203:kummer}
\noindent The map
\begin{equation}\label{root170:kummer-iso}
 \nu:E^\times/E^{\times n}\xrightarrow{\sim}
 \Hom_{\rm cont}(G_E,\Z_n),\qquad[\beta]\longmapsto\nu_\beta
\end{equation}
\noindent is an isomorphism. If $[\beta]$ has order $m$, then
$E(\alpha)=\overline\Q^{\ker\nu_\beta}$ and $[E(\alpha):E]=m$.
\end{proposition}

\begin{definition}
% Former title: The reciprocity condition
\label{v203:Rdefinition}
\noindent Fix a ray class character $\chi:G_{\mathfrak f}^+\to\Z_n$,
inflated to $G_K$ as above. An element $\beta\in E^\times$
represents its restriction to $E$ when
\begin{equation}\label{root170:Rchi}
 \makebox[0.88\linewidth]{%
  $(\mathrm R_\chi)$\hfill
  $\displaystyle\nu_\beta=\chi|_{G_E}.$\hfill
  \phantom{$(\mathrm R_\chi)$}}
\end{equation}
\end{definition}
\noindent Equivalently, $[\beta]=\nu^{-1}(\chi|_{G_E})$.
We calculate this class from Fermat values. Distinct characters
can have the same restriction; descent uses the original $\chi$.
The isomorphism $\nu$ sends products of Kummer classes to sums of characters.

\begin{lemma}
% Former title: The power-residue test
\label{v203:residue-character}
\noindent Let $\beta\in E^\times$ be any candidate representative, and let
$\mathfrak P$ be a prime of $E$ such that $\mathfrak P\nmid n\mathfrak f\OO_E$
and $v_{\mathfrak P}(\beta)=0$, and put $q_{\mathfrak P}=|\OO_E/\mathfrak P|$. Then
\begin{align}
 \zeta_n^{\nu_\beta(\Frob_{\mathfrak P})}
   &\equiv\beta^{(q_{\mathfrak P}-1)/n}\pmod{\mathfrak P},
                                      \label{root170:residue-law}\\
 (\chi|_{G_E})(\Frob_{\mathfrak P})
   &=\chi([\Nm_{E/K}\mathfrak P]_{\mathfrak f}).
                                      \label{root170:norm-frob}
\end{align}
\noindent Consequently $(\mathrm R_\chi)$ is equivalent to
\begin{equation}\label{root170:residue-target}
 \beta^{(q_{\mathfrak P}-1)/n}
      \equiv\zeta_n^{\chi([\Nm_{E/K}\mathfrak P]_{\mathfrak f})}\pmod{\mathfrak P}
\end{equation}
\noindent for all but finitely many primes $\mathfrak P$.
\end{lemma}
\begin{proof}
\noindent The first formula is the standard Kummer Artin-symbol calculation
\cite[Section 3]{FiekerArtin}\cite[Proposition 4]{FHS}; reduction is
injective on $\mu_n$, so $n\mid q_{\mathfrak P}-1$. The second is
norm compatibility of reciprocity
\cite[Part III, Section 6, p.~153]{NeukirchBonn}.
Chebotarev makes equality on almost every Frobenius equivalent
to equality of characters \cite[Chapter VII, Theorem 13.4]{NeukirchANT}.
\end{proof}
\noindent Section~\ref{root170:selection} makes this test finite. Once
$\nu^{-1}(\chi|_{G_E})$ belongs to a known free Kummer
module, evaluation
of the Kummer pairing at separating Frobenius elements distinguishes
all its classes. Matching the prescribed Artin values there determines
the target class; we call this \emph{selection by Artin reciprocity}.

\noindent The next standard compatibility explains the equations for the
coefficients used in descent. It is an immediate consequence of
the Kummer correspondence \cite[Theorem 5.30]{MilneFT}.

\begin{lemma}
% Former title: Galois action on the Kummer representative
\label{v203:covariance}
\noindent Suppose $\beta\in E^\times$ satisfies $(\mathrm R_\chi)$.
For $\delta\in\Gal(E/K)$, let $1\le u_\delta<n$ be determined by
$\delta(\zeta_n)=\zeta_n^{u_\delta}$. Then
\begin{equation}\label{root170:covariance}
 \delta([\beta])=[\beta]^{u_\delta}
                   \quad\text{in }E^\times/E^{\times n}.
\end{equation}
\noindent Equivalently, there is $c_\delta\in E^\times$ such that
\[
                   c_\delta^n=\delta(\beta)/\beta^{u_\delta}.
\]
\end{lemma}
\begin{proof}
\noindent For a lift $\widetilde\delta\in G_K$ of $\delta$ and $g\in G_E$,
Kummer equivariance gives
$\nu_{\delta(\beta)}(g)=u_\delta\chi(\widetilde\delta^{-1}g\widetilde\delta)
=u_\delta\chi(g)$, since $\chi$ is a homomorphism to an abelian
group. Now apply the isomorphism~\eqref{root170:kummer-iso}.
\end{proof}
\noindent This is the power equation needed for a lift
$d_\delta(\alpha)=c_\delta\alpha^{u_\delta}$ with
$d_\delta|_E=\delta$. The descent calculation in
Sections~\ref{root170:ordinary-descent} and
\ref{root170:general-algebra} chooses the roots so that these maps
compose correctly and give the prescribed character $\chi$.

\subsection{Quadratic conjugation and the cyclotomic part}
\noindent Conjugation compatibility of the Artin map
\cite[Section 5.1, Proposition 15]{FHS} gives
\[
 \tau_{\mathfrak f}:G_{\mathfrak f}^+
      \xrightarrow{\sim}G_{\tau_K\mathfrak f}^+,
          \quad[\mathfrak a]\longmapsto[\tau_K\mathfrak a],
 \qquad
 \chi^\tau(g)=\chi(\tau_{\mathfrak f}(g))
                 \quad(\tau_K\mathfrak f=\mathfrak f).
\]
\begin{lemma}
% Former title: Invariant characters are cyclotomic

\label{v206:invariant-cyclotomic}
\noindent For every $n$, if $\tau_K\mathfrak f=\mathfrak f$ and
$\chi^\tau=\chi$, then $H_\chi/\Q$ is abelian. In particular
$H_\chi$ is contained in a cyclotomic field.
\end{lemma}
\begin{proof}
\noindent The stable kernel makes $H_\chi/\Q$ Galois, and invariance
makes $\Gal(H_\chi/K)$ central. A central subgroup with cyclic
quotient generates an abelian group together with any lift of the
quotient generator. Apply Kronecker--Weber
\cite[Part III, Theorem 7.11]{NeukirchBonn}.
\end{proof}

\begin{lemma}
% Former title: Reduction of an arbitrary odd-order character

\label{root170:assembly}
\noindent Suppose $n$ is odd. Put
$\mathfrak f^*=\operatorname{lcm}(\mathfrak f,\tau_K\mathfrak f)$,
and pull $\chi$ back along the natural surjection
$G_{\mathfrak f^*}^+\to G_{\mathfrak f}^+$. On this stable modulus set
\[
 \chi^+=2^{-1}(\chi+\chi^\tau),\qquad
 \chi^-=2^{-1}(\chi-\chi^\tau).
\]
\noindent Then $\chi=\chi^++\chi^-$ and
$(\chi^\pm)^\tau=\pm\chi^\pm$. The field $H_{\chi^+}$ is
cyclotomic. Constructing $H_{\chi^-}$ and their compositum
$L=H_{\chi^+}H_{\chi^-}$ determines the original field by
\[
 H_\chi=L^{\ker s},\qquad
 s:\Gal(L/K)\longrightarrow \Z_n,
        \quad g\longmapsto\chi^+(g|_{H_{\chi^+}})
                          +\chi^-(g|_{H_{\chi^-}}).
\]
\end{lemma}
\begin{proof}
\noindent The ray-group surjection follows from approximation. The character
identities follow by division by two. Lemma~\ref{v206:invariant-cyclotomic}
applies to $\chi^+$, and the fixed-field formula is Galois
correspondence applied to $\chi=\chi^++\chi^-$.
\end{proof}
\noindent A zero component contributes only $K$; each nonzero component
may be constructed at its actual order.

\begin{remark}
% Former title: A finite cyclotomic search
\label{v208:cyclotomic-bound}
\noindent One can show that the odd-order invariant field satisfies
\[
 H_{\chi^+}\subseteq\Q(\zeta_{\Delta N^*}),\qquad
 N^*=|\OO_K/\mathfrak f^*|\in\Z_{>0}.
\]
\end{remark}

%\noindent At two, the general construction treats the entire character,
%without division by two. Conjugation invariance is distinct from
%\emph{total evenness}, meaning vanishing on complex conjugation
%at both real places. Every odd-order character is totally even;
%at two this condition is arranged in Section~\ref{exc170:overview}.

\section{Kummer modules}
\label{root170:ordinary-space}

\noindent We define the central Kummer modules $W_\Sigma$ and $W_\Sigma^{(1)}$. Everything in this section is standard, see, e.g., 
\cite[Sections~5.2.2 and~7.4]{CohenAdvanced}.

\subsection{Supported classes and the choice of primes}
\noindent For the ideal calculations in this subsection, $E$ is any number
field, $n\ge2$, and $\Sigma$ is any finite set of its finite
primes $\mathcal P_E$. No class-group assumption is made. Define
\[
\begin{split}
 R_\Sigma&=\OO_E[\Sigma^{-1}]
       =\{x\in E:v_{\mathfrak p}(x)\ge0\text{ for }{\mathfrak p}\notin\Sigma\},\\
 R_\Sigma^\times
   &=\{x\in E^\times:v_{\mathfrak p}(x)=0\text{ for every }{\mathfrak p}\notin\Sigma\},\\
 \Div(E)&=\bigoplus_{{\mathfrak p}\in\mathcal P_E}\Z[{\mathfrak p}],\qquad
 \Z[\Sigma]=\bigoplus_{{\mathfrak p}\in\Sigma}\Z[{\mathfrak p}]\subset\Div(E),\\
 \divi:E^\times&\longrightarrow\Div(E),
                  \qquad x\longmapsto\sum_{\mathfrak p}v_{\mathfrak p}(x)[{\mathfrak p}].
\end{split}
\]
\noindent The elements of $R_\Sigma^\times$ are the \emph{$\Sigma$-units},
and $\Z[\Sigma]$ is the group of divisors supported on $\Sigma$. In particular
$\divi(R_\Sigma^\times)\subset\Z[\Sigma]$. For an automorphism
$\sigma$ of $E$ preserving $\Sigma$, the actions are
\[
\begin{gathered}
 x\longmapsto\sigma(x),\qquad
 [x]\longmapsto[\sigma(x)]\text{ in }E^\times/E^{\times n},\\
 \sum_{\mathfrak p}m_{\mathfrak p}[{\mathfrak p}]\longmapsto\sum_{\mathfrak p}m_{\mathfrak p}[\sigma {\mathfrak p}],
 \qquad\divi(\sigma x)=\sigma\divi(x).
\end{gathered}
\]

\begin{definition}
% Former title: The $n$-Selmer group
\label{v203:Selmer}
\noindent The $n$-Selmer group with support $\Sigma$ is
\begin{equation}\label{root170:supported}
 \Sel_n(E,\Sigma)=
 \{[x]\in E^\times/E^{\times n}:
                  v_{\mathfrak p}(x)\in n\Z\text{ for }{\mathfrak p}\notin\Sigma\}.
\end{equation}
\end{definition}
\noindent An $R_\Sigma$-unit has zero divisor outside $\Sigma$; a Selmer
representative need only have an $n$-divisible divisor there.
Write $\Cl(R_\Sigma)$ for the group of fractional
$R_\Sigma$-ideals modulo principal ideals.
\begin{proposition}
% Former title: The Kummer--Selmer exact sequence

\label{v203:Selmer-exact}
\noindent There is an exact sequence
\begin{equation}\label{root170:selmer-exact}
 1\longrightarrow R_\Sigma^\times/(R_\Sigma^\times)^n
 \xrightarrow{u\mapsto[u]}\Sel_n(E,\Sigma)
 \xrightarrow{\partial}\Cl(R_\Sigma)[n]\longrightarrow0,
\end{equation}
\noindent where
\[
 \partial[x]=[\mathfrak a_x],\qquad
 \mathfrak a_x=\prod_{{\mathfrak p}\notin\Sigma}
                   (\mathfrak pR_\Sigma)^{v_{\mathfrak p}(x)/n},\qquad
                         xR_\Sigma=\mathfrak a_x^n.
\]
\noindent The obstruction to a unit representative of $[x]$ is exactly
$\partial[x]$. For $n=p^a$ it belongs to
$\Cl(R_\Sigma)[p^a]$, the subgroup of ideal classes killed by
$p^a$. Every supported class has a unit representative if and
only if this subgroup is zero, equivalently
$\Cl(R_\Sigma)[p]=0$.
\end{proposition}
\begin{proof}
\noindent This is the usual Selmer sequence
\cite[Proposition~5.2.8]{CohenAdvanced}, with the same ideal proof
over $R_\Sigma$ for arbitrary $n$: $\partial[x]=0$ precisely when
$\mathfrak a_x=yR_\Sigma$, in which case $x/y^n$ is a unit.
Every $n$-torsion ideal class has $\mathfrak a^n=xR_\Sigma$;
an $n$th root of a unit has zero valuations outside $\Sigma$.
\end{proof}

\begin{lemma}
% Former title: Enlarging the support to kill the class group

\label{v203:class-support}
\noindent For any number field $E$ and disjoint finite sets
$\Sigma_0,\mathcal B\subset\mathcal P_E$, there is a finite
$\Sigma\supseteq\Sigma_0$ disjoint from $\mathcal B$ such that
$\Cl(R_\Sigma)=0$. If a group of automorphisms of $E$ preserves
$\Sigma_0$ and $\mathcal B$, then $\Sigma$ can be chosen stable under it.
Thus one support removes the ideal-class obstruction for every
prime-power exponent, including powers of two.
\end{lemma}
\begin{proof}
\noindent The localization formula
\[
 \Cl(R_\Sigma)=\Cl(\OO_E)/\langle[\mathfrak p]:\mathfrak p\in\Sigma\rangle
\]
\noindent is \cite[Proposition~7.4.4]{CohenAdvanced}; adjoining the supports
of finitely many class generators kills the quotient
\cite[Corollary~7.4.5]{CohenAdvanced}.
Approximation chooses these representatives prime to $\mathcal B$
\cite[Corollary~1.2.11]{CohenAdvanced}. Their automorphism orbits
still avoid the stable set $\mathcal B$ and give the required support.
\end{proof}

\begin{convention}
% Former title: Choice of support
\label{v207:support-convention}
\noindent In all field constructions, including $p=2$, choose a Galois-stable
support with $\Cl(R_\Sigma)=0$, avoiding primes above $p\Delta$
in the ordinary case. The lemma permits this auxiliary choice
without restricting the target. Abstract statements use their
stated weaker hypotheses. Class-group vanishing removes the
ideal-class obstruction, but the dyadic norm obstruction remains
(Section~\ref{v206:coverage}).
\end{convention}

\subsection{The ordinary norm conditions}
\label{v218:ordinary-notation}
\noindent For the ordinary construction let $n=p^a$, $p\nmid2\Delta$, and
$E=K(\zeta_n)$, with $\Sigma$ stable under $\Gal(E/\Q)$.
The quadratic field $K$ is disjoint from
$\Q(\zeta_n)$: a common nontrivial subfield would force the
conductor $\Delta$ of $K$ to divide $n$. We therefore have
\[
 \Gal(E/\Q)=\langle\tau\rangle\times\Gal(E/K),
 \qquad\Gal(E/K)\simeq(\Z_{n})^\times.
\]
\noindent The automorphisms are specified by
\[
\begin{gathered}
 \tau:E\xrightarrow{\sim}E,\qquad
       \tau|_K=\tau_K,\quad\tau(\zeta_n)=\zeta_n,\\
 \delta_u:E\xrightarrow{\sim}E,\qquad
       \delta_u|_K=\id_K,\quad\delta_u(\zeta_n)=\zeta_n^u
                       \quad(u\in\Z_n^\times).
\end{gathered}
\]
\noindent and complex conjugation is $j=\delta_{-1}$.
All these maps commute. Define
\begin{equation}\label{root170:U}
\begin{aligned}
 U_\Sigma&=\{x\in R_\Sigma^\times:x\tau(x)=xj(x)=1\},
                 &T_\Sigma&=U_\Sigma\cap\mu(E),\\
 \Lambda_\Sigma&=\{\lambda\in\Z[\Sigma]:\tau\lambda=j\lambda=-\lambda\},
                 &L_\Sigma&=\divi(U_\Sigma).
\end{aligned}
\end{equation}
\noindent Thus $L_\Sigma$ consists of divisors with generators satisfying
both norm equations. The ordinary \emph{Kummer module} is
\begin{equation}\label{root170:W}
 W_\Sigma=\{[x]\in\Sel_n(E,\Sigma):
                       \tau([x])=j([x])=[x]^{-1}\}.
\end{equation}
\noindent The notation suppresses $E,n$ and the involutions $\tau,j$.
For these fixed data, $W_\Sigma$ is independent of the character
$\chi$. The field and support chosen for an application may depend
on the target or on a family of targets.
We will prove $[\Lambda_\Sigma:L_\Sigma]<\infty$; equality need not hold.

\subsection{The general norm condition}
\noindent In the general construction let $n=p^a$ for any prime $p$, and
let $E/\Q$ be an abelian CM field containing $K$ and $\mu_n$:
$E$ is totally imaginary and $E^+=E^{\langle j\rangle}$ is
totally real. Supports are now $\Gal(E/\Q)$-stable.
The group $\Gal(E/K)$ and the map
$\delta\mapsto u_\delta$ are as in Lemma~\ref{v203:covariance}.
Complex conjugation $j\in\Gal(E/K)$ satisfies $u_j=-1$.
There need not be an automorphism $\tau$ fixing $\mu_n$ and
restricting to $\tau_K$. We impose only the condition for $j$:
\[
\begin{aligned}
 U^{(1)}_\Sigma&=\{x\in R_\Sigma^\times:xj(x)=1\},
                       &T^{(1)}_\Sigma&=\mu(E),\\
 \Lambda^{(1)}_\Sigma&=\{\lambda\in\Z[\Sigma]:j\lambda=-\lambda\},
                       &L^{(1)}_\Sigma&=\divi(U^{(1)}_\Sigma).
\end{aligned}
\]
\noindent Define the general \emph{Kummer module} by
\begin{equation}\label{v204:general-W}
 W^{(1)}_\Sigma=
        \{[x]\in E^\times/E^{\times n}:x\in U^{(1)}_\Sigma\}
                                      \subset\Sel_n(E,\Sigma).
\end{equation}
\noindent This module depends on $E,n,\Sigma,j$, again not on a particular
$\chi$ once these data are fixed. The condition $xj(x)=1$ is exact:
anti-invariance $j([x])=[x]^{-1}$ only makes $xj(x)$ an $n$th power and
is weaker at two.

\noindent These are $\Z_n$-modules, with the multiplicative notation fixed
in Section~\ref{root170:target}. When the module is free, a
\emph{Kummer basis} consists of classes $[z_j]$; the
$z_j\in E^\times$ are \emph{basis representatives}. In
\[
                  [\beta]=\prod_j[z_j]^{e_j},
\]
\noindent the exponents $e_j\in\Z_n$ are the \emph{Kummer coordinates}.
For $n=p$ the module is an $\F_p$-vector space.

\subsection{Divisors and Kummer classes}
\begin{lemma}
% Former title: Kernel and rank of the divisor map

\label{v203:divisor-kernel}
\noindent The divisor maps give exact sequences
\begin{equation}\label{root170:lattice-diagram}
\begin{aligned}
 1&\longrightarrow T_\Sigma\longrightarrow U_\Sigma
           \xrightarrow{\divi}L_\Sigma\longrightarrow0
                         &&\text{(ordinary)},\\
 1&\longrightarrow\mu(E)\longrightarrow U^{(1)}_\Sigma
           \xrightarrow{\divi}L^{(1)}_\Sigma\longrightarrow0
                         &&\text{(general)}.
\end{aligned}
\end{equation}
\noindent The images $L_\Sigma$ and $L^{(1)}_\Sigma$ have finite index in
$\Lambda_\Sigma$ and $\Lambda^{(1)}_\Sigma$, respectively.
In the ordinary case $p\nmid|T_\Sigma|$.
\end{lemma}
\begin{proof}
\noindent For a norm-one unit with zero divisor, every complex absolute value
is one. Kronecker's theorem identifies the kernels with the stated
roots of unity \cite[Chapter~I, Proposition~7.1]{NeukirchANT}.
Put $h_E=|\Cl(\OO_E)|$ and choose $x\in E^\times$ with
$\divi(x)=h_E\lambda$. Then
\[
 \begin{array}{c|c|c}
 \lambda&\text{norm-one element}&\text{divisor}\\ \hline
 \lambda\in\Lambda_\Sigma^{(1)}&x/j(x)&2h_E\lambda\\[2pt]
 \lambda\in\Lambda_\Sigma&x\tau(j(x))/(\tau(x)j(x))&4h_E\lambda
 \end{array}
\]
\noindent so both images have finite index. In the ordinary case the
$p$-power roots of unity of $E$ are exactly $\mu_n$, since the
next cyclotomic degree $p$ cannot divide $[E:\Q(\zeta_n)]=2$.
As $\tau$ fixes $\mu_n$, an element of $T_\Sigma\cap\mu_n$
satisfies $\xi^2=1$; thus $p\nmid|T_\Sigma|$.
\end{proof}

\begin{proposition}
% Former title: From norm-one units to Kummer classes

\label{root170:U-W}\label{v203:unitary-Selmer}
\noindent Let $\Sigma$ be finite and $\Gal(E/\Q)$-stable, with $E$ and
$n=p^a$ as in the respective construction. The class map
$x\mapsto[x]$ induces isomorphisms
\[
\begin{aligned}
 U_\Sigma/U_\Sigma^n&\xrightarrow{\sim}W_\Sigma
           &&\text{if }\Cl(R_\Sigma)[p]=0,\\
 U^{(1)}_\Sigma/(U^{(1)}_\Sigma)^n
                  &\xrightarrow{\sim}W^{(1)}_\Sigma
           &&\text{without a class-group hypothesis}.
\end{aligned}
\]
\noindent In the general case, if $n$ is odd and $\Cl(R_\Sigma)[p]=0$, then
\[
 W^{(1)}_\Sigma=
 \{[x]\in\Sel_n(E,\Sigma):j([x])=[x]^{-1}\}.
\]
\end{proposition}
\noindent The hypothesis $\Cl(R_\Sigma)[p]=0$ supplies a unit representative
for every Selmer class in the ordinary module. The general module
is defined as the image of the exact norm-one units, so its quotient
isomorphism needs no such hypothesis. The stronger convention
$\Cl(R_\Sigma)=0$ is used only in the field constructions.
\begin{proof}
\noindent \emph{General case.} Surjectivity is the definition of
$W^{(1)}_\Sigma$. If $y=c^n\in U^{(1)}_\Sigma$ with $c\in E^\times$,
valuations put $c$ in $R_\Sigma^\times$. The element $cj(c)$ is
totally positive and its $n$th power is one, so $cj(c)=1$.
Thus the kernel is $(U^{(1)}_\Sigma)^n$.
For the equality when $n$ is odd, the Selmer exact sequence
provides a unit representative $x$ of any supported class with
$j([x])=[x]^{-1}$. Choose $b$ with $2b\equiv1\pmod n$; then
$(x/j(x))^b\in U^{(1)}_\Sigma$ and its class is $[x]^{2b}=[x]$.

\noindent \emph{Ordinary case.} The same exact sequence provides a unit
representative $x$ of any class in $W_\Sigma$. Then
\[
 \Pi(x)=\frac{x\tau(j(x))}{\tau(x)j(x)}\in U_\Sigma,
                    \qquad[\Pi(x)]=[x]^4.
\]
\noindent An exponent inverse to four modulo the odd number $n$ gives a
norm-one representative of $[x]$. If $y=c^n\in U_\Sigma$,
put $\xi=c\tau(c)\in\mu_n$ and $c_0=c\xi^{(n-1)/2}$.
Since $\tau$ fixes $\mu_n$, we have $c_0^n=y$ and
$c_0\tau(c_0)=1$. Also $c_0j(c_0)$ is totally positive and
$(c_0j(c_0))^n=1$, so $c_0j(c_0)=1$. Valuations put $c_0$
in $R_\Sigma^\times$. Thus the kernel is $U_\Sigma^n$.
\end{proof}

\begin{corollary}
% Former title: Ranks of the Kummer modules

\label{v215:Kummer-ranks}
\noindent In the ordinary case assume $\Cl(R_\Sigma)[p]=0$;
no class-group hypothesis is needed for the general quotient. Then
\[
 W_\Sigma\simeq(\Z_n)^r,\quad r=\rank\Lambda_\Sigma,
 \qquad
 W^{(1)}_\Sigma\simeq(\Z_n)^{r_1+1},\quad
                         r_1=\rank\Lambda^{(1)}_\Sigma.
\]
\end{corollary}
\begin{proof}
\noindent Use the exact sequences of Lemma~\ref{v203:divisor-kernel} and
Proposition~\ref{root170:U-W}. The free parts contribute $r$ and
$r_1$ copies of $\Z_n$. Ordinary torsion contributes nothing;
$\mu(E)/\mu(E)^n\simeq\Z_n$ contributes the extra general factor.
\end{proof}

\section{The target classes in the Kummer modules}
\label{v206:coverage}\label{exc170:overview}\label{exc170:dyadic}

\noindent To use a Kummer basis, the target class must belong to its module.
The ordinary case follows from its two anti-invariances. In the
general case we first obtain an exact norm-one representative,
then remove its valuations outside $\Sigma$. For powers of two,
the first step requires signature correction and a quadratic norm
equation; the coefficient enlargement below solves that equation.

\subsection{Ordinary anti-invariant characters}
\noindent Use the ordinary notation of Section~\ref{v218:ordinary-notation}
and the support convention~\ref{v207:support-convention}.

\begin{theorem}
% Former title: The ordinary target belongs to the Kummer module

\label{v203:ordinary-support}
\noindent Let $\chi:G_{\mathfrak m}^+\to\Z_n$ be anti-invariant, with
$\tau_K\mathfrak m=\mathfrak m$. If $\Sigma$ contains the primes
above $\mathfrak m$ away from $p\Delta$, then
\[
                    \nu^{-1}(\chi|_{G_E})\in W_\Sigma.
\]
\noindent Equivalently, there is $\beta_\chi\in U_\Sigma$ with
$\nu_{\beta_\chi}=\chi|_{G_E}$.
\end{theorem}
\begin{proof}
\noindent Choose $\rho\in E^\times$ with $\nu_\rho=\chi|_{G_E}$ by
\eqref{root170:kummer-iso}. Kummer covariance
(Lemma~\ref{v203:covariance}) gives $j([\rho])=[\rho]^{-1}$.
Since $\tau$ fixes $\mu_n$, Galois equivariance of reciprocity
identifies its action on the restricted character with
$[\mathfrak a]\mapsto[\tau_K\mathfrak a]$; hence
$\tau([\rho])=[\rho]^{-1}$ as well.

\noindent We verify the supported valuations. At primes above $p$, the
involution $j$ fixes each prime: the cyclotomic extension is
totally ramified at $p$, whereas $K/\Q$ is unramified there.
At a prime above a rational divisor of $\Delta$, the involution
$\tau$ fixes that prime: the quadratic extension is ramified
and the cyclotomic extension is unramified there. In either case,
the corresponding anti-invariance gives
$v_{\mathfrak p}(\rho)\equiv-v_{\mathfrak p}(\rho)\pmod n$.
Since $n$ is odd, $v_{\mathfrak p}(\rho)\in n\Z$. At every other
prime outside $\mathfrak m$, the extension
$E(\sqrt[n]{\rho})/E$ is unramified. Its normalized valuations
therefore restrict to those of $E$, giving
$v_{\mathfrak p}(\rho)\in n\Z$ there as well.
Thus $[\rho]\in W_\Sigma$.

\noindent The support convention and Proposition~\ref{root170:U-W} now give
a representative $\beta_\chi\in U_\Sigma$.
\end{proof}

\noindent For an arbitrary ordinary character, apply the theorem to $\chi^-$
and recover $\chi$ by Lemma~\ref{root170:assembly}.
The prime $p$ may divide the strict class number.

\subsection{General coefficients and support}
\label{exc170:coefficientfields}
\noindent Now let $n=p^a$ and $\chi:G_{\mathfrak f}^+\to \Z_n$ be
arbitrary. Choose the coefficients and support as follows.

\noindent For odd $p$, put
\[
                    \psi=\chi,\qquad
       \mathfrak g=\mathfrak f,\qquad E=K(\zeta_n).
\]
\noindent For $p=2$, choose $h=\sqrt d$ positive at $\sigma_1$ and
negative at $\sigma_2$. Let
\[
 \theta_{-1},\theta_h:G_K\longrightarrow \Z_2
\]
\noindent be the quadratic characters of $K(i)$ and $K(\sqrt h)$, and
let $\mathfrak d_{-1},\mathfrak d_h$ be their relative
discriminant ideals. Let $c_i\in G_K$ be complex conjugation
at the real place $\sigma_i$, and write
$\chi(c_i)=(n/2)a_i$, $a_i\in\{0,1\}$. Define
\begin{equation}\label{exc170:signcorrection}
 \theta=a_1\theta_{-1}+(a_1+a_2)\theta_h,
               \qquad\psi=\chi-\frac n2\theta,
\end{equation}
\noindent where the coefficients of $\theta$ are taken modulo two and
$\Z_2\to \Z_n$ sends $1$ to $n/2$. Let
$a_0\in\{1,-1,h,-h\}$ have sign $(-1)^{a_i}$ at $\sigma_i$;
then $\theta$ is the character of $K(\sqrt{a_0})/K$. Use
\begin{equation}\label{exc170:enlargedmodulus}
 \mathfrak g=\operatorname{lcm}
                    (\mathfrak f,\mathfrak d_{-1},\mathfrak d_h)
\end{equation}
\noindent and the following sufficient level:
\begin{equation}\label{exc170:dyadlevel}
 \begin{gathered}
 b=\max\bigl(\{4,a+2\}\cup
       \{v_2(\ell+1)+1:\ell\mid\Nm(\mathfrak g),
                                      \ell\equiv3\pmod4\}\bigr),\\
                         E=K(\zeta_{2^b}).
 \end{gathered}
\end{equation}
\noindent The terms $a+2$ and $4$ handle primes above two; the remaining
terms handle odd conductor primes, as proved below.

\noindent In either case $E/\Q$ is a known abelian CM field. Write
$N_{\mathfrak g}=|\OO_K/\mathfrak g|$ for the positive integer
norm of $\mathfrak g$. Choose the support in two parts:
\begin{equation}\label{exc170:support}
 \begin{gathered}
 \Sigma_{\rm ram}=\{{\mathfrak p}\in\mathcal P_E:
            {\mathfrak p}\cap\Z=(\ell),\ \ell\mid pN_{\mathfrak g}\},\\
                  \Sigma=\Sigma_{\rm ram}\cup\Sigma_{\rm cl}.
 \end{gathered}
\end{equation}
\noindent Here $\ell\mid pN_{\mathfrak g}$ is divisibility of positive
integers; taking all primes above $\ell$ ensures Galois stability.
Choose ideals whose classes generate $\Cl(\OO_E)$, and let
$\Sigma_{\rm cl}$ contain their prime factors and all their Galois
conjugates. Then
\begin{equation}\label{exc170:classkill}
                         \Cl(R_\Sigma)=0,
\end{equation}
\noindent as in Lemma~\ref{v203:class-support}.

\noindent An explicit alternative, requiring no class-group generators, is
\begin{equation}\label{v211:minkowski-support}
 \begin{gathered}
 \Sigma_{\rm cl}=\{{\mathfrak p}\in\mathcal P_E:\Nm({\mathfrak p})\le B_E\},\qquad [E:\Q]=2r_E,\\
 B_E=\left(\frac4\pi\right)^{r_E}
       \frac{(2r_E)!}{(2r_E)^{2r_E}}
       \sqrt{|\operatorname{disc}(E)|}.
 \end{gathered}
\end{equation}
\noindent Minkowski's theorem gives an integral ideal of norm at most $B_E$
in each class, whose prime factors lie in this set
\cite[Chapter I, Section 6, Exercise 3]{NeukirchANT}.
This larger support is used in the Hilbert-field application.

\begin{theorem}
% Former title: Norm-one representatives of the general target

\label{exc170:allcoverage}
\noindent For the preceding coefficient field $E$ and character $\psi$,
let $\Sigma$ be any finite Galois-stable support containing
$\Sigma_{\rm ram}$ with $\Cl(R_\Sigma)=0$, as in
\eqref{exc170:support}. Then
\[
                    \nu^{-1}(\psi|_{G_E})\in W^{(1)}_\Sigma.
\]
\noindent Equivalently, there exists $\beta_\psi\in R_\Sigma^\times$ with
\[
                  \beta_\psi j(\beta_\psi)=1,
                 \qquad\nu_{\beta_\psi}=\psi|_{G_E}.
\]
\noindent For $p=2$, $\psi$ is totally even with conductor dividing
$\mathfrak g$. There is no restriction on $p$ relative to
$\Delta$ or any class number.
\end{theorem}
\noindent The next three subsections prove the ingredients: exact norm one
at odd order, the dyadic norm equation after signature correction,
and support normalization. We combine them in the proof at the
end of the section.

\subsection{Odd order cases}
\begin{lemma}
% Former title: Norm one at odd order
\label{v203:odd-normone}
\noindent If $n$ is odd and $j([\rho])=[\rho]^{-1}$ in
$E^\times/E^{\times n}$, then $[\rho]$ has a norm-one
representative. In particular this applies to
$\nu^{-1}(\chi|_{G_E})$ for every character
$\chi:G_K\to\Z_n$, even when $p\mid\Delta$.
\end{lemma}
\begin{proof}
\noindent Choose $k\in\Z$ with $2k\equiv1\pmod n$. Then
\[
 \beta=(\rho/j(\rho))^k,
       \qquad\beta j(\beta)=1,
       \qquad[\beta]=[\rho]^{2k}=[\rho].
\]
\noindent For a restricted character, Lemma~\ref{v203:covariance} gives
the required anti-invariance under $j$.
\end{proof}

\begin{lemma}
% Former title: Removal of real signatures
\label{v203:signatures}
\noindent The character $\psi$ in \eqref{exc170:signcorrection} vanishes
at both real inertia groups, and its finite conductor divides
$\mathfrak g$.
\end{lemma}
\begin{proof}
\noindent At the two real places,
\[
 (\theta_{-1}(c_1),\theta_{-1}(c_2))=(1,1),\qquad
 (\theta_h(c_1),\theta_h(c_2))=(0,1),
\]
\noindent so \eqref{exc170:signcorrection} gives $\psi(c_1)=\psi(c_2)=0$.
The conductor assertion follows from the quadratic
conductor--discriminant formula and the conductor bound for a sum
\cite[Chapter VI, Section 3]{SerreLocal}.
\end{proof}
\noindent Thus $H_\psi$ is totally real, which supplies the positivity
needed for the norm equation.

\subsection{The dyadic cases}
\noindent Write $E_0=K(\zeta_{2n})=k_0(i)$, where $k_0=E_0^+$.
The obstruction below is the solvability of
\[
              N_{E_0/k_0}(b)=A,
              \qquad A\in k_0^{\times,+}.
\]
\noindent By the Hasse norm theorem it suffices to solve this equation at
every completion \cite[Part III, Corollary 4.8]{NeukirchBonn}.
We use the following local criteria: an unramified quadratic
extension has norm group consisting of the elements of even
valuation, and its quadratic norm obstruction is killed by an
even local extension degree. The latter follows from
\[
                  (-1,A)_w=(-1,A)_v^{[k_w:(k_0)_v]},
\]
\noindent where $(-1,A)_v=1$ means that $A$ is a norm from $k_{0,v}(i)$
\cite[Chapter XIII, Section 3, Proposition 7]{SerreLocal}.

\begin{lemma}
% Former title: The norm obstruction
\label{v203:norm-obstruction}
\noindent Let $H=H_\psi$ for a totally even character
$\psi:G_K\to\Z_n$, where $n=2^a$, with finite conductor
dividing $\mathfrak g$. Put $E_0=K(\zeta_{2n})$,
$k_0=E_0^+$, and $L=HE_0$.
There is $\alpha\in L^\times$ such that
\[
 \alpha^n=\rho\in E_0^\times,\qquad
                         \nu_\rho=\psi|_{G_{E_0}}.
\]
\noindent If $J$ is complex conjugation on $L$, then
\[
                         A=\alpha J(\alpha)\in k_0^{\times,+}.
\]
\noindent The class $[\rho]$ has a norm-one representative over $E_0$
if and only if
\begin{equation}\label{exc170:normequation}
                b_0j(b_0)=A\qquad(b_0\in E_0^\times).
\end{equation}
\noindent This equation is locally solvable at every real place and at
every finite place not lying above $2\mathfrak g$.
\end{lemma}
\begin{proof}
\noindent Since $H$ and $k_0$ are totally real, $L=Hk_0(i)$ is CM.
Its conjugation $J$ fixes $H$ and restricts to $j$ on $E_0$.
Kummer theory supplies $\alpha$ in $L$, even if the restriction has
order smaller than $n$ \cite[Theorem 5.30]{MilneFT}.
For $g\in\Gal(L/E_0)$, write $g(\alpha)=\zeta_n^{\psi(g)}\alpha$.
The maps $g$ and $J$ commute, and $J$ inverts $\zeta_n$.
Hence $g(\alpha J(\alpha))=\alpha J(\alpha)$, so $A\in E_0$. It is fixed by $J$
and positive at every real embedding of $k_0$.

\noindent Every representative of $[\rho]$ is $(\alpha/b_0)^n$ for some
$b_0\in E_0^\times$. Its norm is
$(A/(b_0j(b_0)))^n$. The ratio is totally positive, so that
norm equals one precisely when \eqref{exc170:normequation} holds.
Here $E_0=k_0(i)$, so it is the quadratic norm equation just
described. Positivity gives a solution at real places, and a
split local extension imposes no condition.

\noindent At an odd nonsplit place outside $\mathfrak g$, both
$E_0/k_0$ and $L/E_0$ are unramified. The valuations of $\alpha$
at all primes of $L$ above a fixed prime of $E_0$ are equal,
since its Galois conjugates differ by roots of unity. That
prime of $E_0$ is fixed by $j$, so $\alpha$ and $J(\alpha)$ have the
same valuation. The ramification indices are one; consequently
$v(A)=2v(\alpha)$ is even in $k_0$. The unramified quadratic norm
criterion proves local solvability there.
\end{proof}

\begin{proposition}
% Former title: A sufficient dyadic coefficient level

\label{exc170:dyadiccoverage}
\noindent For $E$ in \eqref{exc170:dyadlevel}, every totally even character
of exponent dividing $2^a$ and finite conductor dividing
$\mathfrak g$ has a norm-one representative after restriction to $E$.
\end{proposition}
\begin{proof}
\noindent Put $k=E^+$ and use $A\in k_0^{\times,+}$ from
Lemma~\ref{v203:norm-obstruction}. We check the remaining local
norm conditions and then apply the Hasse norm theorem.

\noindent \emph{Odd primes.} At a rational prime $\ell\equiv1\pmod4$, $i\in\Q_\ell$,
so the extension is locally split. For
$\ell\equiv3\pmod4$, put $r=v_2(\ell+1)$. Its residue degree
in $\Q(\zeta_{2^b})^+$ is
\[
 \operatorname{ord}_{\Z_{2^b}^\times/\{\pm1\}}(\ell)
                                      =2^{\max(0,b-r)}.
\]
\noindent Indeed odd powers of $\ell$ are congruent to $-1$ modulo
exactly the same power $2^r$, whereas for even $m$,
$v_2(\ell^m-1)=r+v_2(m)$. If $\ell\mid\Nm(\mathfrak g)$,
our choice $b\ge r+1$ makes the residue degree even. The
completion of this real cyclotomic field is unramified at
$\ell$, and contains the unramified quadratic extension
$\Q_\ell(i)$. Thus $E/k$ is locally split at those primes,
also after adjoining $K$. At all other odd primes the local
norm equation already had a solution over $k_0$.

\noindent \emph{Primes above two.} Here $k/k_0$ has even local degree. Write
$F_t=\Q_2(\zeta_{2^t}+\zeta_{2^t}^{-1})$. For $a\ge2$,
$F_{a+2}/F_{a+1}$ is quadratic and $F_{a+2}/\Q_2$ is cyclic;
its unique quadratic subfield $\Q_2(\sqrt2)$ already lies in
$F_{a+1}$. A completion of the quadratic field $K$ cannot
consume this step: its intersection with $F_{a+2}$ is either
$\Q_2$ or that quadratic subfield, and in both cases is already
contained in $F_{a+1}$. For $a=1$, $k_0=K$ and
$[F_4:\Q_2]=4$, so adjoining $F_4$ to a completion of $K$
has even degree. The bound $b\ge\max(4,a+2)$ includes these
steps, and further enlargement preserves evenness. The local
degree formula for the quadratic Hilbert symbol therefore makes
$A$ a norm at the primes above two. Positivity handles infinity.

\noindent \emph{Global norm.} The cyclic Hasse norm theorem gives $b_1\in E^\times$
with $b_1j(b_1)=A$ \cite[Part III, Corollary 4.8]{NeukirchBonn}.
In $HE$, put $z=\alpha/b_1$. Then $zJ(z)=1$ and
$\beta=z^n=\rho/b_1^n\in E^\times$ represents the same
restricted character and satisfies $\beta j(\beta)=1$.
\end{proof}

\begin{example}
% Former title: Why total evenness and class number one do not suffice

\label{v204:dyadic-Selmer-example}\label{v205:dyadic-enlargement}
\noindent Take the totally real extension $H=K(\sqrt{11})$ of
$K=\Q(\sqrt3)$. Over $E_0=K(i)=\Q(\zeta_{12})$ its class
is $[11]$, with $j([11])=[11]=[11]^{-1}$ modulo squares.
The two primes over $11$ are $j$-fixed, since their Frobenius
is $11=-1\pmod{12}$. At each, every representative $11c^2$
has odd valuation, while a norm-one element has valuation zero.
Thus no norm-one representative exists. Nevertheless
$\Cl(\OO_{E_0})=0$: the discriminant is $144$ and Minkowski's
bound is $18/\pi^2<2$.

\noindent The product of the three quadratic discriminants gives
$\operatorname{disc}(H)=12\cdot44\cdot33=12^2\cdot11^2$.
Both primes over $11$ ramify, so the relative discriminant and
finite conductor are $(11)$. The sufficient level is $b=4$.
In $E=K(\zeta_{16})$, take
\[
 b_1=\sqrt3+2i\sqrt2,\qquad b_1j(b_1)=11,
 \qquad \beta=\frac{11}{b_1^2}
                  =\frac{-5-4\sqrt{-6}}{11}.
\]
\noindent Then $\beta j(\beta)=1$ and $[\beta]=[11]$. This solution
already lies in $K(\zeta_8)$, below the uniform sufficient level.
\end{example}

\subsection{Removing the remaining valuations}
\noindent We can now remove valuations outside the chosen support while
preserving norm one. Only the vanishing of the localized class
group is used here.

\begin{lemma}
% Former title: Support normalization
\label{v203:supported-normone}
\noindent Let $\Sigma$ be any finite $j$-stable support in the CM coefficient
field $E$ with $\Cl(R_\Sigma)=0$. If $\rho j(\rho)=1$ and
$[\rho]\in\Sel_n(E,\Sigma)$, then its class belongs to
$W^{(1)}_\Sigma$.
\end{lemma}
\begin{proof}
\noindent Since $\rho j(\rho)=1$, its divisor vanishes at $j$-fixed primes.
Choosing one prime from each other pair and using the Selmer
condition, write
\[
              \sum_{{\mathfrak p}\notin\Sigma}v_{\mathfrak p}(\rho)[{\mathfrak p}]=n(D_0-jD_0)
\]
\noindent for a divisor $D_0$ outside $\Sigma$. Triviality of $\Cl(R_\Sigma)$
gives $c\in E^\times$ with
$\sum_{{\mathfrak p}\notin\Sigma}v_{\mathfrak p}(c)[{\mathfrak p}]=D_0$. Then
\begin{equation}\label{exc170:supportnormalize}
 \beta=\rho/(c/j(c))^n\in R_\Sigma^\times,
              \qquad\beta j(\beta)=1,
              \qquad[\beta]=[\rho].
\end{equation}
\noindent This gives the required norm-one unit.
\end{proof}

\begin{proof}[Proof of Theorem~\ref{exc170:allcoverage}]
\noindent \emph{1. Exact norm one.}
For odd $p$, Lemma~\ref{v203:odd-normone} gives
$\rho\in E^\times$ with $\nu_\rho=\chi|_{G_E}$ and
$\rho j(\rho)=1$.
For $p=2$, Lemma~\ref{v203:signatures} makes $\psi$ totally even;
Proposition~\ref{exc170:dyadiccoverage} then gives such a $\rho$
over the specified field $E$.

\noindent \emph{2. The Selmer condition.}
The extension $E(\sqrt[n]{\rho})/E$ is unramified outside
the primes above $p\mathfrak g$. For $\mathfrak p\notin\Sigma$,
choose a normalized valuation $w$ above $\mathfrak p$ in this
extension. Since its ramification index is one,
\[
 v_{\mathfrak p}(\rho)=n\,w(\sqrt[n]{\rho})\in n\Z,
 \qquad [\rho]\in\Sel_n(E,\Sigma).
\]

\noindent \emph{3. A supported representative.}
Lemma~\ref{v203:supported-normone} uses $\Cl(R_\Sigma)=0$ to
replace $\rho$ by a norm-one $\Sigma$-unit in the same class.
Thus $[\rho]\in W_\Sigma^{(1)}$.
\end{proof}

\section{Ordinary Fermat Jacobians}
\label{geo170:curves}

\noindent Use the ordinary notation of Section~\ref{v218:ordinary-notation}:
\[
 n=p^a>1,\qquad p\nmid2\Delta,\qquad
 m=\Delta n,\quad F=\Q(\zeta_m),\quad E=K(\zeta_n).
\]
\noindent Choose $\zeta_m=e^{2\pi i/m}$, $\zeta_n=\zeta_m^\Delta$ and
$\zeta_\Delta=\zeta_m^n$. Write $\chi_\Delta$ for the primitive
quadratic character modulo $\Delta$; $\tau,\delta_u,j$ have
their earlier meanings. For an arbitrary ordinary character,
Lemma~\ref{root170:assembly} first separates the cyclotomic
part; we construct the Fermat input for its anti-invariant part.

\noindent The classical Fermat decomposition gives one-dimensional character
spaces and CM by $F$ on the factor primitive in both indices
\cite[Section~1]{KoblitzRohrlich}\cite[Section~4.2]{OY}.
Restriction to $K\subset\Q(\zeta_\Delta)$ groups these spaces
according to its two real embeddings. Their Frobenius determinant
quotient is the input $A_{\mathfrak q}(v)$: CM supplies the
eigenvalues, and RM determines their exponents $+1$ or $-1$.

\subsection{The curve, its endomorphisms, and its determinant line}
\noindent Our geometric input is the smooth projective model of the affine curve
\begin{equation}\label{v229:ordinary-curve}
                   C_{n,\Delta}:\quad Y^n+Z^\Delta=1.
\end{equation}
\noindent Write $C=C_{n,\Delta}$ and $J=\Jac(C)$. Its two coordinate
actions will give the cyclotomic and quadratic structures below.
The Fermat curve $\mathcal F_m:X_1^m+X_2^m=1$ provides the
cover used to compute its character spaces and Frobenius:
\begin{equation}\label{geo170:fermat-dominates}
 \mathcal F_m\longrightarrow C,\qquad
 (X_1,X_2)\longmapsto(X_1^\Delta,X_2^n).
\end{equation}
\noindent The curve $C$ has a unique point $\infty$ over infinity. The deck transformations
$S(Y,Z)=(\zeta_nY,Z)$ and $T(Y,Z)=(Y,\zeta_\Delta Z)$ act on $J_F$ by
$[P-\infty]\mapsto[S(P)-\infty]$ and similarly for $T$.
They commute, and $S^n=T^\Delta=1$. On the open subset $YZ\ne0$ they
are the deck transformations of the finite \'etale cover
\begin{equation}\label{geo170:cover}
 C^\circ\longrightarrow\mathbf P^1_F\setminus\{0,1,\infty\},
 \qquad x=Y^n,
 \qquad\text{with group }\mu_n\times\mu_\Delta.
\end{equation}

\noindent Fix an auxiliary prime $\ell_0$ and an embedding
$\iota:F\hookrightarrow\overline\Q_{\ell_0}$. The Tate module and rational Tate module are
\[
 T_{\ell_0}(J)=\varprojlim_r J[\ell_0^r](\overline F),\qquad
 V_{\ell_0}(J)=T_{\ell_0}(J)\otimes_{\widehat{\Z}_{\ell_0}}\Q_{\ell_0}.
\]
\noindent The transition maps are multiplication by $\ell_0$. We work on
$V=V_{\ell_0}(J)\otimes_{\Q_{\ell_0}}{\overline\Q_{\ell_0}}$.
The prime $\ell_0$ is used only to form the Tate module; it must
differ from the residue characteristic of a prime where Frobenius
is evaluated. The same convention applies in
Section~\ref{geo170:universal}.
The actions of $S,T$ on torsion points induce commuting
${\overline\Q_{\ell_0}}$-linear endomorphisms of $V$.
For $v\in\Z_{n}$ and $c\in\Z_{\Delta}$, define
\begin{equation}\label{geo170:lines}
 V_{v,c}=\{w\in V:Sw=\iota(\zeta_n^v)w,
                         \ Tw=\iota(\zeta_\Delta^c)w\}.
\end{equation}
\noindent The usual character projectors in $\End_{\overline\Q_{\ell_0}}(V)$ are
\begin{equation}\label{geo170:projectors}
 e_v^Y=\frac1n\sum_{b=0}^{n-1}\iota(\zeta_n^{-vb})S^b,
 \qquad
 e_c^Z=\frac1\Delta\sum_{b=0}^{\Delta-1}\iota(\zeta_\Delta^{-cb})T^b.
\end{equation}
\noindent Character orthogonality gives idempotents with
$e_v^Ye_c^ZV=V_{v,c}$; their scalar denominators lie in ${\overline\Q_{\ell_0}}$.
To select the primitive $\Delta$-characters, set
\begin{equation}\label{geo170:primitive}
 e_\Delta^{\rm prim}=\frac1\Delta\sum_{b=0}^{\Delta-1}r_\Delta(b)T^b,
 \quad r_\Delta(b)=\sum_{c\in(\Z_{\Delta})^\times}\zeta_\Delta^{-cb}\in\Z,
 \quad J^{\rm prim}=\operatorname{im}(\Delta e_\Delta^{\rm prim}:J\to J).
\end{equation}
\noindent Here $r_\Delta(b)=\sum_{d\mid(\Delta,b)}d\mu(\Delta/d)$ is the
Ramanujan sum, with $\mu$ the M\"obius function
\cite[Chapter~8]{ApostolANT}. Its integral coefficients make
$\Delta e_\Delta^{\rm prim}$ an endomorphism of $J$.
Below $r=\varphi(\Delta)/2$.
Define the Gaussian endomorphism
\begin{equation}\label{geo170:gaussian}
 \mathsf G_\Delta=\sum_{b\in(\Z_{\Delta})^\times}\chi_\Delta(b)T^b.
\end{equation}
\noindent For $v\in(\Z_{n})^\times$ put
\begin{equation}\label{geo170:packet}
 M_v=e_v^Ye_\Delta^{\rm prim}V,
 \qquad M_v^\pm=\ker(\mathsf G_\Delta\mp\iota(\sqrt \Delta)\mid M_v).
\end{equation}
\noindent The Gaussian sum realizes $\sqrt \Delta\in K$ on these spaces:
$M_v^+$ collects the lines with $\chi_\Delta(c)=1$, and $M_v^-$
those with $\chi_\Delta(c)=-1$. Both have dimension $r$, as proved
below. Their determinant quotient is the line
\begin{equation}\label{geo170:relative-det}
 \mathscr L_v=\bigwedge^rM_v^+
             \otimes\left(\bigwedge^rM_v^-\right)^\vee.
\end{equation}

\noindent Let $\mathfrak q$ be a prime ideal of $F$ over a rational prime $\ell\nmid \Delta n$,
put $k_{\mathfrak q}=\mathcal O_F/\mathfrak q$ and $q=|k_{\mathfrak q}|$, and now require
$\ell_0\ne\ell$. The curve has good reduction: normalization of
the tame cover with disjoint branch sections $0,1,\infty$ is
smooth. The residue field contains the primitive root
$\bar\zeta_m=\zeta_m\bmod\mathfrak q$, so $m\mid q-1$.
Choose a generator $g$ of $k_{\mathfrak q}^\times$ and
$s_{\mathfrak q}\in\Z_m^\times$ satisfying
\begin{equation}\label{v211:finite-character}
 g^{(q-1)/m}=\bar\zeta_m^{\,s_{\mathfrak q}}.
\end{equation}
\noindent The power-residue characters are
\[
\begin{aligned}
 \theta_{\mathfrak q}:k_{\mathfrak q}^\times&\longrightarrow\mu_n(F),
       &g^t&\longmapsto\zeta_n^{s_{\mathfrak q}t},\\
 \xi_{\mathfrak q}:k_{\mathfrak q}^\times&\longrightarrow\mu_\Delta(F),
       &g^t&\longmapsto\zeta_\Delta^{s_{\mathfrak q}t},
\end{aligned}
 \qquad t\in\Z_{q-1}.
\]
\noindent For multiplicative characters
$\lambda_1,\lambda_2:k_{\mathfrak q}^\times\to\overline\Q^\times$,
extend each by zero at zero and define the Jacobi sum as follows
\cite[Section~2.1]{BerndtEvansWilliams}. In particular, for
$v\in(\Z_n)^\times$ and $c\in(\Z_\Delta)^\times$,
\begin{equation}\label{geo170:jacobi}
 \begin{gathered}
 J(\lambda_1,\lambda_2)=\sum_{x\in k_{\mathfrak q}}\lambda_1(x)\lambda_2(1-x),\\
 J_{\mathfrak q,c}(v)=J(\theta_{\mathfrak q}^{v},\xi_{\mathfrak q}^{c})\in F.
 \end{gathered}
\end{equation}
\noindent Finally, $\Frob_{\mathfrak q}$ denotes the $q$-power Frobenius endomorphism of the
reduction of $J$. Our Frobenius convention is
\begin{equation}\label{geo170:frob-convention}
 \Frob_{\mathfrak q}\text{ on the Tate module acts as arithmetic
 Frobenius }x\longmapsto x^q.
\end{equation}
\noindent This convention determines the sign and character indices in the
following theorem. Define the RM determinants $B_{\mathfrak q}^{\pm}(v)\in F$ by
\[
 \iota(B_{\mathfrak q}^{\pm}(v))
       =\det_{\overline\Q_{\ell_0}}(\Frob_{\mathfrak q}\mid M_v^\pm).
\]
\noindent The Jacobi formulas below identify these determinants as algebraic
numbers in $F$, independent of $\ell_0$. Thus $\ell_0$ is not a
parameter of $B_{\mathfrak q}^\pm(v)$ or $A_{\mathfrak q}(v)$.

\begin{theorem}
% Former title: The ordinary determinant coefficient

\label{geo170:detquotient}
\noindent The endomorphism $\mathsf G_\Delta$ gives real multiplication by $K$ on
$J^{\rm prim}$, and $\dim_{\overline\Q_{\ell_0}} M_v^+=\dim_{\overline\Q_{\ell_0}} M_v^-=r$.
The Frobenius scalar on $\mathscr L_v$ is $\iota(A_{\mathfrak q}(v))$, where
\begin{equation}\label{geo170:A}
 A_{\mathfrak q}(v)=\prod_{c\in(\Z_{\Delta})^\times}J_{\mathfrak q,c}(v)^{\chi_\Delta(c)}
       \in E^\times.
\end{equation}
\noindent The two determinants lie in $E$, and
\[
 A_{\mathfrak q}(v)=
             \frac{B_{\mathfrak q}^{+}(v)}{B_{\mathfrak q}^{-}(v)}.
\]
\noindent These coefficients satisfy
\begin{equation}\label{geo170:symmetries}
 \tau A_{\mathfrak q}(v)=A_{\mathfrak q}(v)^{-1},\qquad
 jA_{\mathfrak q}(v)=A_{\mathfrak q}(v)^{-1},\qquad
 \delta_u A_{\mathfrak q}(v)=A_{\mathfrak q}(uv).
\end{equation}
\noindent For every prime ideal $\mathfrak p$ of $E$ not over $\ell$,
$v_{\mathfrak p}(A_{\mathfrak q}(v))=0$.
\end{theorem}

\noindent The proof combines the classical Fermat decomposition with the
Frobenius--Jacobi correspondence
\cite[Theorem~A]{ColemanFermat}\cite[Theorem~1.2]{OY}.
The next three lemmata fix the multiplicities, RM action and sign
convention needed for the determinant quotient.

\begin{lemma}
% Former title: Deck-character multiplicities
\label{geo170:multiplicity}
\noindent The curve has genus
\begin{equation}\label{geo170:genus}
 g(C)=\frac{(n-1)(\Delta-1)}2.
\end{equation}
\noindent Moreover $\dim_{\overline\Q_{\ell_0}} V_{v,c}=1$ if $v,c$ are both nonzero, and it
is zero otherwise.
\end{lemma}
\begin{proof}
\noindent Riemann--Hurwitz for the degree-$n$ map to the $Z$-line, ramified
at the $\Delta$ roots of unity and infinity, gives the genus.
The Fermat character decomposition
\cite[Section~1]{KoblitzRohrlich}, or
\cite[Propositions~4.8 and 4.11]{OY}, assigns a one-dimensional
space to each pair with $v/n$, $c/\Delta$, and $v/n+c/\Delta$
nonintegral. Coprimality makes the last condition automatic when
$v,c\ne0$. The two trivial-character parts vanish, equivalently
\begin{equation}\label{geo170:normrelations}
 \sum_{b=0}^{n-1}S^b=0,\qquad
 \sum_{b=0}^{\Delta-1}T^b=0,
\end{equation}
\noindent because both complete deck quotients have genus zero.
\end{proof}

\subsection{The RM structure and its eigenspaces}
\label{geo170:rm}\label{geo170:intrinsic}
\noindent We write $\End_K^0(A)=\End_K(A)\otimes_\Z\Q$ for endomorphisms up to isogeny.
\begin{lemma}
% Former title: Real multiplication and its eigenspaces
\label{geo170:RM}
\noindent There is an injective homomorphism
\begin{equation}\label{geo170:RM-action}
 \mathsf r:K\hookrightarrow\End_K^0(J^{\rm prim}),
 \qquad\sqrt \Delta\longmapsto\mathsf G_\Delta|_{J^{\rm prim}}.
\end{equation}
\noindent For primitive $c$, its action on $V_{v,c}$ is through the embedding
of $K$ that sends $\sqrt \Delta$ to $\chi_\Delta(c)\iota(\sqrt \Delta)$.
\end{lemma}
\begin{proof}
\noindent The primitive quadratic Gauss sum for the chosen complex roots is
\begin{equation}\label{geo170:gauss-identity}
 \sum_{b\in(\Z_{\Delta})^\times}\chi_\Delta(b)\zeta_\Delta^b=\sqrt \Delta.
\end{equation}
\noindent This is the quadratic Gauss-sum evaluation
\cite[Corollary 4.6]{Washington}; in particular its square
is $\Delta$, with the positive sign for $\Delta>0$ in the standard embedding.
Replacing $b$ by $bc^{-1}$ gives
\begin{equation}\label{geo170:gauss-colour}
 \mathsf G_\Delta|_{V_{v,c}}=\chi_\Delta(c)\iota(\sqrt \Delta),
 \qquad c\in(\Z_{\Delta})^\times.
\end{equation}
\noindent The idempotent $e_\Delta^{\rm prim}$ is defined over $\Q$, since every
Galois automorphism permutes its primitive characters. On its image,
$T$ satisfies $\Phi_\Delta(T)=0$, and the action of $\Q(\zeta_\Delta)$ is
faithful by Lemma~\ref{geo170:multiplicity}. Equation
\eqref{geo170:gauss-colour} therefore gives
$\mathsf G_\Delta^2=\Delta$ and an embedding of $K$.
An automorphism sending $\zeta_\Delta$ to $\zeta_\Delta^b$ sends
$\mathsf G_\Delta$ to $\chi_\Delta(b)\mathsf G_\Delta$. Hence this endomorphism is
fixed by $\Gal(\overline\Q/K)$, proving \eqref{geo170:RM-action}.
There are $\varphi(\Delta)/2$ primitive indices of each sign, so
\[
 M_v^\pm=\bigoplus_{\chi_\Delta(c)=\pm1}V_{v,c},\qquad
 \dim M_v^\pm=r.
\]
\noindent Their coefficient projectors are
\begin{equation}\label{geo170:RM-projectors}
 e_\pm=\frac12\left(1\pm\frac{\mathsf G_\Delta}{\iota(\sqrt \Delta)}\right)
 \quad\text{on }M_v.
\end{equation}
\end{proof}
\noindent These are coefficient projectors for the two embeddings of
$K\otimes{\overline\Q_{\ell_0}}$, not an asserted isogeny decomposition of $J$.
For composite $\Delta$, the primitive projector remains essential:
\eqref{geo170:normrelations} removes only the trivial character.

\subsection{Evaluation of Frobenius}
\label{geo170:determinant}
\begin{lemma}
% Former title: The Jacobi eigenvalue
\label{geo203:Jacobi-Frobenius}
\noindent On the one-dimensional space $V_{v,c}$,
\begin{equation}\label{geo170:line-eigenvalue}
 \Frob_{\mathfrak q}|_{V_{v,c}}=-\iota(J_{\mathfrak q,c}(v))\id.
\end{equation}
\noindent Moreover
\begin{equation}\label{geo170:weight}
 J_{\mathfrak q,c}(v)jJ_{\mathfrak q,c}(v)=q,
 \qquad |\sigma J_{\mathfrak q,c}(v)|=\sqrt q
 \quad(\sigma:F\hookrightarrow\C).
\end{equation}
\end{lemma}
\begin{proof}
\noindent This is the classical Frobenius--Jacobi correspondence
\cite[Theorem~A]{ColemanFermat}, in the character-projector
form of \cite[Theorem~1.2 and Section~6]{OY}.
On the fibre above $x$, Frobenius has deck character
$\theta_{\mathfrak q}(x)^v\xi_{\mathfrak q}(1-x)^c$. Its degree-one trace is the
negative of the sum; thus our unsigned Jacobi sum gives
$-\iota(J_{\mathfrak q,c}(v))$. This records the convention needed here.
The norm and absolute values follow from
$J(\alpha,\beta)\overline{J(\alpha,\beta)}=q$ when
$\alpha,\beta,\alpha\beta$ are nontrivial
\cite[Section~2.1]{BerndtEvansWilliams}.
\end{proof}

\begin{proof}[Proof of Theorem~\ref{geo170:detquotient}]
\noindent The lemmata give \eqref{geo170:A}: both determinants have sign
$(-1)^r$, which cancels in their quotient. The group $\Gal(F/E)$
preserves each sign of $\chi_\Delta$, whereas $\tau$
interchanges the two signs; this proves membership in $E$ and
$\tau A=A^{-1}$. Equation~\eqref{geo170:weight} gives $jA=A^{-1}$,
and coefficient conjugation gives $\delta_uA(v)=A(uv)$.
Integrality and $JjJ=q$ restrict every Jacobi divisor, hence the
quotient's divisor, to primes above $\ell$.
\end{proof}

\subsection{The determinant as a quotient of norms}
\noindent The automorphisms of $F/E$ fix $\zeta_n$ and send $\zeta_\Delta$
to $\zeta_\Delta^c$ for the indices $c$ with $\chi_\Delta(c)=1$.
The indices of each RM sign occur in pairs $c,-c$, so $r$ is even
and the Jacobi signs multiply to one. Thus
\begin{equation}\label{geo170:norm-ratio}
 \begin{gathered}
 B_{\mathfrak q}^{+}(v)=N_{F/E}J_{\mathfrak q,1}(v)
       =\prod_{\chi_\Delta(c)=1}J_{\mathfrak q,c}(v),\qquad
 B_{\mathfrak q}^{-}(v)=\tau(B_{\mathfrak q}^{+}(v)),\\
 A_{\mathfrak q}(v)=\frac{B_{\mathfrak q}^{+}(v)}{B_{\mathfrak q}^{-}(v)}.
 \end{gathered}
\end{equation}
\noindent Thus the determinant quotient is a quotient of coefficient norms.
For comparison with normalization by complex conjugation,
$B_{\mathfrak q}^{\pm}(v)j(B_{\mathfrak q}^{\pm}(v))=q^r$ gives
\begin{equation}\label{cmp202:det-bridge}
 \begin{gathered}
 \frac{B_{\mathfrak q}^{+}(v)j(B_{\mathfrak q}^{-}(v))}{j(B_{\mathfrak q}^{+}(v))B_{\mathfrak q}^{-}(v)}=A_{\mathfrak q}(v)^2,\\
 [A_{\mathfrak q}(v)]=
 \left[\frac{B_{\mathfrak q}^{+}(v)j(B_{\mathfrak q}^{-}(v))}{j(B_{\mathfrak q}^{+}(v))B_{\mathfrak q}^{-}(v)}\right]^{(n+1)/2}
 \quad\text{in }E^\times/E^{\times n}.
 \end{gathered}
\end{equation}
\noindent Here $(n+1)/2$ is the inverse of $2$ modulo the odd integer $n$.
This identity describes the extra quadratic anti-invariance of the
ordinary determinant. Proposition~\ref{v215:ordinary-general-comparison}
uses it to compare the ordinary and uniform RM constructions.

\section{Ordinary RM divisors}
\label{geo170:ordinary-spanning}
\noindent We compute the divisors of the ordinary determinant coefficients
and prove that they span the required anti-invariant part.

\noindent Keep $n=p^a$, with $p$ odd and $(p,\Delta)=1$. Let $\Sigma$ be a finite
$\Gal(E/\Q)$-stable set of prime ideals of $E$ not over $p\Delta$, and
let $\Sigma_F$ be all primes of $F$ above it. Recall the ordinary
divisor lattice of Section~\ref{root170:ordinary-space}:
\[
 \Lambda_\Sigma=
 \{Z\in\Z[\Sigma]:\tau Z=-Z,\ jZ=-Z\},
\]
\noindent where $\Z[\Sigma]=\bigoplus_{\mathfrak p\in\Sigma}\Z[\mathfrak p]$ and the map used
in this section is
\[
 \operatorname{div}_{E,\Sigma}:R_\Sigma^\times\longrightarrow\Z[\Sigma],
 \qquad x\longmapsto\sum_{\mathfrak p\in\Sigma}v_{\mathfrak p}(x)[\mathfrak p].
\]
\noindent We abbreviate it to $\operatorname{div}_E$ on supported elements.
In particular its restriction $U_\Sigma\to\Lambda_\Sigma$ is the
map defined in Section~\ref{root170:ordinary-space}.

\begin{theorem}
% Former title: Spanning by the determinant divisors

\label{geo170:ordinary-complete}
\noindent One has
\[
 \operatorname{span}_\Q
 \{\operatorname{div}_EA_{\mathfrak q}(v):\mathfrak q\in\Sigma_F,
                                  \ v\in(\Z_{n})^\times\}
       =\Lambda_\Sigma\otimes\Q.
\]
\noindent There is no restriction on the class number of $K$.
\end{theorem}

\noindent The proof has three steps. Stickelberger's formula gives the
divisors of the Jacobi values. Fourier decomposition reduces their
span to eigenvalues proportional to odd Dirichlet $L$-values at
zero. Their nonvanishing gives every required component, including
at nonsplit primes. This is the classical minus Stickelberger
argument \cite{Sinnott1978}\cite[Section~8]{OY}, applied to our
quadratic determinant. For the associated twisted Bernoulli determinant,
see \cite[Theorem~1]{Endo2008}, with
$(m,f,k,H,\psi)=(n,\Delta,1,\{1\},\chi_\Delta)$.

\noindent The norm and support identities are those of
Theorem~\ref{geo170:detquotient}. We state the classical valuation
formula at arbitrary conductor for use also in the general case.

\begin{lemma}
% Former title: Jacobi valuations
\label{geo203:jacobi-valuations}
\noindent Let $M>2$, $L=\Q(\zeta_M)$, and let $\mathfrak q$ be a prime of
$L$ over $\ell\nmid M$, and put $k_{\mathfrak q}=\OO_L/\mathfrak q$,
$q_0=|k_{\mathfrak q}|$. Choose $g,s_{\mathfrak q}$ as in
\eqref{v211:finite-character}, with $m=M$ and $q=q_0$.
The inverse character, in the coordinates of
\eqref{v211:finite-character}, is
\[
 \lambda:k_{\mathfrak q}^\times\longrightarrow\mu_M(L),\qquad
 \lambda(g^t)=\zeta_M^{-s_{\mathfrak q}t}.
\]
\noindent Character powers are formed on $k_{\mathfrak q}^\times$ and then
extended by zero. Write $\{x\}=x-\lfloor x\rfloor$ and assume
\[
             r\ne0,\qquad t\ne0,\qquad r+t\ne0
                             \quad\text{in }\Z_M.
\]
\noindent Put
\[
 J_{\mathfrak q}(r,t)=\sum_{x\in k_{\mathfrak q}}\lambda^r(x)\lambda^t(1-x),
 \qquad f_{r,t}(b)=
 \geoFrac{\frac{rb}{M}}+\geoFrac{\frac{tb}{M}}
                         -\geoFrac{\frac{(r+t)b}{M}}.
\]
\noindent For $\sigma_b(\zeta_M)=\zeta_M^b$, the factorization of its
principal ideal, or equivalently its divisor, is
\[
 \operatorname{div}_LJ_{\mathfrak q}(r,t)
       =\sum_{b\in(\Z_{M})^\times}f_{r,t}(b)[\sigma_b^{-1}\mathfrak q].
\]
\noindent If $h$ is the order of $\ell$ modulo $M$, the valuation at
$\sigma_b^{-1}\mathfrak q$ is consequently
$\sum_{i=0}^{h-1}f_{r,t}(b\ell^i)$.
\end{lemma}
\begin{proof}
\noindent Apply Stickelberger's Jacobi factorization
\cite[Theorem~11.2.3]{BerndtEvansWilliams}.
For positive residue exponent $d$, the convention in
\cite[Section~8, equation~(8.5)]{OY} is
$\sum_b\{-db/M\}\sigma_b^{-1}$; substitution of $-r,-t$
gives the displayed coefficients. Since the decomposition group
is $\langle\sigma_\ell\rangle$, the indices $b\ell^i$ represent
the same prime and their coefficients add. This includes
nonsplit primes.
\end{proof}

\subsection{The divisor operator and its Fourier eigenvalues}
\noindent Let $G=\Gal(F/\Q)$, with $\sigma_b(\zeta_m)=\zeta_m^b$.
Define
\begin{equation}\label{geo170:ordinary-carry}
 c_{\Delta,n}(b)=\geoFrac{\frac{-\Delta b}{m}}+
 \geoFrac{\frac{-nb}{m}}-\geoFrac{\frac{-(\Delta+n)b}{m}},
 \qquad
 \mathcal C_{\Delta,n}=\sum_{b\in(\Z_{m})^\times}
                         c_{\Delta,n}(b)\sigma_b^{-1}.
\end{equation}
\noindent The coefficients are zero or one. In Lemma~\ref{geo203:jacobi-valuations},
the characters of \eqref{geo170:jacobi} have exponents $-\Delta,-n$
relative to the inverse character of order $m$. Hence
\begin{equation}\label{geo170:ordinary-divisor}
 \operatorname{div}_FJ_{\mathfrak q,1}(1)=\mathcal C_{\Delta,n}[\mathfrak q].
\end{equation}
\noindent This explains both the minus signs in the fractional parts and the
inverse Galois indices: $b$ records the conjugate prime
$\sigma_b^{-1}\mathfrak q$, with repeated primes summed as in the lemma.

\noindent Let $\widetilde H=\Gal(F/E)$ and
$N_{\widetilde H}=\sum_{h\in\widetilde H}h$. Choose
$\upsilon\in G$ fixing $\zeta_n$ and acting nontrivially on $K$.
The norm expression \eqref{geo170:norm-ratio} gives
\begin{equation}\label{geo170:A-divisor}
 \operatorname{div}_FA_{\mathfrak q}(1)=
 (1-\upsilon)N_{\widetilde H}\mathcal C_{\Delta,n}[\mathfrak q].
\end{equation}
\noindent For a number field $E$, write
$\Div(E)=\bigoplus_{\mathfrak p}\Z[\mathfrak p]$, over its nonzero
prime ideals. For a finite extension $F/E$, the standard extension
of ideals, written additively on divisors, is the map we denote by
\[
 \begin{aligned}
 \operatorname{ext}_{F/E}:\Div(E)&\longrightarrow\Div(F),\\
 [\mathfrak p]&\longmapsto
       \sum_{\mathfrak q\mid\mathfrak p}
                   e(\mathfrak q/\mathfrak p)[\mathfrak q].
 \end{aligned}
\]
\noindent Here $[\mathfrak p]$ is a prime divisor, not an ideal class, and
$e(\mathfrak q/\mathfrak p)$ is the ramification index. This is the
additive form of the ideal factorization
$\mathfrak p\OO_F=\prod_{\mathfrak q\mid\mathfrak p}
\mathfrak q^{e(\mathfrak q/\mathfrak p)}$
\cite[Chapter I, Section 8]{NeukirchANT}. It satisfies
\[
 \operatorname{ext}_{F/E}\bigl(\divi_E(x)\bigr)=\divi_F(x)
                         \qquad(x\in E^\times),
\]
\noindent so it compares the divisor of $A_{\mathfrak q}(1)\in E^\times$
with its computation in $F$ above. Here all primes in $\Sigma$
are unramified in $F/E$, so
$\operatorname{ext}_{F/E}[\mathfrak p]=
\sum_{\mathfrak q\mid\mathfrak p}[\mathfrak q]$; after extending
scalars to $\Q$, this identifies $\Q[\Sigma]$ with
$\Q[\Sigma_F]^{\widetilde H}$.
On the subspace where $\upsilon$ and $j$ both act by $-1$, the
factor $(1-\upsilon)N_{\widetilde H}$ is multiplication by
$2|\widetilde H|=\varphi(\Delta)$.

\begin{lemma}
% Former title: Fourier calculation
\label{geo203:ordinary-Fourier}
\noindent Let $\psi$ be an odd character modulo $n$, and regard
$\vartheta=\chi_\Delta\psi$ as a character of $G$. The eigenvalue of
$\mathcal C_{\Delta,n}$ on the $\vartheta$-eigenspace is
\begin{equation}\label{geo170:ordinary-eigenvalue}
 \sum_b c_{\Delta,n}(b)\vartheta(b)^{-1}
       =\vartheta(-(\Delta+n))L(0,\vartheta^{-1}).
\end{equation}
\noindent The eigenvalue of $(1-\upsilon)N_{\widetilde H}\mathcal C_{\Delta,n}$
on this space is
\begin{equation}\label{geo170:ordinary-full-eigenvalue}
 \varphi(\Delta)\vartheta(-(\Delta+n))L(0,\vartheta^{-1}).
\end{equation}
\end{lemma}
\begin{proof}
\noindent The first two terms vanish by summing respectively over
$\chi_\Delta$ and the nontrivial odd character $\psi$.
For any unit $d\bmod m$, substitute $y=-db$ to obtain
\[
 \sum_{b\in(\Z_{m})^\times}\geoFrac{\frac{-db}{m}}\vartheta(b)^{-1}
 =\vartheta(-d)\sum_{y=1}^m\frac ym\vartheta(y)^{-1}
 =-\vartheta(-d)L(0,\vartheta^{-1}).
\]
\noindent The Bernoulli identity applies also to imprimitive nonprincipal
characters \cite[Theorem~12.20]{ApostolANT}.
Take $d=\Delta+n$, which is prime to $m$; the additional
factor is $(1-\upsilon)N_{\widetilde H}=\varphi(\Delta)$
on the required eigenspace.
\end{proof}

\begin{lemma}
% Former title: Classical nonvanishing at zero
\label{geo203:L0-nonvanishing}
\noindent If $\lambda$ is a primitive odd Dirichlet character, then
$L(0,\lambda)\ne0$. If $\psi$ is odd modulo $n=p^a$ and
$(p,\Delta)=1$, the same holds for $\lambda=\chi_\Delta\psi$ and its inverse,
whether or not they are primitive modulo $\Delta n$.
\end{lemma}
\begin{proof}
\noindent For a primitive odd character of conductor $f$, let
$\mathfrak g(\lambda)=\sum_{b\bmod f}\lambda(b)e^{2\pi ib/f}$.
The functional equation and Dirichlet's nonvanishing theorem give
\begin{equation}\label{geo170:L0-equation}
 L(0,\lambda)=\frac{\mathfrak g(\lambda)}{i\pi}
                         L(1,\overline\lambda)\ne0;
\end{equation}
\noindent see \cite[equations~(25.15.5) and (25.15.9)]{DLMFDirichlet}.
The primitive character inducing $\psi$ has conductor $p^b$,
$1\le b\le a$, so the primitive character inducing
$\chi_\Delta\psi$ has conductor $p^b\Delta$. Enlarging to
$p^a\Delta$ introduces no new prime divisor and therefore no
Euler factor \cite[equation~(25.15.4)]{DLMFDirichlet}.
The primitive nonvanishing consequently applies unchanged.
\end{proof}

\begin{proof}[Proof of Theorem~\ref{geo170:ordinary-complete}]
\noindent Decompose
$\Q[\Sigma_F]$ into its prime orbits. The orbit of $\mathfrak q$ is the
permutation representation $\Q[G/G_{\mathfrak q}]$, where $G_{\mathfrak q}$ is the
decomposition group. Since $G$ is abelian, its complex characters
are exactly those trivial on $G_{\mathfrak q}$, each with multiplicity one.
Requiring $\widetilde H$ to act trivially and $\upsilon,j$ to act
by $-1$ selects exactly the characters $\chi_\Delta\psi$ with $\psi$
odd that occur in this permutation representation.
Lemmata~\ref{geo203:ordinary-Fourier} and
\ref{geo203:L0-nonvanishing} show that the operator in
\eqref{geo170:A-divisor} is invertible on this subspace.
Its image on the full permutation representation lies in that
subspace by the two anti-invariance identities already proved.
Since its restriction there is invertible, its image is exactly the
required subspace, and is spanned by its values on the prime basis.
Those values are precisely
$\operatorname{div}_FA_{\mathfrak q}(1)$. Taking the direct sum over the prime
orbits and descending through $\operatorname{ext}_{F/E}$ proves
the claim. In particular the assertion holds at nonsplit rational
primes; their decomposition groups merely remove some characters
from the representation.
\end{proof}

\begin{corollary}
% Former title: Primes above the product $p\Delta$
\label{geo170:ordinary-bad-support}
\noindent For any finite $\Gal(E/\Q)$-stable set $\Sigma$, put
$\Sigma^\circ=\{\mathfrak p\in\Sigma:\mathfrak p\nmid p\Delta\}$.
Then $\Lambda_\Sigma=\Lambda_{\Sigma^\circ}$, so the determinant
values above $\Sigma^\circ$ span $\Lambda_\Sigma\otimes\Q$.
\end{corollary}
\begin{proof}
\noindent As in Theorem~\ref{v203:ordinary-support}, $j$ fixes every prime
above $p$ and $\tau$ fixes every prime above a divisor of $\Delta$.
The corresponding anti-invariance forces its coefficient to vanish.
Thus deleting these primes leaves $\Lambda_\Sigma$ unchanged.
\end{proof}

\noindent Theorem~\ref{geo170:ordinary-complete} gives rational spanning.
Section~\ref{root170:saturation} defines the generated subgroup and
obtains a basis modulo $n$th powers by internal $p$-saturation.

\section{A basis of the ordinary Kummer module}
\label{root170:saturation}

\noindent A subgroup of finite index need not supply every Kummer class:
$\langle u^p\rangle$ has trivial image in
$\langle u\rangle/\langle u\rangle^p$.
The standard remedy is \emph{$p$-saturation}
\cite[Section 4.3, Lemma 4.3, Algorithm 4.9 and Theorem 4.10]{BFHP2022}.
All roots extracted here already lie in $E$; the extension
$\alpha^n=\beta$ comes later.

\begin{proposition}
% Former title: Elementary saturation at a prime

\label{v212:elementary-saturation}\label{v203:prime-to-p-index}
\label{v203:saturation-loop}
\noindent Let $U$ be a finitely generated abelian group, written multiplicatively,
and let $\Gamma\subset U$ have finite index $I$. Fix a prime $p$.
Throughout, the notation $\Gamma^{\mathrm{sat}}$ means saturation
at this fixed $p$, in the specified ambient group $U$; put
\begin{equation}\label{root170:Vsat}
 \Gamma^{\mathrm{sat}}
   =\{x\in U:x^{p^b}\in\Gamma\text{ for some }b\ge0\}.
\end{equation}
\noindent Then
\[
 \Gamma^{\mathrm{sat}}/\Gamma=(U/\Gamma)[p^\infty],\qquad
 p\nmid[U:\Gamma^{\mathrm{sat}}],\qquad
 [\Gamma^{\mathrm{sat}}:\Gamma]=p^{v_p(I)}.
\]
\noindent For every $a\ge1$, inclusion induces
\[
 \Gamma^{\mathrm{sat}}/(\Gamma^{\mathrm{sat}})^{p^a}
                       \xrightarrow{\sim}U/U^{p^a}.
\]
\noindent Starting from $\Gamma$, repeatedly adjoining an element $x\in U$
with $x^p$ in the current subgroup but $x$ outside it reaches
$\Gamma^{\mathrm{sat}}$ after $v_p(I)$ such enlargements.
\end{proposition}
\begin{proof}
\noindent The first quotient is the $p$-primary subgroup of $U/\Gamma$.
Multiplication by $p^a$ is bijective on the remaining quotient,
giving the isomorphism modulo powers. Each enlargement has index
$p$; a nontrivial $p$-primary quotient supplies another element
of order $p$. This is the standard finite-index argument
\cite[Section 4.3, Lemma 4.3]{BFHP2022}.
\end{proof}

\paragraph{The ordinary application.}
Use the ordinary data of Sections~\ref{root170:ordinary-space}
and~\ref{geo170:ordinary-spanning}, with $\Cl(R_\Sigma)=0$.
Enumerate $\Sigma_F\times\Z_n^\times$ as $((\mathfrak q_i,v_i))_i$ and define
\begin{equation}\label{v220:ordinary-input}
 A_i=A_{\mathfrak q_i}(v_i)
     =\prod_{c\in\Z_\Delta^\times}
           J_{\mathfrak q_i,c}(v_i)^{\chi_\Delta(c)}\in U_\Sigma.
\end{equation}
\noindent These are the Jacobi determinant quotients of
Theorem~\ref{geo170:detquotient}, including their coefficient conjugates.
Set
\[
 \begin{aligned}
 \Gamma&=\langle T_\Sigma,A_i\ (\text{all }i)\rangle\subset U_\Sigma,\\
 L_\Gamma&=\divi_E(\Gamma)
     =\sum_i\Z\,\divi_E(A_i)\subseteq L_\Sigma\subseteq\Lambda_\Sigma,\\
 B&=[\Lambda_\Sigma:L_\Gamma],\qquad r=\rank\Lambda_\Sigma.
 \end{aligned}
\]
\noindent Theorem~\ref{geo170:ordinary-complete} makes $L_\Gamma$ a full-rank
lattice. Smith normal form of its valuation matrix computes the
finite index $B$. Since $\ker(\divi_E|_{U_\Sigma})=T_\Sigma\subset\Gamma$,
\begin{equation}\label{v212:ordinary-indices}
 U_\Sigma/\Gamma\simeq L_\Sigma/L_\Gamma,\qquad
 [U_\Sigma:\Gamma]\mid B.
\end{equation}
\begin{corollary}
% Former title: A Kummer basis
\label{root170:saturation-theorem}
\noindent The class map induces
\begin{equation}\label{root170:V-W}
 \Gamma^{\mathrm{sat}}/(\Gamma^{\mathrm{sat}})^n
       \xrightarrow{\sim}U_\Sigma/U_\Sigma^n
       \xrightarrow{\sim}W_\Sigma.
\end{equation}
\noindent If $z_1,\ldots,z_r$ lift a $\Z$-basis of
$\Gamma^{\mathrm{sat}}/T_\Sigma$, then
\[
 (\Z_n)^r\xrightarrow{\sim}W_\Sigma,
       \qquad(e_j)_j\longmapsto\left[\prod_jz_j^{e_j}\right].
\]
\end{corollary}
\begin{proof}
\noindent Apply Propositions~\ref{v212:elementary-saturation} and~\ref{root170:U-W}.
The torsion contributes nothing because $p\nmid|T_\Sigma|$.
\end{proof}

\begin{corollary}
% Former title: When saturation is unnecessary
\label{v224:no-saturation}
\[
 \Gamma^{\mathrm{sat}}=\Gamma
   \quad\Longleftrightarrow\quad p\nmid[U_\Sigma:\Gamma].
\]
\noindent In particular, $p\nmid B$ suffices. In this case a Kummer basis
can be represented by products of the $A_i$.
\end{corollary}
\begin{proof}
\noindent Use Proposition~\ref{v212:elementary-saturation} and
$[U_\Sigma:\Gamma]\mid B$.
\end{proof}
\noindent Sections~\ref{v221:ordinary-quintic} and~\ref{v224:ordinary-cubic}
give examples with $B=1$ and $B=2$, respectively.

\paragraph{Computing the basis representatives.}
To compute the $z_j$ in Corollary~\ref{root170:saturation-theorem},
start with $\Gamma'=\Gamma$ and test representatives
$\gamma\in\Gamma'/\Gamma'^p$ by factoring $X^p-\gamma$ over $E$.
If $b\in E$ is a root, put
\[
 \xi=b\tau(b)\in\mu_p,\qquad c=b\xi^{(p-1)/2}.
\]
\noindent Then $c^p=\gamma$, $c\tau(c)=1$ and $cj(c)=1$: the first norm
identity follows because $\tau$ fixes $\mu_p$, and the second
because $cj(c)$ is positive and has $p$th power one.
Valuations put $c$ in $U_\Sigma$. It is the unique root there,
since $U_\Sigma\cap\mu_p=1$.

\noindent Adjoin $c$ if its divisor is outside $L_{\Gamma'}$.
Membership is tested on divisors because $T_\Sigma\subset\Gamma'$.
The tests are independent of the chosen representatives of
$\Gamma'/\Gamma'^p$ and find every enlargement.
Proposition~\ref{v212:elementary-saturation}
bounds their number by $v_p(B)$.
Integer basis changes, together with the root equations, give
\begin{equation}\label{root170:z-relations}
 z_j^{p^{b_j}}=\xi_j\prod_i A_i^{m_{ji}},\qquad
 b_j\ge0,\quad\xi_j\in T_\Sigma,\quad m_{ji}\in\Z.
\end{equation}
\noindent Polynomial factorization and integer basis changes are standard
\cite[Chapters~3--6]{CohenANT}.

\section{Selection by reciprocity in the ordinary case}
\label{root170:selection}

\noindent Selection by Artin reciprocity determines the Kummer coordinates
of the target from its prescribed Frobenius action. This is the
classical Artin map method \cite[Section~3]{FiekerArtin},
\cite[Section~5.4]{CohenAdvanced},
\cite[Section~4.2]{FHS}. We use separating primes rather than
compute another ray class group over $E$; compare \cite[Remark~6]{FHS}.
Target membership from Section~\ref{v206:coverage} and completeness
of the Kummer basis are both essential: together they ensure that
these tests determine the prescribed class.

\noindent Let $n=p^a$, $p\nmid2\Delta$, and $E=K(\zeta_n)$. Inflate the
character to $\mathfrak m=\operatorname{lcm}(\mathfrak f,\tau_K\mathfrak f)$
and assume $\chi^\tau=-\chi$; Lemma~\ref{root170:assembly}
handles its invariant part separately. Choose the ordinary support
$\Sigma$ as in Theorem~\ref{v203:ordinary-support}, put
$r=\rank\Lambda_\Sigma$, and choose representatives
$z_1,\ldots,z_r$ from Corollary~\ref{root170:saturation-theorem}.
The basis and target membership mean
\[
 \begin{gathered}
 (\Z_n)^r\xrightarrow{\sim}W_\Sigma,
       \qquad(e_j)_j\longmapsto\prod_{j=1}^r[z_j]^{e_j},\\
                    \nu^{-1}(\chi|_{G_E})\in W_\Sigma.
 \end{gathered}
\]
\noindent Here $\nu([x])=\nu_x$ is the Kummer isomorphism
\eqref{root170:kummer-iso}. Thus the unknown is the
coordinate vector of a specified class in a known finite module.
For $N=E(\sqrt[n]{z_1},\ldots,\sqrt[n]{z_r})$, the perfect
Kummer pairing is
\begin{equation}\label{v225:kummer-pairing}
 \begin{aligned}
 \langle\ ,\ \rangle:\Gal(N/E)\times W_\Sigma&\longrightarrow\mu_n,\\
 (g,[x])&\longmapsto g(\alpha_x)/\alpha_x.
 \end{aligned}
\end{equation}
\noindent Here $\alpha_x\in N$ satisfies $\alpha_x^n=x$; the pairing is independent
of the chosen root and representative. The target
is specified by
\[
 \langle\Frob_{\mathfrak q},\nu^{-1}(\chi|_{G_E})\rangle
       =\zeta_n^{\chi([\Nm_{E/K}\mathfrak q]_{\mathfrak m})}
\]
\noindent at primes where the residue symbols are defined. Call a prime of $E$
\emph{admissible} if it lies outside
\[
 \mathcal T=\Sigma\cup\{\mathfrak p:\mathfrak p\mid n\mathfrak m\OO_E\}.
\]
\noindent For a prime $\mathfrak q\notin\mathcal T$ and $x\in R_\Sigma^\times$,
write $q=|\OO_E/\mathfrak q|$ and define the additive residue symbol
$r_{\mathfrak q}(x)\in\Z_n$ by
\[
 x^{(q-1)/n}\equiv\zeta_n^{r_{\mathfrak q}(x)}\pmod{\mathfrak q}.
\]
\noindent Reduction is injective on $\mu_n$, so $n\mid q-1$ and this exponent
is unique. The homomorphism $r_{\mathfrak q}:R_\Sigma^\times\to\Z_n$
depends only on the Kummer class whenever that class has a supported
unit representative. In particular it defines
$r_{\mathfrak q}:W_\Sigma\to\Z_n$. It evaluates the Kummer
character at arithmetic Frobenius, by
Lemma~\ref{v203:residue-character}.

\begin{lemma}
% Former title: Separating primes
\label{root170:test-group}
\noindent Let $n=p^a$ for any prime $p$, and let $E$ be any number field
containing a fixed primitive $n$th root of unity $\zeta_n$.
Let $V\subset E^\times/E^{\times n}$ be free of rank $d$ over
$\Z_n$, with specified basis $[z_1],\ldots,[z_d]$ represented by
supported units.
Outside any prescribed finite set containing the primes above $n$
and the supports of the $z_j$, there are primes
$\mathfrak q_1,\ldots,\mathfrak q_d$ such that
\[
 V\xrightarrow{\sim}(\Z_n)^d,\qquad
 [x]\longmapsto(r_{\mathfrak q_1}(x),\ldots,r_{\mathfrak q_d}(x)).
\]
\noindent They can be found by enumerating primes and keeping a residue row
whenever it increases the rank modulo $p$.
\end{lemma}
\begin{proof}
\noindent The Kummer pairing identifies
\[
 \Gal\bigl(E(\sqrt[n]{z_1},\ldots,\sqrt[n]{z_d})/E\bigr)
               \simeq\Hom(V,\mu_n)\simeq(\Z_n)^d.
\]
\noindent Chebotarev supplies every Frobenius vector outside the excluded set
\cite[Chapter~VII, Theorem~13.4]{NeukirchANT}. Consequently the
residue rows span modulo $p$, and the procedure reaches rank $d$.
A square matrix invertible modulo $p$ is invertible modulo $p^a$.
For $d=0$ the empty matrix is the unique isomorphism of the zero module.
\end{proof}
\noindent These auxiliary Kummer extensions need not be constructed:
all rows are computed in finite residue fields. The lemma gives
termination, not a bound asserting that the required primes are small.

\begin{definition}
% Former title: The matrix of the Kummer pairing
\label{v203:matrix-definition}
\noindent Enumerate admissible primes of $E$. After examining $s$ primes
$\mathfrak q_1,\ldots,\mathfrak q_s$, the residue-testing map and
its matrix in the basis $[z_1],\ldots,[z_r]$ are
\[
 \begin{gathered}
 W_\Sigma\longrightarrow(\Z_n)^s,\qquad
 [x]\longmapsto(r_{\mathfrak q_k}(x))_{k=1}^s,\\
 R_s:(\Z_n)^r\longrightarrow(\Z_n)^s,\qquad
 R_s=(r_{\mathfrak q_k}(z_j))_{\substack{1\le k\le s\\1\le j\le r}}.
 \end{gathered}
\]
\noindent Continue until $R_s\bmod p$ has rank $r$, and choose row indices
$i_1,\ldots,i_r$ independent modulo $p$. Lemma~\ref{root170:test-group}
guarantees termination. Define the resulting square matrix and target vector by
\begin{equation}\label{root170:exponent-matrix}
 \begin{gathered}
 C=(r_{\mathfrak q_{i_\ell}}(z_j))_{1\le\ell,j\le r}
                         \in\operatorname{GL}_r(\Z_n),\qquad
 b_\ell=\chi([\Nm_{E/K}\mathfrak q_{i_\ell}]_{\mathfrak m}),\\
 \boxed{\;\boldsymbol e=C^{-1}\boldsymbol b\in(\Z_n)^r.\;}
 \end{gathered}
\end{equation}
\noindent Thus $C$ is the matrix of the Kummer pairing at the selected primes:
$\langle\Frob_{\mathfrak q_{i_\ell}},[z_j]\rangle=\zeta_n^{C_{\ell j}}$.
Choose representatives $0\le e_j<n$ and define
\begin{equation}\label{root170:beta-computed}
                \beta(\boldsymbol e)=\prod_{j=1}^r z_j^{e_j}.
\end{equation}
\noindent For $r=0$, the product is $1$ and the coordinate vector is empty.
\end{definition}
\noindent The chosen elements $z_j$ determine $\beta(\boldsymbol e)$ itself.
Their recorded identities~\eqref{root170:z-relations} also express
a power of it in the original Jacobi values:
\begin{equation}\label{v229:ordinary-input-identity}
 \beta(\boldsymbol e)^{p^b}
 =\prod_{j=1}^r\xi_j^{e_jp^{b-b_j}}
    \prod_i A_i^{\sum_{j=1}^r e_jm_{ji}p^{b-b_j}},
 \qquad b=\max_j b_j.
\end{equation}
\noindent Set $b=0$ if $r=0$. This identity records the relations between
the chosen basis representatives and the $A_i$; the product
\eqref{root170:beta-computed} defines $\beta(\boldsymbol e)$ itself.
The ideals and units defining $\chi$ are computed over $K$;
an equation for $H_\chi$ is not an input.

\begin{theorem}
% Former title: Selection by Artin reciprocity

\label{root170:reciprocity-proved}
\noindent The element $\beta_\chi=\beta(\boldsymbol e)$ from
\eqref{root170:exponent-matrix} satisfies
\[
                         \nu_{\beta_\chi}=\chi|_{G_E}.
\]
\noindent Its coordinate vector is the unique solution. If $\chi$ has exact
order $n$, then $[\beta_\chi]$ has order $n$ and
\[
 E(\sqrt[n]{\beta_\chi})=EH_\chi,\qquad H_\chi\cap E=K.
\]
\end{theorem}
\begin{lemma}
% Former title: No loss of an ordinary anti-invariant character

\label{v203:ordinary-intersection}
\noindent For $H=H_\chi$, one has $H\cap E=K$.
\end{lemma}
\begin{proof}
\noindent The stable kernel makes $H/\Q$ Galois. On $\Gal(H\cap E/K)$,
conjugation by $\tau_K$ acts by inversion because $\chi$ is
anti-invariant, and trivially because $E/\Q$ is abelian. This group
has odd order and is therefore trivial.
\end{proof}

\begin{proof}[Proof of Theorem~\ref{root170:reciprocity-proved}]
\noindent Target membership and the basis give a unique vector
$\boldsymbol e_\chi$ representing $\nu^{-1}(\chi|_{G_E})$.
By the residue law and norm compatibility of Artin reciprocity,
$C\boldsymbol e_\chi=\boldsymbol b$. Invertibility therefore gives
$\boldsymbol e=\boldsymbol e_\chi$ and equality of the entire
characters, not merely their values at the chosen primes.
The final assertions follow from Kummer theory and
Lemma~\ref{v203:ordinary-intersection}.
\end{proof}

\begin{corollary}
% Former title: A direct certificate for the raw values

\label{v224:raw-residue-certificate}
\noindent Let $\Gamma=\langle T_\Sigma,A_i\rangle$ and
$r=\rank\Lambda_\Sigma$. The following are equivalent:
\begin{enumerate}[label=\textup{(\roman*)}]
\item $p\nmid[U_\Sigma:\Gamma]$;
\item the classes $[A_i]$ generate $W_\Sigma$;
\item some $r$ values $A_{i_1},\ldots,A_{i_r}$ and $r$
admissible primes have residue matrix invertible modulo $p$.
\end{enumerate}
\noindent Whenever \textup{(iii)} holds, these $r$ values may be used as the
basis representatives $z_j$ in \eqref{root170:exponent-matrix}.
No saturation is needed.
\end{corollary}
\begin{proof}
\noindent The quotient $U_\Sigma/\Gamma$ is finite; multiplication by $p$
is bijective on it precisely when \textup{(i)} holds. As
$p\nmid|T_\Sigma|$, this is equivalent to \textup{(ii)}.
A generating set of the free $\Z_n$-module $W_\Sigma$ contains
$r$ elements whose reductions form a basis modulo $p$; their
classes then form a basis over $\Z_n$. Lemma~\ref{root170:test-group}
gives \textup{(iii)}. Conversely, an invertible residue matrix
makes the map $(\Z_n)^r\to W_\Sigma$ injective. Both groups have
order $n^r$, proving \textup{(ii)}.
\end{proof}
\begin{example}
% Former title: Selection in a two-element Kummer basis
\label{v215:small-residue-matrix}
\noindent The quintic construction in Section~\ref{v221:ordinary-quintic}
uses $K=\Q(\sqrt2)$, $n=5$ and $(z_1,z_2)=(A_1,A_2)$.
The primes
\[
 \mathfrak q_1=(31,\sqrt2-8,\zeta_5-2),\qquad
 \mathfrak q_2=(31,\sqrt2-8,\zeta_5-4)
\]
\noindent of $E$ both lie above $(7+3\sqrt2)$ and give
\[
 C=\begin{pmatrix}2&0\\0&4\end{pmatrix},\qquad
 \boldsymbol b=\binom44,\qquad
 C^{-1}\boldsymbol b=\binom21\quad\text{in }\F_5,
 \qquad\beta_\chi=A_1^2A_2.
\]
\noindent Their rows separate the rank-two module $W_\Sigma$; they need not
generate the auxiliary ray class group. The example also shows that
these Jacobi values already represent a Kummer basis.
\end{example}

\section{Ordinary generator by equivariant descent}
\label{root170:ordinary-descent}

\noindent Let $n=p^a$, $p\nmid2\Delta$, and let $\chi$ be an
anti-invariant ray class character of exact order $n$.
Section~\ref{root170:selection} gives
$\beta_\chi=\prod_jz_j^{e_j}\in U_\Sigma$ with
$\nu_{\beta_\chi}=\chi|_{G_E}$, where $E=K(\zeta_n)$.
Put $\alpha^n=\beta_\chi$ and $L=E(\alpha)=EH_\chi$.
By Lemma~\ref{v203:ordinary-intersection}, $[L:E]=n$.
For $b\in\Z_n$ define
\[
 r_b:L\longrightarrow L,\qquad
 r_b|_E=\id,\quad r_b(\alpha)=\zeta_n^b\alpha.
\]
\noindent Thus $\Spec L$ over $E$ is a torsor under the constant cyclic group $\Z_n$:
on its $n$ geometric points this group acts freely and transitively.
To recover the desired field over $K$, we descend this torsor
\emph{equivariantly}, meaning that the semilinear descent maps commute
with every $r_b$. This is standard Galois descent
\cite[Section~2.3]{GilleSzamuely}; its computational Kummer form is
\cite[Lemma~4.1]{FiekerArtin}, \cite[Section~4.3]{FHS}.

\noindent Represent $u\in\Z_n^\times$ by $1\le u<n$, and let
$[uv]_n$ denote the representative of its product. Recall
$\delta_u|_K=\id$, $\delta_u(\zeta_n)=\zeta_n^u$.

\begin{theorem}
% Former title: The ordinary cyclic generator

\label{root170:ordinary-descent-theorem}\label{root170:trace-generator}
\noindent There is a unique $c_u\in E^\times$ satisfying
\[
 c_u^n=\frac{\delta_u(\beta_\chi)}{\beta_\chi^u},
                  \qquad c_u\tau(c_u)=1.
\]
\noindent The maps
\begin{equation}\label{root170:du}
 d_u:L\longrightarrow L,\qquad
 d_u|_E=\delta_u,\quad d_u(\alpha)=c_u\alpha^u
\end{equation}
\noindent are the equivariant descent to $H_\chi$. In particular,
\begin{equation}\label{root170:trace}
 \boxed{\theta_\chi=\sum_{u\in\Z_n^\times}c_u\alpha^u
       =\Tr_{L/H_\chi}(\alpha),\qquad K(\theta_\chi)=H_\chi.}
\end{equation}
\end{theorem}
\noindent The exponent $u$ in \eqref{root170:du} is forced by equivariance:
$d_ur_b(\alpha)=r_bd_u(\alpha)=\zeta_n^{ub}c_u\alpha^u$.
The extra norm under $\tau$ chooses the required descent uniquely;
in the general case we will choose it by its Artin character.

\subsection{The proof}
\begin{lemma}
% Former title: The normalized root
\label{root170:unique-root}
\noindent If $x\in E^{\times n}$ and $\tau(x)=x^{-1}$, there is a unique
$c\in E^\times$ with $c^n=x$ and $c\tau(c)=1$.
If $xj(x)=1$, then $cj(c)=1$ as well.
\end{lemma}
\begin{proof}
\noindent For $b^n=x$, put $\omega=b\tau(b)\in\mu_n$ and
$c=b\omega^{(n-1)/2}$. This gives the required power and norm.
Two such roots differ by $\xi\in\mu_n$ with $\xi^2=1$, so
are equal. Finally $cj(c)$ is totally positive and its $n$th
power is one.
\end{proof}
\noindent Covariance (Lemma~\ref{v203:covariance}) makes
$\delta_u(\beta_\chi)/\beta_\chi^u$ an $n$th power inverted by
$\tau$. Thus polynomial factorization over $E$ and this lemma
compute $c_u$ without knowing $H_\chi$.

\noindent For representatives $z_j$ lifting a $\Z$-basis of
$\Gamma^{\rm sat}/T_\Sigma$ as in
Corollary~\ref{root170:saturation-theorem},
$\Gal(E/\Q)$-stability of $\Gamma^{\rm sat}$ gives an alternative:
\[
 \delta_u(z_j)=\omega_{u,j}\prod_\ell z_\ell^{S_{u,\ell j}},
 \qquad S_u\in\operatorname{GL}_r(\Z),\quad\omega_{u,j}\in T_\Sigma.
\]
\noindent Valuations give $S_u$, and division gives $\omega_{u,j}$.
Covariance implies
$\boldsymbol q_u=(S_u\boldsymbol e-u\boldsymbol e)/n\in\Z^r$.
Since $p\nmid|T_\Sigma|$, there is a unique $\omega_u\in T_\Sigma$
with $\omega_u^n=\prod_j\omega_{u,j}^{e_j}$. Hence
\begin{equation}\label{root170:coefficient-word}
 c_u=\omega_u\prod_jz_j^{q_{u,j}}.
\end{equation}
\noindent Thus both $\beta_\chi$ and $c_u$ are explicit products of the
same basis representatives $z_j\in E^\times$.

\begin{lemma}
% Former title: Composition
\label{v203:ordinary-cocycle}
\noindent The normalized coefficients satisfy $c_1=1$ and
\[
 c_{[uv]_n}=\delta_u(c_v)c_u^v
                 \beta_\chi^{(uv-[uv]_n)/n}.
\]
\end{lemma}
\begin{proof}
\noindent The right side has $n$th power
$\delta_{uv}(\beta_\chi)/\beta_\chi^{[uv]_n}$ and norm one
under $\tau$. Apply Lemma~\ref{root170:unique-root}.
\end{proof}

\begin{proof}[Proof of Theorem~\ref{root170:ordinary-descent-theorem}]
\noindent The power and composition equations give a group
$\mathcal D=\{d_u\}$ commuting with $\{r_b\}$, and
$\Gal(L/K)=\langle r_1\rangle\times\mathcal D$.
To identify its fixed field, extend $\tau$ by
$\widehat\tau(\alpha)=\alpha^{-1}$. The normalization gives
$\widehat\tau d_u=d_u\widehat\tau$ and
$\widehat\tau r_b\widehat\tau^{-1}=r_{-b}$.
The subgroup fixing $H_\chi$ is a complement, of the form
$\{r_{b(u)}d_u\}$. Since anti-invariance makes $H_\chi/\Q$
Galois, this complement is $\widehat\tau$-stable. Thus
$b(u)=-b(u)$; oddness of $n$ gives $b(u)=0$. Its fixed field is
therefore exactly $H_\chi$.

\noindent The sum in \eqref{root170:trace} is invariant under all $d_u$.
If $r_b$ fixes it, compare its coefficient of $\alpha$ in the
$E$-basis $1,\alpha,\ldots,\alpha^{n-1}$. Since $c_1=1$, this
coefficient gives $\zeta_n^b=1$. Thus the trace has trivial
stabilizer in $\Gal(H_\chi/K)$ and is primitive.
\end{proof}
\noindent For an arbitrary ordinary character, recombine its cyclotomic
and anti-invariant parts by Lemma~\ref{root170:assembly}.
The general descent in Section~\ref{root170:general-algebra}
uses the same equivariant principle when the norm normalization
or the equality $H_\chi\cap E=K$ is unavailable.

\section{General Fermat Jacobians}
\label{geo170:universal}\label{v216:change-of-curve}

\noindent As in Section~\ref{geo170:curves}, we take Frobenius determinants on
RM eigenspaces of a Fermat quotient. The curve is now
$C_{M,\Delta}:Y^M+Z^\Delta=1$, with $\Delta\mid M$.
The ordinary curve is a quotient of this enlarged curve, and the
same Gaussian endomorphism supplies multiplication by $K$.
The essential change is the normalization:
\[
 \frac{B_v^+}{B_v^-}\quad\longrightarrow\quad
 \frac{B_v^+}{j(B_v^+)}.
\]
\noindent The first quotient compares the two RM signs; the second compares
one RM determinant with its complex conjugate. Only the second
always compares equal-dimensional spaces on the enlarged curve.

\noindent We define these determinants geometrically, evaluate them by
Jacobi sums at good primes, and give explicit CM-ideal inputs at
bad primes. The latter are necessary for the uniform computational
procedure given here; they avoid evaluating the exact bad-prime
Frobenius values. Section~\ref{v222:general-spanning} proves that
the resulting divisors span, just as
Section~\ref{geo170:ordinary-spanning} does in the ordinary case.
The basis, reciprocity and descent steps then have the same form.

\subsection{The curve, its endomorphisms, and its determinant line}
\label{v223:uniform-curve}\label{v223:uniform-RM}
\noindent Let $E/\Q$ be the abelian CM field chosen in
Section~\ref{exc170:coefficientfields}, and choose
\begin{equation}\label{v222:RM-level}
 \Delta\mid M,\qquad E\subseteq F=\Q(\zeta_M),\qquad t=M/\Delta.
\end{equation}
\noindent One may take $M=\lcm(\Delta,n)$ at odd $p$, and
$M=\lcm(\Delta,2^b)$ at two, with $b$ as in
\eqref{exc170:dyadlevel}. Fix $\zeta_\Delta=\zeta_M^t$.
Our curve is the smooth projective model of
\begin{equation}\label{v229:uniform-affine}
                  C=C_{M,\Delta}:\quad Y^M+Z^\Delta=1,
\end{equation}
\noindent and $J=\Jac(C)$. Its Fermat cover is
\begin{equation}\label{v223:uniform-cover}
 \mathcal F_M:U^M+V^M=1\longrightarrow C,
                   \qquad(U,V)\longmapsto(U,V^t).
\end{equation}
\noindent The maps
\[
 S:C_F\to C_F,\ (Y,Z)\mapsto(\zeta_MY,Z),\qquad
 T:C_F\to C_F,\ (Y,Z)\mapsto(Y,\zeta_\Delta Z)
\]
\noindent act on $J_F$ by pushing forward degree-zero divisors. No point at
infinity is chosen. Over $\mathbf P^1\setminus\{0,1,\infty\}$,
$x=Y^M$ is a cover with deck group $\mu_M\times\mu_\Delta$.

\noindent Choose $\ell_0$ outside the rational primes below $\Sigma$ and fix
$\iota:F\hookrightarrow\overline\Q_{\ell_0}$.
On $V=V_{\ell_0}(J)\otimes_{\Q_{\ell_0}}{\overline\Q_{\ell_0}}$, for
$(v,b)\in\Z_M\times\Z_\Delta$, put
\[
 \begin{split}
 V_{v,b}&=\{w\in V:Sw=\iota(\zeta_M^v)w,
                    \ Tw=\iota(\zeta_\Delta^b)w\},\\
 e_v^Y&=\frac1M\sum_{a=0}^{M-1}\iota(\zeta_M^{-va})S^a,
 \qquad e_b^Z=\frac1\Delta\sum_{a=0}^{\Delta-1}
                             \iota(\zeta_\Delta^{-ba})T^a.
 \end{split}
\]
\noindent These are the projectors of \eqref{geo170:projectors} with the
first level enlarged: $e_v^Ye_b^ZV=V_{v,b}$. The spaces are
$G_F$-stable because $S,T$ are defined over $F$.
Select primitive second indices with \eqref{geo170:primitive}:
\[
 J_0=J_{M,\Delta}^{\rm prim}
       =\operatorname{im}(\Delta e_\Delta^{\rm prim}:J\to J).
\]
\noindent The rational projector is fixed by coefficient Galois action,
so $J_0$ is defined over $\Q$. A pair $(v,b)$ is \emph{admissible} if
\begin{equation}\label{v222:RM-admissible}
 \boxed{\quad v\ne0\text{ in }\Z_M,\qquad
 b\in\Z_\Delta^\times,\qquad v+tb\ne0\text{ in }\Z_M.\quad}
\end{equation}
\noindent Precisely these pairs give the one-dimensional character spaces of
$J_0$ (Lemma~\ref{v223:uniform-multiplicity}). All nonzero first
indices are used; restricting $v$ to units can lose a divisor
component.% (Example~\ref{v223:nonprimitive-example}).

\noindent The Gaussian endomorphism is unchanged:
\begin{equation}\label{v223:uniform-K-action}
 \mathsf G_\Delta=\sum_{a\in\Z_\Delta^\times}\chi_\Delta(a)T^a,
 \qquad K\hookrightarrow\End_K^0(J_0),\quad
             \sqrt\Delta\longmapsto\mathsf G_\Delta.
\end{equation}
\noindent Indeed Lemma~\ref{geo170:RM} gives
$\mathsf G_\Delta|_{V_{v,b}}=\chi_\Delta(b)\iota(\sqrt\Delta)$
and $\mathsf G_\Delta^2=[\Delta]$. The primitive
$\Q(\zeta_\Delta)$-action is faithful, and coefficient equivariance
descends the quadratic endomorphism to $K$, exactly as there.
For $1\le v<M$, define
\begin{equation}\label{v222:RM-blocks}
 \begin{gathered}
 I_v^\pm=\{b\in\Z_\Delta^\times:
       \chi_\Delta(b)=\pm1,\ v+tb\ne0\pmod M\},\qquad r_v^\pm=|I_v^\pm|,\\
 M_v^\pm=\ker(\mathsf G_\Delta\mp\iota(\sqrt\Delta)
                    \mid e_v^Ye_\Delta^{\rm prim}V)
             =\bigoplus_{b\in I_v^\pm}V_{v,b}.
 \end{gathered}
\end{equation}
\noindent These are subspaces of the rational Tate module after extending
scalars; they are not asserted to be individual abelian subvarieties.
Their projectors are $\tfrac12(1\pm\mathsf G_\Delta/\iota(\sqrt\Delta))$.

\noindent Put $r=\varphi(\Delta)/2$. Both dimensions equal $r$ unless
$v+tb=0$ for a primitive $b$; then one equals $r-1$ and the other
$r$. For example, on $C_{8,8}$ with $v=1$, the positive index
$b=7$ is excluded, so $(r_1^+,r_1^-)=(1,2)$. The ordinary
quotient of determinants would have absolute value $q^{-1/2}$.
In contrast, $\chi_\Delta(-1)=1$ gives
$I_{-v}^\pm=-I_v^\pm$ and $r_{-v}^\pm=r_v^\pm$. We therefore use
\begin{equation}\label{v223:uniform-det-line}
 \mathscr L_v^{(1)}=
       \det M_v^+\otimes(\det M_{-v}^+)^\vee,
 \qquad \det M_v^+=\bigwedge^{r_v^+}M_v^+.
\end{equation}
\noindent Here $\det0={\overline\Q_{\ell_0}}$. This one-dimensional representation has
weight zero. The negative RM sign gives another valid family;
either sign alone suffices for spanning. Complex conjugation here
acts on the algebraic coefficients, not as an ${\overline\Q_{\ell_0}}$-linear
operator on $V$.

\subsection{The Frobenius determinant}
\label{v223:uniform-Frobenius}
\noindent For each rational prime $\ell$ below $\Sigma$, choose
$L_\ell\supseteq F$ over which $J_0$ has good reduction and a prime
$\mathfrak P_\ell\mid\ell$. Set
\[
 \mathfrak q_\ell=\mathfrak P_\ell\cap\OO_F,
 \qquad q_\ell=|\OO_{L_\ell}/\mathfrak P_\ell|.
\]
\noindent At $\ell\nmid M$ one may take $L_\ell=F$; at every prime one
may use Lemma~\ref{v222:potential-good}. Write $\Frob_{\mathfrak P_\ell}$ for the
$q_\ell$-power Frobenius endomorphism of the reduction, acting on
Tate modules by arithmetic Frobenius as in
\eqref{geo170:frob-convention}.
At a bad prime this means Frobenius on the restriction to
$G_{L_\ell}$ after killing inertia. It is not a canonical
Frobenius lift on the full ramified $G_F$-representation, nor
necessarily the operator on inertia invariants defining its local
Euler factor over $F$.

\begin{definition}
% Former title: General RM determinants
\label{v222:RM-values}
\noindent Define $B_{\mathfrak q_\ell}^\pm(v)\in F$ by
\begin{equation}\label{v229:RM-determinant-definition}
 \boxed{\quad
 \iota(B_{\mathfrak q_\ell}^\pm(v))=
       \det_{\overline\Q_{\ell_0}}(\Frob_{\mathfrak P_\ell}\mid M_v^\pm).
 \quad}
\end{equation}
\noindent Thus $\Frob_{\mathfrak P_\ell}$ acts on $\det M_v^\pm$ by
$\iota(B_{\mathfrak q_\ell}^\pm(v))$. The prime subscript is in $F$,
as in Section~\ref{geo170:curves}; the chosen pair
$(L_\ell,\mathfrak P_\ell)$ is understood throughout. The algebraicity of these determinants
follows from Lemma~\ref{v222:potential-good} and the CM algebra
below. For the positive sign put
\begin{equation}\label{v222:RM-input}
 A_{\mathfrak q_\ell}(v)=N_{F/E}\!\left(\frac{B_{\mathfrak q_\ell}^+(v)}{j(B_{\mathfrak q_\ell}^+(v))}\right),
\end{equation}
\noindent and
\begin{equation}\label{v217:general-index-set}
 \mathcal A_\Sigma=\{A_{\mathfrak q_\ell}(v):1\le v<M,
                    \ \ell\text{ lies below }\Sigma\}.
\end{equation}
\noindent An empty determinant is $1$.
\end{definition}
\noindent The quotient is the Frobenius scalar on
$\mathscr L_v^{(1)}$; $A_{\mathfrak q_\ell}(v)$ is its field norm.
In Section~\ref{geo170:curves} the ordinary determinants already
lie in $E$, so no further norm is needed. Here $\Gal(F/E)$ may
move the first index, and $N_{F/E}$ puts the input in $E$.
This is a coefficient norm, not a CM reflex norm.
Equivalently, identifying the determinant scalars with elements of
$F$ through $\iota$,
\[
 A_{\mathfrak q_\ell}(v)=N_{F/E}\!\left(
 \frac{\det(\Frob_{\mathfrak P_\ell}\mid M_v^+)}
      {\det(\Frob_{\mathfrak P_\ell}\mid M_{-v}^+)}\right).
\]

\begin{theorem}
% Former title: The uniform determinant coefficient

\label{v223:uniform-determinant}
\noindent On the determinant line,
\[
 \Frob_{\mathfrak P_\ell}\big|_{\mathscr L_v^{(1)}}
    =\iota\!\left(\frac{B_{\mathfrak q_\ell}^+(v)}{j(B_{\mathfrak q_\ell}^+(v))}\right)\id.
\]
\noindent For every $1\le v<M$ and support prime $\ell$,
\[
 B_{\mathfrak q_\ell}^+(-v)=j(B_{\mathfrak q_\ell}^+(v)),\qquad
 B_{\mathfrak q_\ell}^+(v)j(B_{\mathfrak q_\ell}^+(v))=q_\ell^{r_v^+},
\]
\noindent and
\[
 A_{\mathfrak q_\ell}(v)\in U_\Sigma^{(1)},\qquad
 j(A_{\mathfrak q_\ell}(v))=A_{\mathfrak q_\ell}(v)^{-1},\qquad
 A_{\mathfrak q_\ell}(-v)=A_{\mathfrak q_\ell}(v)^{-1}.
\]
\noindent Its divisor is supported above $\ell$. These assertions hold also
when $\ell\mid M$, using the specified good-reduction extension.
\end{theorem}
\noindent The proof, like that of Theorem~\ref{geo170:detquotient}, uses
the Fermat character decomposition. CM supplies a common
Frobenius at every prime; the individual eigenvalues below are
used to evaluate this determinant and its divisor.

\subsection{Character decomposition and proof of the determinant theorem}
\label{v229:uniform-proof}\label{v223:uniform-CM}
\begin{lemma}
% Former title: The character spaces
\label{v223:uniform-multiplicity}
\noindent The space $V_{v,b}$ is one-dimensional precisely when
\[
 v\ne0\text{ in }\Z_M,\qquad b\ne0\text{ in }\Z_\Delta,
                 \qquad v+tb\ne0\text{ in }\Z_M;
\]
\noindent otherwise it is zero. These spaces exhaust $V$, and
\[
 2g(C_{M,\Delta})=\dim_{\overline\Q_{\ell_0}} V=(M-2)(\Delta-1),
 \qquad 2\dim J_0=(M-2)\varphi(\Delta).
\]
\end{lemma}
\begin{proof}
\noindent The cover \eqref{v223:uniform-cover} identifies $(v,b)$ with
Fermat indices $(v,tb)$. The classical multiplicity-one theorem
applies exactly when both indices and their sum are nonzero
\cite[Section~1]{KoblitzRohrlich}
\cite[Propositions~4.8 and 4.11]{OY}. For each $b\ne0$ the two
excluded first indices are $0$ and $-tb$. Counting all second
indices, or only primitive ones, gives the two dimensions.
\end{proof}

\noindent The CM algebra is $\mathcal D=\Q[S,T]|_{J_0}$.
For an admissible pair, its character field is
\begin{equation}\label{v222:RM-CM-field}
 F_{v,b}=\Q(\zeta_M^v,\zeta_\Delta^b)
        =\Q\bigl(\zeta_{M/\gcd(v,t)}\bigr)\subseteq F,
\end{equation}
\noindent since $\lcm(M/\gcd(M,v),\Delta)=M/\gcd(v,t)$.
Writing $\sigma_u(\zeta_M)=\zeta_M^u$, coefficient conjugation
sends $(v,b)$ to $(uv,ub)$. Choose one pair in each such orbit
$\mathcal O$, and put $F_{\mathcal O}=F_{v,b}$. Character
orthogonality gives
\[
 \mathcal D\simeq\prod_{\mathcal O}F_{\mathcal O},\qquad
 \dim_\Q\mathcal D=2\dim J_0,
 \qquad V_{\ell_0}(J_0)\text{ is free of rank one over }
                      \mathcal D\otimes\Q_{\ell_0}.
\]
\noindent Every field $F_{\mathcal O}$ contains $\Q(\zeta_\Delta)$ and
hence $K$. The isogeny factors need not be simple; only the
rank-one CM action is used.

\begin{lemma}
% Former title: A common good-reduction field
\label{v222:potential-good}
\noindent The variety $J_0$ has good reduction at every finite prime of
$L=F(J_0[15])$. Its residue-field Frobenius belongs to the
specialized algebra $\mathcal D$.
\end{lemma}
\begin{proof}
\noindent CM gives potential good reduction \cite[Theorem~7]{SerreTate1968}.
At residue characteristic $\ell$, choose
$r\in\{3,5\}\setminus\{\ell\}$. Over $L$, finite inertia acts
trivially modulo $r$ on $T_r(J_0)$. The kernel
$1+rM_{2\dim J_0}(\widehat\Z_r)$ is torsion-free by the $r$-adic
logarithm, so inertia is trivial. N\'eron--Ogg--Shafarevich gives
good reduction \cite[Theorem~1]{SerreTate1968}.

\noindent The deck endomorphisms specialize faithfully and are defined over
the residue field. Frobenius therefore centralizes $\mathcal D$.
The rank-one property makes its Tate-module centralizer
$\mathcal D\otimes\Q_r$; hence the centralizer in the rational
endomorphism algebra has dimension at most $\dim_\Q\mathcal D$.
Since it contains $\mathcal D$, it equals $\mathcal D$.
\end{proof}
\noindent This also applies at primes dividing $2\Delta$:
$\mathsf G_\Delta^2=[\Delta]$ survives specialization and
$[\Delta]$ is still a nonzero isogeny. Neither smooth reduction
of the written affine curve nor separability of $[\Delta]$ is
required.

\noindent For an admissible pair, define the one-dimensional determinant
$B_{\mathfrak q_\ell}(v,b)$ by evaluating
$\Frob_{\mathfrak P_\ell}\in\mathcal D$ on $V_{v,b}$:
\begin{equation}\label{v222:RM-eigenvalues}
 \begin{gathered}
 B_{\mathfrak q_\ell}(v,b)\in\OO_{F_{v,b}},\qquad
 \Frob_{\mathfrak P_\ell}|_{V_{v,b}}
      =\iota(B_{\mathfrak q_\ell}(v,b))\id,\\
 B_{\mathfrak q_\ell}(v,b)j(B_{\mathfrak q_\ell}(v,b))=q_\ell.
 \end{gathered}
\end{equation}
\noindent All complex absolute values are $q_\ell^{1/2}$, and
$\sigma_u(B_{\mathfrak q_\ell}(v,b))=B_{\mathfrak q_\ell}(uv,ub)$, with the reduction
prime fixed. Consequently
\[
                 B_{\mathfrak q_\ell}^\pm(v)=
                         \prod_{b\in I_v^\pm}B_{\mathfrak q_\ell}(v,b).
\]
\noindent This is the character-space evaluation of
\eqref{v229:RM-determinant-definition}.
For fixed $(L_\ell,\mathfrak P_\ell)$, these algebraic values are
independent of the auxiliary prime $\ell_0$: they evaluate the
intrinsic endomorphism $\Frob_{\mathfrak P_\ell}\in\mathcal D$
at the algebraic deck characters. The embedding $\iota$ realizes
those same characters on the Tate module. Consequently
$B_{\mathfrak q_\ell}^\pm(v)$ and $A_{\mathfrak q_\ell}(v)$
do not acquire an additional $\ell_0$-parameter.

\begin{proof}[Proof of Theorem~\ref{v223:uniform-determinant}]
\noindent Coefficient conjugation sends $I_v^+$ to $I_{-v}^+$, giving
$j(B_{\mathfrak q_\ell}^+(v))=B_{\mathfrak q_\ell}^+(-v)$. The product of the Weil
identities in \eqref{v222:RM-eigenvalues} gives
$B_{\mathfrak q_\ell}^+(v)j(B_{\mathfrak q_\ell}^+(v))=q_\ell^{r_v^+}$. Thus
\[
 \frac{B_{\mathfrak q_\ell}^+(v)}{j(B_{\mathfrak q_\ell}^+(v))}=\frac{B_{\mathfrak q_\ell}^+(v)^2}{q_\ell^{r_v^+}}.
\]
\noindent This quotient has norm one under $j$. Taking determinants on \eqref{v223:uniform-det-line} gives its
Frobenius scalar. Integrality and the Weil norm restrict every
individual divisor to primes above $\ell$. Norm to $E$ preserves
that support and commutes with $j$, proving the assertions.
\end{proof}
\noindent More generally
$\sigma_u(B_{\mathfrak q_\ell}^\varepsilon(v))
=B_{\mathfrak q_\ell}^{\varepsilon\chi_\Delta(u)}(uv)$.
Enlarging the residue field by degree $e$ raises every determinant
to its $e$th power. The normalized coefficients may consequently
change by powers; their saturated Kummer module will not.

\subsection{Jacobi evaluation at good primes}
\noindent At $\ell\nmid M$, take $L_\ell=F$, put
$k=\OO_F/\mathfrak q_\ell$, $q=|k|$, and define
$\lambda:k^\times\to\mu_M(F)$ by
$\overline{\lambda(x)}=x^{(q-1)/M}$. Extend its character powers
by zero at zero. With the convention \eqref{geo170:jacobi},
\begin{equation}\label{v222:RM-Jacobi}
 J_{\mathfrak q_\ell}(v,tb)
 =J(\lambda^v,\lambda^{tb})
 =\sum_{x\in k}\lambda(x)^v\lambda(1-x)^{tb}.
\end{equation}
\noindent For an admissible pair,
\[
                  B_{\mathfrak q_\ell}(v,b)=-J_{\mathfrak q_\ell}(v,tb)
\]
\noindent by the same Frobenius--Jacobi formula as in
Lemma~\ref{geo203:Jacobi-Frobenius}
\cite[Theorem~A]{ColemanFermat}\cite[Theorem~1.2 and Section~6]{OY}.
Thus the determinant input is
\begin{equation}\label{v227:uniform-Jacobi-input}
 A_{\mathfrak q_\ell}(v)=N_{F/E}\left(
     \prod_{b\in I_v^+}
       \frac{J_{\mathfrak q_\ell}(v,tb)}
            {j(J_{\mathfrak q_\ell}(v,tb))}\right).
\end{equation}
\noindent The determinant signs cancel. For explicit evaluation, choose a
generator $g$ of $k^\times$ and $a\in\Z_M^\times$ with
$g^{(q-1)/M}=\bar\zeta_M^a$; then
$\lambda(g^e)=\zeta_M^{ae}$. The first index need not be
primitive: the two characters and their product need only be
nontrivial.

\subsection{Explicit input at bad primes}
\label{v222:principalized-input}
\noindent At $\ell\mid M$, the finite-field character of order $M$ used
above does not exist. This does not rule out smaller-conductor
Jacobi descriptions or ramified local formulas. Our uniform
algorithm uses the CM ideal of the determinant and computes
an equivalent input in $E$; it does not compute its exact
Frobenius eigenvalue.

\noindent Put $G=(\Z/M\Z)^\times$. For an admissible pair, define
\begin{equation}\label{v222:RM-CM-type}
 h_{v,b}(u)=
 \begin{cases}
 1,&\{uv/M\}+\{ub/\Delta\}<1,\\
 0,&\{uv/M\}+\{ub/\Delta\}>1,
 \end{cases}\qquad u\in G.
\end{equation}
\noindent Equality cannot occur by \eqref{v222:RM-admissible}. These are
the CM-type indicators; Section~\ref{v223:general-CM-divisors}
explains their geometric origin and proves the ideal formula.
For the chosen $\mathfrak q_\ell\mid\ell$, compute
\begin{equation}\label{v228:principalized-ideal}
 \begin{aligned}
 H_{v,\ell}&=\sum_{u\in G}\sum_{b\in I_v^+}
            h_{v,b}(u)[\sigma_u^{-1}\mathfrak q_\ell],\\
 \mathfrak a_{v,\ell}&=\prod_{u\in G}\prod_{b\in I_v^+}
            (\sigma_u^{-1}\mathfrak q_\ell)^{h_{v,b}(u)}.
 \end{aligned}
\end{equation}
\noindent Repeated primes are combined. Proposition~\ref{v222:RM-divisor}
gives the determinant formula
\[
 \divi_F B_{\mathfrak q_\ell}^+(v)=m_\ell H_{v,\ell},\qquad
 m_\ell=f(\mathfrak P_\ell/\mathfrak q_\ell).
\]
\noindent The divisor $H_{v,\ell}$ is computable without knowing $m_\ell$
or constructing $L_\ell$.

\noindent Compute the order $h$ of $[\mathfrak a_{v,\ell}]$ in
$\Cl(\OO_F)$ and a generator of its principal power
\cite[Chapters~6--7]{CohenANT}. The explicit input is
\begin{equation}\label{v227:principalized-input-data}
 \boxed{\quad
 (x)=\mathfrak a_{v,\ell}^{h},\qquad
 A'_{\mathfrak q_\ell}(v)=N_{F/E}\left(\frac{x}{j(x)}\right)
                  \in U_\Sigma^{(1)}.
 \quad}
\end{equation}
\noindent It is an exact algebraic number. Ideal factorization, class-group
and principal-ideal arithmetic in $F$ suffice; no root-of-unity determination,
good-reduction model or residue-degree multiplier is needed.
The procedure also works at good primes.

\begin{proposition}
% Former title: Replacement by a principalized CM ideal

\label{v223:principalized-input-proposition}
\noindent Let $h>0$ and $x\in F^\times$ satisfy $\divi_F(x)=hH_{v,\ell}$,
and define $A'_{\mathfrak q_\ell}(v)$ by \eqref{v227:principalized-input-data}.
Then, for some $\xi\in\mu(F)$,
\begin{equation}\label{v222:principalized-RM}
 \left(\frac{x}{j(x)}\right)^{m_\ell}
           =\xi \left(\frac{B_{\mathfrak q_\ell}^+(v)}{j(B_{\mathfrak q_\ell}^+(v))}\right)^h,\qquad
 (A'_{\mathfrak q_\ell}(v))^{m_\ell}=N_{F/E}(\xi)A_{\mathfrak q_\ell}(v)^{h}.
\end{equation}
\noindent In particular,
$\divi_E A'_{\mathfrak q_\ell}(v)=(h/m_\ell)\divi_E A_{\mathfrak q_\ell}(v)$.
If $h$ is the ideal-class order, then
\begin{equation}\label{v228:principalized-classorder}
 h\mid m_\ell,\qquad
 \frac{B_{\mathfrak q_\ell}^+(v)}{j(B_{\mathfrak q_\ell}^+(v))}=\xi_0\left(\frac{x}{j(x)}\right)^{m_\ell/h},
                       \qquad\xi_0\in\mu(F).
\end{equation}
\noindent Replacing any or all of the $A_{\mathfrak q_\ell}(v)$ by $A'_{\mathfrak q_\ell}(v)$
preserves their rational divisor span. After adjoining $\mu(E)$
and taking internal $p$-saturation, each complete family gives
$W_\Sigma^{(1)}(E,p^a)$ for every $a\ge1$, including $p=2$.
\end{proposition}
\begin{proof}
\noindent The two sides of the first identity, omitting $\xi$, have divisor
$m_\ell h(H_{v,\ell}-jH_{v,\ell})$. Their quotient is a unit
with absolute value one at every embedding, since $F$ is CM and
both normalized elements have norm one under $j$. Kronecker's
theorem makes it a root of unity
\cite[Chapter~I, Proposition~7.1]{NeukirchANT}. Taking the field
norm gives the second identity. For the class order,
$(B_{\mathfrak q_\ell}^+(v))=\mathfrak a_{v,\ell}^{m_\ell}$ implies
$h\mid m_\ell$; the same argument proves
\eqref{v228:principalized-classorder}.

\noindent Each replaced divisor is a nonzero rational multiple of the
original one. Theorem~\ref{geo170:all-directions} therefore gives
finite index in $U_\Sigma^{(1)}$ for the group $\Gamma$ generated
by the chosen family and $\mu(E)$. Proposition~\ref{v212:elementary-saturation}
gives
\[
 \Gamma^{\rm sat}/(\Gamma^{\rm sat})^{p^a}
 \simeq U_\Sigma^{(1)}/(U_\Sigma^{(1)})^{p^a}
 \simeq W_\Sigma^{(1)}(E,p^a).
\]
\noindent This is saturation inside $U_\Sigma^{(1)}$, not division by $h$
or $m_\ell$ modulo $p^a$.
\end{proof}
\noindent Changing $x$ by a unit changes $x/j(x)$ only by a root of unity.
The generated groups, and even their saturated subgroups of
$E^\times$, need not agree; their classes give the same Kummer module. In particular, a $p$-factor in $m_\ell/h$ may
change the index removed by saturation. The full group $\mu(E)$
must be included, also at two. Section~\ref{v222:bad-prime-example}
works out this input for $K=\Q(\sqrt6)$ at $\ell=p=3$.

\begin{remark}
% Former title: Exact eigenvalues and the role of local formulas

\label{v228:exact-CM-values}
\noindent The CM ideal formula and the principalization procedure answer
different questions. The former computes the ideal of the actual
Frobenius; \eqref{v227:principalized-input-data} computes an
explicit substitute whose divisor is a rational multiple of that
of the normalized Frobenius value.
It does not claim to recover the actual eigenvalue.

\noindent For a chosen good-reduction residue cardinality $q_\ell$, the
ideal and Weil relation determine each
$B=B_{\mathfrak q_\ell}(v,b)$ up to $\mu(F_{v,b})$. If $(a)=(B)$,
solve the unit equation
\[
 u\in\OO_{F_{v,b}}^\times,\qquad
 uj(u)=q_\ell/(aj(a)).
\]
\noindent A solution exists because $B/a$ is a unit; unit-group generators
reduce the equation to a calculation with integer relations.
The finite candidate set is $\{\zeta au:\zeta\in\mu(F_{v,b})\}$.
The remaining root of unity requires the local CM character or
an actual Frobenius computation; the ideal alone does not select it.

\noindent Exact ramified formulas do exist. To relate them to our determinants,
let $\rho_{v,b,\mathfrak q_\ell}:F_{\mathfrak q_\ell}^{\times}
\to{\overline\Q_{\ell_0}}^\times$ be the character on $V_{v,b}$ composed with
local reciprocity, normalized to send a uniformizer to arithmetic
Frobenius in the unramified quotient. For a uniformizer $\varpi$ of
$(L_\ell)_{\mathfrak P_\ell}$, compatibility with local norms gives
\[
 \iota(B_{\mathfrak q_\ell}(v,b))
   =\rho_{v,b,\mathfrak q_\ell}
       \bigl(N_{(L_\ell)_{\mathfrak P_\ell}/F_{\mathfrak q_\ell}}
                    (\varpi)\bigr).
\]
\noindent The restriction over $L_\ell$ is unramified, so this value does not
depend on $\varpi$; multiplying over $b\in I_v^+$ gives
$\iota(B_{\mathfrak q_\ell}^+(v))$. Thus a fully normalized local-character
formula computes the exact eigenvalue, including its root of unity.

\noindent Coleman--McCallum give such local formulas at odd residue
characteristic in \cite[Theorem~5.3, pp.~83--85]{ColemanMcCallum1988},
with the additional $3$-adic unit cases treated in
\cite[Theorem~7.2, pp.~97--98]{ColemanMcCallum1988};
Tsushima gives an elementary proof \cite{Tsushima2010}.
For residue characteristic greater than three, crystalline Weil-group
matrices are computed in
\cite[Theorems~1.7 and~3.13]{ColemanMatrices1990}.
Those formulas are exact; the odd-prime period reformulation in
\cite[Section~2.2, Theorem~2.4, pp.~1235--1237]{Kashio2024}
is instead modulo roots of unity.
At two, \cite{Prapavessi1991} gives Hilbert-symbol formulas for
local components at levels $2^r$.
These results must be applied with their stated hypotheses and
reciprocity conventions; a Weil element which permutes character
spaces does not supply an individual eigenvalue merely by taking
one matrix entry. We do not implement these local evaluations at
every mixed level. Our uniform algorithm, including at two, uses
the principal-ideal construction 
\eqref{v227:principalized-input-data}; it computes exact Kummer
inputs without first determining the exact Frobenius values.
\end{remark}

\subsection{Comparison with the ordinary determinant}
\label{v223:uniform-comparison}
\noindent Suppose $p\nmid2\Delta$, $n=p^a$, and $M=\Delta n$.
The ordinary curve is a quotient of the uniform one:
\[
\begin{tikzcd}[column sep=large]
 \mathcal F_M \arrow[r,"\rho"]\arrow[dr,swap,"q\circ\rho"]
      & C_{M,\Delta}\arrow[d,"q"]\\
      & C_{n,\Delta}.
\end{tikzcd}
\]
\noindent Here $\rho(U,V)=(U,V^n)$ and
\begin{equation}\label{v222:ordinary-quotient}
 q:C_{M,\Delta}\longrightarrow C_{n,\Delta},
            \qquad(Y,Z)\longmapsto(Y^{M/n},Z).
\end{equation}
\noindent Since $q_*q^*=[M/n]$, pullback on Jacobians is injective up to
isogeny. It commutes with $T$ and $\mathsf G_\Delta$ and identifies
the ordinary line $(v,b)$ with $(\Delta v,b)$.
For $v\in\Z_n^\times$ no such pair is excluded. In terms of the
uniform RM blocks, the two determinant lines are
\[
 \mathscr L_v=\det M_{\Delta v}^+
                  \otimes(\det M_{\Delta v}^-)^\vee,
 \qquad
 \mathscr L_{\Delta v}^{(1)}=\det M_{\Delta v}^+
                  \otimes(\det M_{-\Delta v}^+)^\vee.
\]
\noindent The uniform construction also uses the other nonzero first
indices and spans all $\Lambda_\Sigma^{(1)}\otimes\Q$; the ordinary
one addresses its additional $\tau$-anti-invariant part.
In this subsection only, superscripts $\mathrm{ord}$ and
$\mathrm{gen}$ distinguish the two inputs $A_{\mathfrak q}(v)$.

\begin{proposition}
% Former title: Recovery of the ordinary coefficients

\label{v215:ordinary-general-comparison}\label{v215:geometry-comparison}
\noindent At a prime $\mathfrak q\mid\ell$ of $F$ with $\ell\nmid M$,
use $\Frob_{\mathfrak q}$ in both constructions. Then
$B_{\mathfrak q}^\pm(\Delta v)\in E$, and
\[
 A_{\mathfrak q}^{\mathrm{ord}}(v)
       =\frac{B_{\mathfrak q}^+(\Delta v)}{B_{\mathfrak q}^-(\Delta v)},
 \qquad
 A_{\mathfrak q}^{\mathrm{gen}}(\Delta v)
       =N_{F/E}\!\left(
        \frac{B_{\mathfrak q}^+(\Delta v)}{j(B_{\mathfrak q}^+(\Delta v))}
                    \right).
\]
\noindent Moreover,
\begin{equation}\label{v222:ordinary-bridge}
 \begin{gathered}
 \frac{B_{\mathfrak q}^+(\Delta v)j(B_{\mathfrak q}^-(\Delta v))}
      {j(B_{\mathfrak q}^+(\Delta v))B_{\mathfrak q}^-(\Delta v)}
       =A_{\mathfrak q}^{\mathrm{ord}}(v)^2,\\
 [A_{\mathfrak q}^{\mathrm{ord}}(v)]
 =\left[
 \frac{B_{\mathfrak q}^+(\Delta v)j(B_{\mathfrak q}^-(\Delta v))}
      {j(B_{\mathfrak q}^+(\Delta v))B_{\mathfrak q}^-(\Delta v)}
        \right]^{(n+1)/2}
       \quad\text{in }E^\times/E^{\times n}.
 \end{gathered}
\end{equation}
\noindent Selecting the same anti-invariant character in the two constructions gives
\begin{equation}\label{v215:same-kummer-class}
 \beta_\chi^{\rm gen}=b^n\beta_\chi^{\rm ord}
                       \quad\text{for some }b\in E^\times;
\end{equation}
\noindent equivariant descent with that character gives the same cyclic
extension of $K$.
\end{proposition}
\begin{proof}
\noindent The group $\Gal(F/E)$ fixes the pulled-back ordinary index
$\Delta v$ and preserves the RM sign, so it fixes
$B_{\mathfrak q}^\pm(\Delta v)$.
The involution $\tau$ interchanges the signs. Both determinants
have degree $r$, and
\[
 B_{\mathfrak q}^\pm(\Delta v)j(B_{\mathfrak q}^\pm(\Delta v))=q^r,
 \qquad q=|\OO_F/\mathfrak q|.
\]
\noindent Dividing these two identities proves the square relation.
Since $n$ is odd, raising to $(n+1)/2$ gives the class identity.
Uniqueness of the class selected by the character and descent
with that full character give the last assertion.
\end{proof}
\noindent The quotient inside the uniform norm already lies in $E$ on
these ordinary blocks. Its norm is therefore its $[F:E]$th power,
so
\begin{equation}\label{v223:comparison-norm-power}
 \frac{A_{\mathfrak q}^{\mathrm{gen}}(\Delta v)}
      {\tau(A_{\mathfrak q}^{\mathrm{gen}}(\Delta v))}
                    =A_{\mathfrak q}^{\mathrm{ord}}(v)^{2[F:E]}.
\end{equation}
\noindent If $p\mid[F:E]$, this power cannot be inverted modulo $n$;
internal saturation removes the possible loss. For example,
$\Delta=13$, $n=3$ gives $[F:E]=6$. At a prime above $79$,
$A_{\mathfrak q}^{\mathrm{ord}}(1)$ has Jacobi valuations $\pm2$,
but its power in \eqref{v223:comparison-norm-power} is a cube.
The uniform input is the sixth power of the normalized determinant;
its internal cube root is
\[
 \left(\frac{B_{\mathfrak q}^+(\Delta v)}
                  {j(B_{\mathfrak q}^+(\Delta v))}\right)^2.
\]
\noindent Dividing this by its $\tau$-conjugate gives
$A_{\mathfrak q}^{\mathrm{ord}}(v)^4$, representing
$[A_{\mathfrak q}^{\mathrm{ord}}(v)]$ modulo cubes.
Thus different initial values yield compatible selected classes.
At two the same RM determinant applies; signatures and the
coefficient enlargement still require their separate treatment.
For exceptional primes, $\tau_K$ need not extend to $E$ while
fixing $\mu_n$; at two, division by two in the character
decomposition is also unavailable. The ordinary hypothesis
remains $p\nmid2\Delta$.

\section{General RM divisors}
\label{v222:general-spanning}

\noindent We prove that the normalized RM determinants of
Section~\ref{geo170:universal} span the required divisor space.
The explicit inputs constructed in
Section~\ref{v222:principalized-input} have the same span by
Proposition~\ref{v223:principalized-input-proposition}; thus the
result applies to the computation at both good and bad primes.
As in the ordinary case, this gives finite index; the next section
passes to a basis modulo $p^a$.

\noindent Keep $E\subseteq F=\Q(\zeta_M)$, $t=M/\Delta$, and
$\mathcal A_\Sigma$ from Section~\ref{geo170:universal}. Recall
\[
 \Lambda_\Sigma^{(1)}
   =\{D\in\Z[\Sigma]:jD=-D\},\qquad
 U_\Sigma^{(1)}=\{x\in R_\Sigma^\times:xj(x)=1\},\qquad
 \divi_E(x)=\sum_{\mathfrak p\in\Sigma}
                    v_{\mathfrak p}(x)[\mathfrak p].
\]
\noindent The valuations are normalized by
$v_{\mathfrak p}(\mathfrak p)=1$. In contrast to the ordinary
lattice, no second anti-invariance condition is imposed.

\begin{theorem}
% Former title: Uniform RM spanning
\label{geo170:all-directions}
\noindent For every finite $\Gal(E/\Q)$-stable support $\Sigma$,
\begin{equation}\label{v222:full-RM-span}
 \operatorname{span}_\Q\{\divi_E A:A\in\mathcal A_\Sigma\}
                    =\Lambda_\Sigma^{(1)}\otimes\Q.
\end{equation}
\noindent Consequently $\langle\mu(E),\mathcal A_\Sigma\rangle$ has finite
index in $U_\Sigma^{(1)}$. Either RM sign may be used, provided
all nonzero first indices $v$ are used. The statement includes
support above $2$ and above every prime divisor of $M$.
\end{theorem}

\noindent Section~\ref{v222:diagonal-alternative} gives an alternative
using $C_c:z^c=x(1-x)$ and only good-reduction Jacobi sums.
The proof below uses the fixed RM curve and its determinant divisors.

\noindent The proof has three steps.
\begin{enumerate}
\item The Shimura--Taniyama formula computes the Frobenius divisors,
including at ramified primes (Proposition~\ref{v222:RM-divisor}).
\item After normalization, their coefficient functions $d_v$ span
all odd functions on $(\Z_M)^\times$. Additive Fourier transform
reduces this to the nonvanishing of $L(1,\eta)$ for odd characters
(Lemmas~\ref{v223:RM-multiplier}--\ref{v222:RM-Fourier}).
\item Projection to prime divisors and norm to $E$ give
\eqref{v222:full-RM-span} (Section~\ref{v223:general-span-proof}).
\end{enumerate}
\noindent The use of $L$-values at one is essential: their extra Euler factors
never vanish, so imprimitive characters are included. The ordinary
argument at zero needs its more restricted family of characters.

\subsection{CM types and the Frobenius divisor}
\label{v223:general-CM-divisors}
\noindent Put $G=(\Z/M\Z)^\times$, with
$\sigma_u(\zeta_M)=\zeta_M^u$, and fix one rational prime $\ell$
below $\Sigma$. Let $\mathfrak P_\ell$ be the chosen good-reduction
prime of Section~\ref{geo170:universal}, lying over
$\mathfrak q_\ell$ in $F$. A good-reduction field containing $F$
will be denoted by $L$ in this subsection; thus it may be $F$
when $\ell\nmid M$. Set
\[
 q_\ell=\Nm\mathfrak P_\ell=\ell^{f_L},\qquad
 f_L=f(\mathfrak P_\ell/\ell),\quad
 f_F=f(\mathfrak q_\ell/\ell),\quad m_\ell=f_L/f_F.
\]
\noindent Thus $m_\ell=f(\mathfrak P_\ell/\mathfrak q_\ell)$ is common
to all eigenspaces of this Frobenius.

\noindent The function $h_{v,b}$ in \eqref{v222:RM-CM-type} is the
indicator of the CM type of $F_{v,b}$, inflated to $G$.
Indeed the Fermat holomorphic differential with positive indices
$(r,s)$ occurs precisely when $r+s<M$; here $(r,s)=(uv,tub)$
modulo $M$ \cite[pp.~1183--1184]{KoblitzRohrlich}.
The deck maps act on the Jacobian by pushing forward divisors.
Their tangent eigenvalues are the eigenvalues of their pullback
on holomorphic differentials; this fixes the CM type used here.
In particular
\[
                    h_{v,b}(u)+h_{v,b}(-u)=1.
\]
\noindent This fixes the sign of the normalized divisor below.

\begin{proposition}
% Former title: The CM Frobenius divisor
\label{v222:RM-divisor}
\noindent The determinant defined in Section~\ref{v223:uniform-Frobenius}
has divisor
\[
             \divi_F B_{\mathfrak q_\ell}^+(v)=m_\ell H_{v,\ell},
\]
\noindent where $H_{v,\ell}$ is given by \eqref{v228:principalized-ideal}.
More precisely, for each admissible $(v,b)$,
\begin{equation}\label{v222:RM-Frobenius-ideal}
 \divi_FB_{\mathfrak q_\ell}(v,b)
       =m_\ell\sum_{u\in G}
                h_{v,b}(u)[\sigma_u^{-1}\mathfrak q_\ell].
\end{equation}
\noindent Repeated occurrences of a prime on the right are added. This
formula is valid whether or not $\ell$ ramifies in $F$.
\end{proposition}
\begin{proof}
\noindent Fix the $\ell$-adic embedding inducing $\mathfrak P_\ell$ and put
$F_c=F_{v,b}$. For a prime $\mathfrak w\mid\ell$ of $F_c$, let
$\Phi_{\mathfrak w}$ be the CM-type embeddings inducing
$\mathfrak w$. The local Shimura--Taniyama formula is
\begin{equation}\label{v223:local-CM-slope}
 \frac{v_{\mathfrak w}(B_{\mathfrak q_\ell}(v,b))}
      {v_{\mathfrak w}(q_\ell)}
       =\frac{|\Phi_{\mathfrak w}|}
                    {[(F_c)_{\mathfrak w}:\Q_\ell]}.
\end{equation}
\noindent It holds without an unramifiedness hypothesis
\cite[Section~4, Lemma~5, formula~(10); Section~5, Theorem~3]{TateCM1969}.
The numerator and denominator are the dimension and height of the
corresponding CM $\ell$-divisible group. The formula therefore also
applies to an imprimitive CM type with the rank-one $F_c$-action
of Section~\ref{v223:uniform-CM}.

\noindent Here is its precise conversion to divisors in $F$. Write $e_c,f_c$
for the local ramification index and residue degree in $F_c$,
and $e_F,f_F$ for those in $F$. Since
$v_{\mathfrak w}(q_\ell)=e_cf_L$, \eqref{v223:local-CM-slope} gives
\[
 v_{\mathfrak w}(B_{\mathfrak q_\ell}(v,b))
                  =\frac{f_L}{f_c}|\Phi_{\mathfrak w}|.
\]
\noindent Extension to $F$ multiplies this coefficient by $e_F/e_c$.
The embeddings of $F_c$ in the relevant local decomposition
coset each have $e_Ff_F/(e_cf_c)$ lifts to $F$.
Consequently the coefficient of that prime in the right-hand side
of \eqref{v222:RM-Frobenius-ideal} is
\[
 \frac{f_L}{f_F}\frac{e_Ff_F}{e_cf_c}|\Phi_{\mathfrak w}|
      =\frac{e_F}{e_c}\frac{f_L}{f_c}|\Phi_{\mathfrak w}|,
\]
\noindent which is the required extended valuation. Summing over
$b\in I_v^+$ proves the determinant formula. The same $m_\ell$
works for every first index, although the orbit fields may differ.
\end{proof}

\noindent For clarity, if $D_\ell\subset G$ is the decomposition group of
$\mathfrak q_\ell$, the coefficient at
$\sigma_a^{-1}\mathfrak q_\ell$ is
\[
        m_\ell\sum_{u\in aD_\ell}h_{v,b}(u).
\]
\noindent Thus complete splitting is unnecessary: the CM type is summed over
each decomposition coset. At $\ell\nmid M$, with $L=F$, this is
exactly the Jacobi ideal formula of Lemma~\ref{geo203:jacobi-valuations}
in the positive-character convention of \eqref{v222:RM-Jacobi}.
At $\ell\mid M$, the same statement comes from CM reduction,
not from a character of order $M$ in characteristic $\ell$.

\subsection{The normalized divisor functions}
\label{v223:normalized-divisor-functions}
\noindent Fix the positive RM sign and write
\[
 S_+=\{b\in\Z_\Delta^\times:\chi_\Delta(b)=1\},\qquad
 s=|S_+|=\varphi(\Delta)/2.
\]
\noindent The set is nonempty and stable under $b\mapsto-b$ because $K$ is
real quadratic. Introduce the odd periodic Bernoulli function
\[
 B_1(x)=\begin{cases}\{x\}-\tfrac12,&x\notin\Z,\\0,&x\in\Z.
 \end{cases}
\]
\noindent Its value at integers matters: it removes the pairs for which the
Fermat eigenspace is zero. For $v\in\Z_M$ define a function on $G$ by
\begin{equation}\label{v222:RM-carry}
 d_v(u)=2\sum_{b\in S_+}\left(
 B_1(vu/M)+B_1(bu/\Delta)-B_1((v+tb)u/M)\right).
\end{equation}
\noindent Each $d_v$ has rational values, satisfies $d_v(-u)=-d_v(u)$,
and $d_0=0$.

\begin{lemma}
% Former title: Divisors of the normalized RM determinants

\label{v223:normalized-RM-divisor}
\noindent For $1\le v<M$,
\begin{equation}\label{v222:RM-normalized-divisor}
 \divi_F\!\left(\frac{B_{\mathfrak q_\ell}^+(v)}{j(B_{\mathfrak q_\ell}^+(v))}\right)
      =-m_\ell\sum_{u\in G}d_v(u)[\sigma_u^{-1}\mathfrak q_\ell].
\end{equation}
\end{lemma}
\begin{proof}
\noindent For nonintegral $a,b,a+b$, put
$f(a,b)=\{a\}+\{b\}-\{a+b\}\in\{0,1\}$. The CM-type
indicator is $h=1-f$, and
\[
 h(a,b)-h(-a,-b)
  =-2\bigl(B_1(a)+B_1(b)-B_1(a+b)\bigr).
\]
\noindent Subtract the divisor of $j(B_{\mathfrak q_\ell}^+(v))$ from that of $B_{\mathfrak q_\ell}^+(v)$
in Proposition~\ref{v222:RM-divisor}. This gives the displayed formula
for the admissible terms. If $v+tb=0$ in $\Z_M$, then
$B_1(vu/M)+B_1(bu/\Delta)=0$ and the third term is zero.
If $v=0$, the first term is zero and the other two cancel.
Thus including all $b\in S_+$ and all $v\in\Z_M$ introduces no
contribution from an absent eigenspace.
\end{proof}

\noindent It remains to show that the coefficient rows $d_v$ span all odd
functions on $G$. This is a rational spanning assertion, distinct
from the integral basis calculation in Section~\ref{exc170:classnumber}.

\subsection{Fourier transform}
\label{v223:additive-Fourier}
\noindent Fourier transform is taken in the additive index $v\in\Z_M$;
$u$ remains in the multiplicative group $G$. Put
\[
 \widehat d_k(u)=\sum_{v\in\Z_M}d_v(u)\zeta_M^{-kv},
                                  \qquad k\in\Z_M.
\]
\noindent The transform is invertible after extension of scalars to $\C$,
so the families $(d_v)_v$ and $(\widehat d_k)_k$ have the same row
span there. The rows indexed by $k\in G$ will suffice.

\begin{lemma}
% Former title: The Fourier multiplier
\label{v223:RM-multiplier}
\noindent For $k,u\in G$,
\begin{equation}\label{v222:RM-DFT}
 \widehat d_k(u)=i\lambda_+(k)\cot(\pi ku^{-1}/M),\qquad
 \lambda_+(k)=s-\sum_{b\in S_+}\zeta_\Delta^{kb}>0.
\end{equation}
\end{lemma}
\begin{proof}
\noindent Pairing $b$ with $-b$ makes
$\sum_{b\in S_+}B_1(bu/\Delta)=0$. For $k\ne0$ in $\Z_M$,
differentiation of a finite geometric sum gives
\[
 \sum_{v\in\Z_M}B_1(v/M)\zeta_M^{-kv}
    =\frac{1+\zeta_M^k}{2(1-\zeta_M^k)}
    =\frac i2\cot(\pi k/M).
\]
\noindent In the first term of \eqref{v222:RM-carry}, replace $v$ by
$u^{-1}v$. In the last term also translate $v$ by $tb$;
that translation contributes $\zeta_\Delta^{kb}$.
This proves the identity. Its multiplier is real, and
\[
 \lambda_+(k)=\sum_{b\in S_+}
               \bigl(1-\cos(2\pi kb/\Delta)\bigr)>0,
\]
\noindent since $k$ and $b$ are units modulo $\Delta$.
\end{proof}

\noindent The multiplier also displays the same quadratic Gauss sum as the
RM endomorphism in Section~\ref{geo170:curves}. Splitting the
indicator of $S_+$ as $(1+\chi_\Delta)/2$ gives
\begin{equation}\label{v223:Gaussian-multiplier}
 \lambda_+(k)=
 \frac{\varphi(\Delta)-\mu(\Delta)-\chi_\Delta(k)\sqrt\Delta}{2}.
\end{equation}
\noindent Here $\mu(\Delta)$ is the M\"obius function: the ordinary character
sum is the Ramanujan sum $c_\Delta(k)=\mu(\Delta)$ for $k$ a unit,
and the quadratic sum is $\chi_\Delta(k)\sqrt\Delta$
\cite[Chapter~1]{BerndtEvansWilliams}.
For the negative RM sign the last minus sign becomes a plus sign.
The elementary positivity above handles both; the remaining
nonvanishing question concerns the cotangent matrix.

\subsection{Character eigenvalues and nonvanishing}
\label{v223:general-L-nonvanishing}
\noindent Consider the endomorphism
\[
 \mathcal T:\C^G\longrightarrow\C^G,\qquad
 (\mathcal Tf)(k)=\sum_{u\in G}
                  \cot(\pi ku^{-1}/M)f(u).
\]
\noindent Multiplicative characters $\eta:G\to\C^\times$ form a basis of
$\C^G$. Substituting $a=ku^{-1}$ gives the exact eigenvalue formula
\begin{equation}\label{v223:cotangent-eigenvalue}
 \mathcal T\eta=\kappa_\eta\eta,
 \qquad\kappa_\eta=
       \sum_{a\in G}\eta(a)^{-1}\cot(\pi a/M).
\end{equation}
\noindent Even characters have eigenvalue zero, by pairing $a$ and $-a$.
For odd characters, the following classical fact is the analytic
input to the spanning theorem
\cite[Proposition~3, pp.~134--135]{Girstmair2001}.

\begin{lemma}
% Former title: Classical nonvanishing at one
\label{v223:L1-nonvanishing}
\noindent For every odd Dirichlet character $\eta$ modulo $M$,
\begin{equation}\label{v222:cotangent-L}
 \sum_{a\in G}\eta(a)\cot(\pi a/M)
                       =\frac{2M}{\pi}L(1,\eta)\ne0.
\end{equation}
\noindent Thus $\kappa_\eta=(2M/\pi)L(1,\eta^{-1})\ne0$.
No primitivity assumption on $\eta$ is needed.
\end{lemma}
\begin{proof}
\noindent Sum the partial-fraction expansion
\[
 \frac\pi M\cot(\pi a/M)
       =\lim_{R\to\infty}\sum_{h=-R}^{R}\frac1{a+hM}
\]
\noindent against $\eta(a)$, extended by zero at nonunits. Oddness pairs
the positive and negative terms, giving $2L(1,\eta)$.
Dirichlet's nonvanishing theorem applies
\cite[p.~149; Chapter~12]{ApostolANT}.
If $\eta^\circ$ of conductor $f$ induces $\eta$, then
\begin{equation}\label{v223:imprimitive-L1}
 L(1,\eta)=L(1,\eta^\circ)
      \prod_{\substack{r\mid M\\r\nmid f}}
                   (1-\eta^\circ(r)/r).
\end{equation}
\noindent Each extra factor is nonzero, since $|\eta^\circ(r)/r|<1$
\cite[p.~262]{ApostolANT}.
\end{proof}

\begin{lemma}
% Former title: One RM sign spans the odd functions
\label{v222:RM-Fourier}
\noindent The functions $d_v$, $1\le v<M$, span
\[
 \Q[G]^-=\{f:G\to\Q:f(-u)=-f(u)\}.
\]
\end{lemma}
\begin{proof}
\noindent By \eqref{v223:cotangent-eigenvalue} and
Lemma~\ref{v223:L1-nonvanishing}, $\mathcal T$ is invertible on
the odd functions and zero on the even functions. Its rank is
therefore $|G|/2$. The matrix
$(\widehat d_k(u))_{k,u\in G}$ is obtained from the matrix of
$\mathcal T$ by multiplying row $k$ by the nonzero scalar
$i\lambda_+(k)$. It has the same rank $|G|/2$.
These rows are linear combinations of all the rows $d_v$.
Hence the latter have rank at least $|G|/2$ over $\C$;
their oddness gives the reverse inequality. Since all their
entries are rational, their rank over $\Q$ is also $|G|/2$.
Finally $d_0=0$, so that row may be discarded.
\end{proof}

\begin{remark}
% Former title: Relation with the ordinary $L$-value calculation

\label{v223:L0-L1-comparison}
\noindent For the ordinary construction, the divisor eigenvalues are
proportional to $L(0,(\chi_\Delta\psi)^{-1})$, where $\psi$ is
odd modulo $p^a$. The primitive conductor contains every prime
dividing $\Delta p^a$, so no additional Euler factor occurs.
The uniform construction must allow every odd character modulo
$M$; its extra factors $1-\eta^\circ(r)$ at zero may vanish.
For example, the odd quadratic character $\chi_{-8}$ of conductor $8$,
inflated to modulus $24$, has
\[
 L(0,\eta)=0,
 \qquad L(1,\eta)=\frac23L(1,\eta^\circ)\ne0.
\]
\noindent Fourier transform in all first indices replaces this calculation
by the cotangent operator, with nonzero factors at one.
For primitive odd characters the functional equation
\eqref{geo170:L0-equation} relates the two descriptions.
All first indices suffice; each need not be indispensable.
\end{remark}

%\begin{example}
% Former title: A component supplied by a nonprimitive first index

%\label{v223:nonprimitive-example}
%\noindent Take $M=\Delta=24$, so $S_+=\{1,5,19,23\}$.
%An odd function on $G$ is determined by its values at
%$u=1,5,7,11$. Since $d_{-v}=-d_v$, the rows with primitive
%first index are generated by $v=1,5,7,11$. Formula
%\eqref{v222:RM-carry} gives their matrix
%\[
% \begin{pmatrix}
% -3&-1&-3&-1\\
% -1&-3&-1&-3\\
% 0&0&-4&-4\\
% 0&0&-4&-4
% \end{pmatrix}.
%\]
%\noindent It has rank three and annihilates $(1,-1,-1,1)^t$, the values
%of the odd quadratic character $\chi_{-8}$ at these four units.
%The additional row $v=3$ is $(-2,0,4,-2)$. Replacing the last
%row by this one gives determinant $256$, hence full rank four.
%Thus restricting the first indices to units would lose a
%component even on the curve $C_{24,24}$. This is the same
%imprimitive character for which the $L$-value at zero vanished above.
%Moreover $\chi_{-8}$ is trivial on the decomposition group
%$D_3=\{1,11,17,19\}$. It is exactly the odd component of the
%prime-divisor representation above $3$. The missing component
%is therefore relevant to arithmetic support, not merely to the
%larger space of functions on $G$; it is used in
%Section~\ref{v222:bad-prime-example}.
%\end{example}

\subsection{Proof of Theorem~\ref{geo170:all-directions}}
\label{v223:general-span-proof}
\begin{proof}[Proof of Theorem~\ref{geo170:all-directions}]
\noindent For the chosen prime $\mathfrak q_\ell$, define
\[
 \mathcal P_\ell:\Q^G\longrightarrow
                  \bigoplus_{\mathfrak q\mid\ell}\Q[\mathfrak q],
 \qquad f\longmapsto
             \sum_{u\in G}f(u)[\sigma_u^{-1}\mathfrak q_\ell].
\]
\noindent It is surjective and intertwines $u\mapsto-u$ with complex
conjugation. Applying $(1-j)/2$ to a preimage therefore proves
surjectivity on odd parts.

\noindent Lemma~\ref{v222:RM-Fourier} and \eqref{v222:RM-normalized-divisor}
now show that the divisors of
\[
 \frac{B_{\mathfrak q_\ell}^+(v)}{j(B_{\mathfrak q_\ell}^+(v))}\qquad(1\le v<M)
\]
\noindent span all odd divisors over $\ell$ in $F$. The factor $m_\ell$ is a nonzero integer and does not
change their rational span. Equivalently, projection selects
exactly the characters of $G$ trivial on $D_\ell$; this is the
same decomposition-group step as in the ordinary proof.
If $j\in D_\ell$, every prime over $\ell$ is fixed by $j$ and
the odd target is zero. Otherwise its dimension is
$[G:D_\ell]/2$. This includes decomposition groups containing
wild inertia.

\noindent For $\mathfrak p$ in $E$ and $\mathfrak q$ above it in $F$, norm
and extension of divisors are given by
\[
 N_{F/E}[\mathfrak q]
        =f(\mathfrak q/\mathfrak p)[\mathfrak p],\qquad
 \operatorname{ext}_{F/E}[\mathfrak p]
        =\sum_{\mathfrak q\mid\mathfrak p}
                  e(\mathfrak q/\mathfrak p)[\mathfrak q].
\]
\noindent They commute with $j$, and
$N_{F/E}\operatorname{ext}_{F/E}=[F:E]$. Hence norm is surjective
on odd divisors after tensoring with $\Q$. Since
$A_{\mathfrak q_\ell}(v)=N_{F/E}\!\left(\frac{B_{\mathfrak q_\ell}^+(v)}{j(B_{\mathfrak q_\ell}^+(v))}\right)$, it gives
\eqref{v222:full-RM-span}, one rational prime at a time.

\noindent By Lemma~\ref{v203:divisor-kernel}, the divisor map on
$U_\Sigma^{(1)}$ has kernel $\mu(E)$ and image a full-rank
sublattice of $\Lambda_\Sigma^{(1)}$. The spanning just proved
therefore gives finite index in the unit group after adjoining
$\mu(E)$.
\end{proof}

\noindent No class-group hypothesis is needed for this theorem; target
membership was proved in Section~\ref{v206:coverage}.
Section~\ref{exc170:classnumber} removes the possible $p$-primary
index. In particular, the rational projection $(1-j)/2$ used above
is not an operation modulo $2^a$.

\subsection{An alternative using diagonal Fermat quotients}
\label{v222:diagonal-alternative}
\noindent The same divisor components can be obtained using only good
reduction, if the curve is allowed to vary with the conductor.
For $c>2$, let $C_c$ be the smooth projective curve with affine
equation $z^c=x(1-x)$. The map
\[
 \mathcal F_c\longrightarrow C_c,\qquad
                     (U,V)\longmapsto(U^c,UV)
\]
\noindent selects the Fermat character pairs $(u,u)$: both character
arguments are equal. For $u\in\Z_c^\times$, the indices
$u,u,2u$ are nonzero in $\Z_c$. At a prime $\ell\nmid c$ these
curves have good reduction and their corresponding Frobenius
values are the signed Jacobi sums, with the convention of
Lemma~\ref{geo203:Jacobi-Frobenius}.

\noindent Put $F_c=\Q(\zeta_c)$ and $E_c=E\cap F_c$. For each divisor
$c>2$ of $M$ and each rational support prime $\ell\nmid c$,
choose $\mathfrak q\mid\ell$ in $F_c$. With the inverse residue
character $\lambda$ of order $c$ from
Lemma~\ref{geo203:jacobi-valuations}, define
\begin{equation}\label{v228:diagonal-input}
 \begin{gathered}
 J_{c,\mathfrak q}=J(\lambda,\lambda)
     =\sum_{x\in\OO_{F_c}/\mathfrak q}\lambda(x)\lambda(1-x),\\
 A_{\mathfrak q}^{\rm diag}(c)
     =N_{F_c/E_c}\left(\frac{J_{c,\mathfrak q}}
                              {j(J_{c,\mathfrak q})}\right)
        \in E_c^\times\subset E^\times,\\
 \mathcal J_\Sigma=
   \{\sigma(A_{\mathfrak q}^{\rm diag}(c)):c\mid M,\ c>2,
          \ \ell\text{ lies below }\Sigma,\ \ell\nmid c,
          \ \sigma\in\Gal(E/\Q)\}.
 \end{gathered}
\end{equation}
\noindent Here $J$ is exactly the sum defined in \eqref{geo170:jacobi}.
All these values lie in $U_\Sigma^{(1)}$.

\begin{theorem}
% Former title: Spanning by diagonal Jacobi sums

\label{v228:diagonal-span}
\noindent One has
\[
 \operatorname{span}_\Q\{\divi_E A:A\in\mathcal J_\Sigma\}
                     =\Lambda_\Sigma^{(1)}\otimes\Q.
\]
\noindent Thus $\langle\mu(E),\mathcal J_\Sigma\rangle$, after internal
$p$-saturation, also gives $W_\Sigma^{(1)}$. Every Jacobi sum
in this construction is evaluated at a prime of good reduction,
including when $p=2$ or the support contains primes dividing $M$.
\end{theorem}
\begin{proof}
\noindent Fix a rational support prime $\ell$. The prime-divisor
representation in $E$, after extension to $\C$, contains exactly
the characters of $\Gal(E/\Q)$ trivial on its decomposition
group $D_\ell$, each with multiplicity one. Let $\eta$ be such
an odd character, and let $\eta^\circ$ be its primitive Dirichlet
character of conductor $c$. It is trivial on inertia at $\ell$;
therefore
\[
                   \ell\nmid c,\qquad \eta^\circ(\ell)=1.
\]
\noindent Also $c\mid M$ and $c>2$. This is the reason to use $C_c$ at
this prime: its level $c$ is prime to the residue characteristic.

\noindent By Lemma~\ref{geo203:jacobi-valuations}, the divisor of
$J_{c,\mathfrak q}$ is obtained by projecting the function
\[
             f_c(u)=2\{u/c\}-\{2u/c\},\qquad u\in\Z_c^\times,
\]
\noindent to prime divisors. Its primitive odd coefficient is
\begin{equation}\label{v228:diagonal-coefficient}
 \sum_{u\in\Z_c^\times}\eta^\circ(u)f_c(u)
   =-L(0,\eta^\circ)
      \begin{cases}
        2-\eta^\circ(2)^{-1},&c\text{ odd},\\
        2,&c\text{ even}.
      \end{cases}
\end{equation}
\noindent Indeed $\sum_u\eta^\circ(u)\{u/c\}=-L(0,\eta^\circ)$.
If $c$ is odd, substitute $2u$ in the second sum. If $c$ is
even, primitivity implies $4\mid c$ and the two lifts
$u,u+c/2$ from $\Z_{c/2}^\times$ have opposite character
values; their contributions to $\{2u/c\}$ cancel.
The right-hand side is nonzero by
Lemma~\ref{geo203:L0-nonvanishing}. Passing to $J/j(J)$
multiplies an odd coefficient by two.

\noindent Projection to divisors does not kill this coefficient because
$\eta^\circ(\ell)=1$. Moreover the field cut out by $\eta$
lies in both $E$ and $F_c$, hence in $E_c$. Norm to $E_c$
therefore preserves its component, up to a nonzero integer
factor; extension of divisors from $E_c$ to $E$ does likewise.
The conjugates of $A_{\mathfrak q}^{\rm diag}(c)$ consequently span a
space containing the $\eta$-component. Every odd component is
obtained in this way. The equality of rational spans follows,
and the finite-index and saturation assertions follow as above.
\end{proof}

\noindent At $\ell=2$, every contributing conductor $c$ is odd and
$\eta^\circ(2)=1$, so the extra factor in
\eqref{v228:diagonal-coefficient} is exactly one. If
$j\in D_\ell$, there is no odd divisor component at that prime.
Thus no bad-reduction Frobenius evaluation or principalization
is required for this alternative. It uses several curves,
choosing their conductors according to the components that
actually occur. The fixed RM construction instead keeps
$\Delta\mid M$, so primes dividing $\Delta$ cannot be avoided
in its defining level; Section~\ref{v222:principalized-input}
handles them while preserving the quadratic multiplication.

\section{A basis of the general Kummer module}
\label{exc170:classnumber}

\noindent The RM divisor theorem gives a finite-index subgroup of
$U_\Sigma^{(1)}$. We now apply the same $p$-saturation as in
Section~\ref{root170:saturation}; the additional basis class comes
from roots of unity. The procedure is standard
\cite[Section~4.3]{BFHP2022}. The steps parallel
Section~\ref{root170:saturation}: compute the divisor lattice,
extract the missing $p$th roots in $E$, and read off a basis modulo
$n$. We give the corresponding statements, including the torsion
and the criterion for omitting root extraction.

\noindent Let $n=p^a$, let $E/\Q$ be abelian CM with $K,\mu_n\subset E$,
and fix a finite Galois-stable support $\Sigma$.
The module is $W_\Sigma^{(1)}(E,n)$ of
Section~\ref{root170:ordinary-space}. Target membership, with the
signature correction at two, is supplied separately by
Theorem~\ref{exc170:allcoverage}.

\subsection{What the divisor calculation supplies}
\noindent For computation, take the Jacobi values at $\ell\nmid M$ and the
explicit principalized inputs of Section~\ref{v222:principalized-input}
at $\ell\mid M$, and enumerate them as $A_1,\ldots,A_m$.
One may also use actual Frobenius values if available, or the
$C_c$ inputs of Section~\ref{v222:diagonal-alternative}. In every
case these are specified elements of $E$, with computable ideals;
no unknown Frobenius scalar occurs in the following algorithms. Put
\[
 \begin{gathered}
 \Gamma=\langle\mu(E),A_1,\ldots,A_m\rangle\subset U_\Sigma^{(1)},
 \qquad L_\Gamma=\sum_{i=1}^m\Z\,\divi_E(A_i),\\
 L_\Gamma\subseteq L_\Sigma^{(1)}\subseteq\Lambda_\Sigma^{(1)},
 \qquad B=[\Lambda_\Sigma^{(1)}:L_\Gamma],\qquad
 I=[U_\Sigma^{(1)}:\Gamma],\qquad r=\rank\Lambda_\Sigma^{(1)}.
 \end{gathered}
\]
\noindent Theorem~\ref{geo170:all-directions} makes $B$ finite. Choose primes
$\mathfrak p_1,\ldots,\mathfrak p_r$, one from each pair
$\{\mathfrak p,j\mathfrak p\}\subset\Sigma$ with
$j\mathfrak p\ne\mathfrak p$. Then
\[
 \Lambda_\Sigma^{(1)}
   =\bigoplus_{k=1}^r\Z\bigl([\mathfrak p_k]-[j\mathfrak p_k]\bigr),
 \qquad
 D=\bigl(v_{\mathfrak p_k}(A_i)\bigr)_{k,i}\in M_{r,m}(\Z).
\]
\noindent Primes fixed by $j$ contribute no basis divisor. Thus
$D:\Z^m\to\Z^r$ is the matrix of
$(a_i)_i\mapsto\divi_E(\prod_iA_i^{a_i})$ in this fixed basis.
Its image is $L_\Gamma$, and its Smith normal form computes $B$
and products of the $A_i$ whose divisors form a basis of $L_\Gamma$.

\noindent Because the divisor kernel is the full group $\mu(E)$, already
included in $\Gamma$, Lemma~\ref{v203:divisor-kernel} gives
\begin{equation}\label{hi171:twoindices}
 U_\Sigma^{(1)}/\Gamma\simeq L_\Sigma^{(1)}/L_\Gamma,
 \qquad
 B=[\Lambda_\Sigma^{(1)}:L_\Sigma^{(1)}]
                  [U_\Sigma^{(1)}:\Gamma].
\end{equation}
\noindent The two factors have different meanings. The first measures which
formal anti-invariant divisors possess a norm-one generator in $E$.
The second measures which such generators are missing from our
chosen subgroup. Internal root extraction changes only the second
factor. In particular, the condition $\Cl(R_\Sigma)=0$ used for
target membership does not make every supported divisor principal
in $\OO_E$, and does not force the first factor to be one.

\subsection{Roots in the coefficient field}
\begin{lemma}
% Former title: Preservation of the norm-one condition

\label{v223:root-closure}
\noindent If $x\in E^\times$ and $x^q\in U_\Sigma^{(1)}$ for an integer
$q\ge1$, then $x\in U_\Sigma^{(1)}$. Consequently
\[
 \Gamma^{\mathrm{sat}}
   =\{x\in E^\times:x^{p^b}\in\Gamma\text{ for some }b\ge0\}
   \subset U_\Sigma^{(1)}
\]
\noindent is the saturation in the intended norm-one group, including at $p=2$.
\end{lemma}
\begin{proof}
\noindent For $\mathfrak p\notin\Sigma$, the equation
$qv_{\mathfrak p}(x)=0$ gives $v_{\mathfrak p}(x)=0$.
Moreover $(xj(x))^q=1$. At every real embedding of $E^+$,
$xj(x)$ is positive, so it equals one. This proves both assertions.
\end{proof}
\noindent Unlike in the ordinary group, no correction by a root of unity is
needed to restore a second norm equation. On the other hand, a
$p$th root is now determined only up to $\mu_p\subset E$.
These alternatives are harmless because $\mu(E)\subset\Gamma$;
roots of unity must nevertheless be kept in the Kummer basis.

\noindent Here is a finite version of the root procedure. At an intermediate
stage let $\Gamma'$ contain $\mu(E)$ and let $g_1,\ldots,g_r$
lift a basis of its divisor lattice. Let $\xi$ generate the
$p$-primary part of $\mu(E)$. Representatives for
$\Gamma'/\Gamma'^p$ are
\begin{equation}\label{v223:root-candidates}
 \gamma_{\boldsymbol a}=\xi^{a_0}\prod_{i=1}^r g_i^{a_i},
 \qquad 0\le a_i<p\quad(0\le i\le r).
\end{equation}
\noindent In the same fixed divisor basis put
\[
 D_{\Gamma'}=\bigl(v_{\mathfrak p_k}(g_i)\bigr)_{1\le k,i\le r}
       \in M_r(\Z).
\]
\noindent This matrix represents the injection
$\Z^r\to\Lambda_\Sigma^{(1)}$ sending the $i$th standard vector
to $\divi_E(g_i)$. A root in $E$ can exist only if
$D_{\Gamma'}(a_1,\ldots,a_r)^{\mathsf t}\equiv0\pmod p$.
For each remaining candidate, factor $X^p-\gamma_{\boldsymbol a}$
over $E$. A linear factor gives an actual element $x\in E$ with
\[
 x^p=\xi^{a_0}\prod_{i=1}^r g_i^{a_i},\qquad
 \divi_E(x)=\frac1p\sum_{i=1}^r a_i\divi_E(g_i).
\]
\noindent Accept it only if
$\divi_E(x)\notin\divi_E(\Gamma')$, and replace $\Gamma'$ by
$\langle\Gamma',x\rangle$. Update the lattice basis by integer
row or column operations, keeping the corresponding products of
actual field elements. Repeat until there is no accepted root.
All these are effective number-field and lattice operations
\cite[Chapters~3--6]{CohenANT}.

\begin{proposition}
% Former title: Termination
\label{exc170:saturation}
\noindent For $I=[U_\Sigma^{(1)}:\Gamma]$, one has
\[
 \Gamma^{\mathrm{sat}}/\Gamma
        =(U_\Sigma^{(1)}/\Gamma)[p^\infty],\qquad
 [\Gamma^{\mathrm{sat}}:\Gamma]=p^{v_p(I)},\qquad
 p\nmid[U_\Sigma^{(1)}:\Gamma^{\mathrm{sat}}].
\]
\noindent The procedure reaches $\Gamma^{\mathrm{sat}}$ after exactly
$v_p(I)\le v_p(B)$ accepted roots. The class map induces
\begin{equation}\label{hi171:modn}
 \Gamma^{\mathrm{sat}}/(\Gamma^{\mathrm{sat}})^n
       \xrightarrow{\sim}
 U_\Sigma^{(1)}/(U_\Sigma^{(1)})^n
       \xrightarrow{\sim}W_\Sigma^{(1)}.
\end{equation}
\end{proposition}
\begin{proof}
\noindent Lemma~\ref{v223:root-closure} puts every tested root in
$U_\Sigma^{(1)}$. Since $\mu(E)\subset\Gamma'$, membership in
$\Gamma'$ is equivalent to membership of its divisor in
$L_{\Gamma'}$, and one root per candidate suffices.
The representatives~\eqref{v223:root-candidates} therefore detect
every index-$p$ enlargement. Apply
Proposition~\ref{v212:elementary-saturation},
\eqref{hi171:twoindices} and Proposition~\ref{root170:U-W}.
\end{proof}

\begin{corollary}
% Former title: The general Kummer basis
\label{v226:general-basis}
\noindent Let $z_1,\ldots,z_r$ lift a $\Z$-basis of
$\Gamma^{\mathrm{sat}}/\mu(E)$. If $p^t$ is the order of the
$p$-primary subgroup of $\mu(E)$, choose its generator $z_0$;
then $t\ge a$ and $[z_0]$ generates $\mu(E)/\mu(E)^n\simeq\Z_n$.
The general module has rank $r+1$, with basis
\begin{equation}\label{root170:general-W}
 (\Z_n)^{r+1}\xrightarrow{\sim}W_\Sigma^{(1)},\qquad
 (e_0,\ldots,e_r)\longmapsto
                 \left[z_0^{e_0}\prod_{j=1}^r z_j^{e_j}\right].
\end{equation}
\noindent For $j\ge1$, recording the root equations and integer basis
changes gives
\begin{equation}\label{v215:general-Jacobi-coordinates}
 z_j^{p^{b_j}}=\xi_j\prod_i A_i^{m_{ji}},\qquad
 b_j\ge0,\quad\xi_j\in\mu(E),\quad m_{ji}\in\Z.
\end{equation}
\end{corollary}
\begin{proof}
\noindent The free part contributes $r$ copies of $\Z_n$ and the cyclic
torsion, which contains $\mu_n$, contributes one. Each adjoining
relation is a $p$th-root equation; substitution gives the stated
identities. The basis consists of the specified classes $[z_j]$,
not merely their displayed powers.
\end{proof}

\noindent The full $p$-primary torsion is essential: $\zeta_3$ is a cube
in $\Q(\zeta_9)$, whereas $\zeta_9$ generates the torsion quotient
modulo cubes. This accounts for the extra class $[z_0]$; ordinary
torsion has order prime to $p$ and contributes no class.

\begin{corollary}
% Former title: When general saturation is unnecessary

\label{v227:general-no-saturation}
\noindent With $\Gamma,I,B$ as above, the following are equivalent:
\[
 \Gamma^{\mathrm{sat}}=\Gamma
 \quad\Longleftrightarrow\quad p\nmid I
 \quad\Longleftrightarrow\quad
 W_\Sigma^{(1)}=\langle[z_0],[A_1],\ldots,[A_m]\rangle.
\]
\noindent In particular $p\nmid B$ suffices. In this case, choose products
\[
                 g_j=\prod_i A_i^{m_{ji}}\quad(1\le j\le r)
\]
\noindent whose divisors form a $\Z$-basis of $L_\Gamma$. Then
$[z_0],[g_1],\ldots,[g_r]$ are a Kummer basis, with no root
extraction and with $b_j=0$, $\xi_j=1$ in
\eqref{v215:general-Jacobi-coordinates}.
\end{corollary}
\begin{proof}
\noindent Apply Proposition~\ref{v212:elementary-saturation} to
$U_\Sigma^{(1)}$, using $I\mid B$ and
$\mu(E)\subset\Gamma$. A finite quotient has no $p$-torsion
exactly when reduction modulo $p^a$ is already surjective.
The divisor basis lifts to a basis of $\Gamma/\mu(E)$.
\end{proof}
\noindent Thus the principalized RM input with $p=3$ and $B=16$ in
Section~\ref{v222:bad-prime-example} needs no internal cube root.
One may also certify the criterion by an invertible residue matrix
on $r+1$ classes selected from $[z_0],[A_i]$, exactly as in
Corollary~\ref{v224:raw-residue-certificate}. By contrast, the
quartic example in Section~\ref{v221:dyadic-example} requires one
internal square root. As in the ordinary case, a particular target
can already lie in the image of $\Gamma$ even when the entire
module does not.

\begin{remark}
% Former title: Divisibility of a divisor is only a root test

\label{v226:lattice-root-remark}
\noindent The congruence on $D_{\Gamma'}$ is necessary for a root in $E$, but does
not prove its existence. The divided divisor must be principal
and its generator must have the required power.
In~\eqref{hi171:twoindices}, saturation changes
$[U_\Sigma^{(1)}:\Gamma]$, not
$[\Lambda_\Sigma^{(1)}:L_\Sigma^{(1)}]$.
Thus the stopping condition is the absence of further roots in
$E$, not an index prime to $p$ in the formal divisor lattice.
\end{remark}

\section{Selection by reciprocity in the general case}
\label{v222:general-selection-section}
\label{v217:general-radicand-construction}

\noindent Selection by Artin reciprocity uses the same Kummer pairing and
separating-prime argument as Section~\ref{root170:selection}.
The general Kummer module includes a torsion basis class.
At two we first select the Kummer class of the signature-corrected
character $\psi$, then restore the signature factor.
The underlying Artin map method is standard \cite[Section~3]{FiekerArtin},
\cite[Sections~5.4.3--5.4.5]{CohenAdvanced},
\cite[Section~4.2]{FHS}; here the input basis is supplied by the RM
inputs and $p$-power roots taken inside $E$.

\noindent Fix $\chi:G_{\mathfrak f}^+\to\Z_n$ of exact order $n=p^a$.
Use $\psi,\mathfrak g,E,\Sigma$ from
Theorem~\ref{exc170:allcoverage}. Put
$r=\rank\Lambda_\Sigma^{(1)}$ and use the basis
\[
 W_\Sigma^{(1)}=\bigoplus_{j=0}^r\Z_n[z_j],\qquad
             \nu^{-1}(\psi|_{G_E})\in W_\Sigma^{(1)}
\]
\noindent of \eqref{root170:general-W}. The class of $z_0$ generates
$\mu(E)/\mu(E)^n$: one may choose a generator of $\mu(E)$ or of
its full $p$-primary subgroup. A generator of $\mu_n$ need not suffice.
Exclude
\[
 \mathcal T=\Sigma\cup
       \{\mathfrak p:\mathfrak p\mid n\mathfrak g\OO_E\}.
\]
\noindent Call the primes outside $\mathcal T$ admissible. Applying
Lemma~\ref{root170:test-group} with $d=r+1$, enumerate these primes
and choose $\mathfrak q_1,\ldots,\mathfrak q_{r+1}$ whose residue
rows are independent modulo $p$. The resulting isomorphism is
\[
 W_\Sigma^{(1)}\longrightarrow(\Z_n)^{r+1},\qquad
 [x]\longmapsto(r_{\mathfrak q_k}(x))_{k=1}^{r+1}.
\]
\noindent Its matrix in the ordered basis $[z_0],\ldots,[z_r]$ and the target
vector are
\begin{equation}\label{root170:general-matrix}
 \begin{gathered}
 C:(\Z_n)^{r+1}\xrightarrow{\sim}(\Z_n)^{r+1},\qquad
 C_{k,j+1}=r_{\mathfrak q_k}(z_j)
       \quad\substack{1\le k\le r+1\\0\le j\le r},\\
 b_k=\psi([\Nm_{E/K}\mathfrak q_k]_{\mathfrak g}),
 \qquad\boxed{\;\boldsymbol e=C^{-1}\boldsymbol b,\quad
                 \boldsymbol e=(e_0,\ldots,e_r)^{\mathsf t}.\;}
 \end{gathered}
\end{equation}
\noindent For odd $p$, this right-hand side uses $\psi=\chi$. For $p=2$,
with $a_0$ from \eqref{exc170:signcorrection}, one computes
\[
 b_k=\chi([\Nm_{E/K}\mathfrak q_k]_{\mathfrak g})
                  -r_{\mathfrak q_k}(a_0^{n/2})
                         \quad\text{in }\Z_n.
\]
\noindent Here $\chi$ is inflated to $\mathfrak g$; the excluded primes
include those in the signature correction. Thus the same finite-field
residue operation handles both the basis and the correction.

\begin{theorem}
% Former title: Selection by Artin reciprocity in the general case

\label{v216:general-radicand}
\noindent For the unique solution of \eqref{root170:general-matrix}, choose
$0\le e_j<n$ and put
\[
                    \beta_\psi=z_0^{e_0}\prod_{j=1}^r z_j^{e_j}.
\]
\noindent Then $\nu_{\beta_\psi}=\psi|_{G_E}$, and
\begin{equation}\label{v222:general-beta}
 \beta_\chi=
 \begin{cases}
   \beta_\psi,&p\ne2,\\
   a_0^{n/2}\beta_\psi,&p=2
 \end{cases}
 \qquad\text{satisfies}\qquad
                  \nu_{\beta_\chi}=\chi|_{G_E}.
\end{equation}
\end{theorem}
\begin{proof}
\noindent Theorem~\ref{exc170:allcoverage} places
$\nu^{-1}(\psi|_{G_E})$ in $W_\Sigma^{(1)}$.
Evaluation at the selected primes is injective
on this entire module because $C$ is invertible. Its unique vector
is therefore $C^{-1}\boldsymbol b$. At two, $a_0^{n/2}$ represents
$(n/2)\theta$, so multiplication by it restores $\chi$.
\end{proof}

\noindent The answer is the displayed product of the specified $z_j$.
The recorded identities \eqref{v215:general-Jacobi-coordinates}
separately express a power of it in the RM inputs. For
$s_0\ge\max_jb_j$, with
$s_0=0$ if $r=0$,
\begin{equation}\label{v215:general-jacobi-radicand}
 \beta_\psi^{p^{s_0}}
 =z_0^{e_0p^{s_0}}\prod_{j=1}^r\xi_j^{e_jp^{s_0-b_j}}
       \prod_i A_i^{\,\sum_{j=1}^r e_jm_{ji}p^{s_0-b_j}}.
\end{equation}
\noindent The power identity alone can leave a root-of-unity ambiguity.
We do not choose a root of that equation: $\beta_\psi$ is already
defined by the product of the computed elements $z_j$ above.
At good reduction the $A_i$ are
Jacobi products. At bad reduction they are CM Frobenius values,
or the principalized replacements of \eqref{v222:principalized-RM}.
For fixed basis and admissible testing primes, one inverse matrix
serves every eligible target. No ray class group of $E$ is needed
for this selection.

\noindent The proof of Corollary~\ref{v224:raw-residue-certificate} applies here with rank
$r+1$: an invertible residue matrix on $z_0$ and $r$ raw RM values
certifies that they already form the required Kummer basis. Omitting
the torsion coordinate would invalidate that rank test.

\noindent The examples in Sections~\ref{v221:dyadic-example} and
\ref{v206:five-example} give the selected Kummer coordinates, including the torsion
coordinate; the second also displays the signature correction. These
calculations determine only the restriction to $E$. Its order may
be smaller than $n$, even equal to one. In particular, the full
algebra $E[X]/(X^n-\beta_\chi)$ must be kept for the equivariant
descent of the next section, which recovers the original character
and its field over $K$.

\section{General generator by equivariant descent}
\label{root170:general-algebra}

\noindent Let $\chi:G_{\mathfrak f}^+\to\Z_n$ have exact order $n=p^a$,
and put $\beta=\beta_\chi$ from \eqref{v222:general-beta}.
The preceding section gives $\nu_\beta=\chi|_{G_E}$. The remaining
step is standard equivariant Galois descent
\cite[Section~2.3]{GilleSzamuely}\cite[Section~4.3]{FHS}:
descend the cyclic torsor, select its Artin character, and take an
invariant generator. Work in the full algebra
\[
 B_\beta=E[X]/(X^n-\beta),\qquad
 \alpha=X\bmod(X^n-\beta),\qquad
 r:B_\beta\longrightarrow B_\beta,\quad
 r|_E=\id,\quad r(\alpha)=\zeta_n\alpha.
\]
\noindent Its spectrum is a torsor under the constant group $\Z_n$: $r$ acts
freely and transitively on its $n$ geometric points. The algebra
need not be a field, since restriction to $E$ can lower the order
of $\chi$.

\subsection{The descent equations and the resulting generator}
\noindent For $\delta\in\Gal(E/K)$ choose $1\le u_\delta<n$ with
$\delta(\zeta_n)=\zeta_n^{u_\delta}$; take $u_1=1$.
Equivariant descent is a semilinear action of $\Gal(E/K)$ on
$B_\beta$ commuting with $r$. Equivariance forces the form
\[
 d_\delta:B_\beta\longrightarrow B_\beta,\qquad
 d_\delta|_E=\delta,\quad
 d_\delta(\alpha)=c_\delta\alpha^{u_\delta}.
\]
\noindent Indeed $d_\delta(\alpha)$ belongs to the $\zeta_n^{u_\delta}$-eigenspace
of $r$. Preserving the relation $\alpha^n=\beta$ and composing the maps gives
\begin{equation}\label{root170:marks}
 \begin{aligned}
 c_1&=1,&
 c_\delta^n&=\delta(\beta)/\beta^{u_\delta},\\
 c_{\delta\delta'}&=
 \delta(c_{\delta'})c_\delta^{u_{\delta'}}
 \beta^{(u_\delta u_{\delta'}-u_{\delta\delta'})/n}.
 \end{aligned}
\end{equation}
\noindent The last power reduces $\alpha^{u_\delta u_{\delta'}}$ using
$\alpha^n=\beta$. Distinct $\delta$ may have the same $u_\delta$,
so the coefficients are indexed by automorphisms of $E/K$.
For a solution $c=(c_\delta)_\delta$, define its fixed algebra
\[
 A_{\beta,c}=\{x\in B_\beta:d_\delta(x)=x
                              \text{ for every }\delta\in\Gal(E/K)\}.
\]
\noindent Its cyclic action determines a character
$\chi_c:G_K\to\Z_n$, as proved in
Proposition~\ref{root170:marked-descent}. Fix a $K$-basis
$b_1,\ldots,b_d$ of $E$, where $d=[E:K]$.

\begin{theorem}
% Former title: The general cyclic generator

\label{root170:general-R}\label{v203:reynolds-generators}
\noindent Exactly one solution of \eqref{root170:marks} has character
$\chi_c=\chi$. Its fixed algebra is isomorphic to $H_\chi$. It is spanned over
$K$ by
\begin{equation}\label{root170:reynolds}
 \theta_{ik}=\sum_{\delta\in\Gal(E/K)}
       \delta(b_i)c_\delta^k\alpha^{ku_\delta},
 \qquad 1\le i\le d,\quad0\le k<n.
\end{equation}
\noindent One can choose integers $\lambda_{ik}$ such that
\begin{equation}\label{v222:general-primitive}
 \boxed{\quad H_\chi=K(\theta_\chi),\qquad
        \theta_\chi=\sum_{i,k}\lambda_{ik}\theta_{ik}.\quad}
\end{equation}
\noindent A proposed combination is primitive precisely when its $n$
translates $r^a(\theta_\chi)$, $a\in\Z_n$, are distinct.
\end{theorem}
\noindent The next two subsections compute this descent and select $\chi$.
The proof at the end constructs the primitive invariant.

\subsection{Computing the equivariant descents}
\label{v218:descent-computation}
\noindent For each nonidentity $\delta$, factor
$X^n-\delta(\beta)/\beta^{u_\delta}$ over $E$. Its linear factors
give the possible $c_\delta$. Set $c_1=1$ and keep exactly
the families satisfying the composition equations in
\eqref{root170:marks}. This is a finite search with at most
$n^{d-1}$ candidates, using only the known field $E$ and
element $\beta\in E^\times$. The following standard descent statement
proves that the list is complete and nonempty.

\begin{proposition}
% Former title: Cyclic torsors and their characters

\label{root170:marked-descent}
\noindent Solutions of \eqref{root170:marks} are precisely the equivariant
semilinear $\Gal(E/K)$-actions on $B_\beta$. They correspond
bijectively to continuous characters
$\chi_c:G_K\to\Z_n$ extending $\nu_\beta$. If
$m=|\im\chi_c|$ and $L=\overline\Q^{\ker\chi_c}$, then
\[
 E\otimes_K A_{\beta,c}\simeq B_\beta,\qquad
 A_{\beta,c}\simeq L^{\,n/m},\qquad [L:K]=m.
\]
\noindent In particular, the fixed algebra is a field exactly when the
character is surjective.
\end{proposition}
\begin{proof}
\noindent Semilinear descent gives the first isomorphism
\cite[Lemma~2.3.8 and Remark~2.3.10]{GilleSzamuely}.
The finite \'etale Galois correspondence
\cite[Theorems~8.21--8.23]{MilneFT} identifies equivariant
descents with actions on the $n$ points commuting with $r$.
Its centralizer is $\Z_n$, giving exactly the characters extending
$\nu_\beta$. Their orbits have stabilizer $\ker\chi_c$, hence
$A_{\beta,c}\simeq L^{n/m}$. This gives actual coefficient families:
two equivariant identifications over $E$ differ by a power of $r$,
which commutes with every descent and therefore leaves it unchanged.
\end{proof}

\begin{lemma}
% Former title: Twists and rescaling the generator

\label{v203:descent-ambiguity}\label{v206:integral-monomial-descent}
\noindent Relative to one solution, all others are indexed by
$b\in\Hom(\Gal(E/K),\Z_n)$:
\[
 \begin{gathered}
 d'_\delta=r^{b(\delta)}d_\delta,\qquad
 c'_\delta=\zeta_n^{u_\delta b(\delta)}c_\delta,\\
 \chi_{c'}(g)=\chi_c(g)-b(g|_E).
 \end{gathered}
\]
\noindent Changing $\alpha$ to $a\alpha$, for $a\in E^\times$, preserves
the character and gives
\begin{equation}\label{root170:gauge}
 \beta\longmapsto a^n\beta,\qquad
 c_\delta\longmapsto\frac{\delta(a)}{a^{u_\delta}}c_\delta.
\end{equation}
\noindent If $(\beta_1,c^{(1)})$ and $(\beta_2,c^{(2)})$ correspond to
$\chi_1,\chi_2:G_K\to\Z_n$, then
\[
 \left(\beta_1\beta_2,
       (c^{(1)}_\delta c^{(2)}_\delta)_\delta\right)
                     \quad\text{corresponds to}\quad\chi_1+\chi_2.
\]
\noindent In particular the signature character $(n/2)\theta$
at two has data
\begin{equation}\label{exc170:signpair}
 \beta_{a_0}=a_0^{n/2},\qquad
 c_{\delta,a_0}=a_0^{(1-u_\delta)/2}.
\end{equation}
\end{lemma}
\begin{proof}
\noindent Two extensions of $\nu_\beta$ differ by a homomorphism from
$\Gal(E/K)$. Replacing the lift $d_\delta$ by
$r^{b(\delta)}d_\delta$ subtracts $b(\delta)$ from the
translation character. The displayed coefficients follow by
applying the maps to $\alpha$. The remaining statements follow
by substitution; \eqref{exc170:signpair} uses the root
$\sqrt{a_0}$ of $X^n-a_0^{n/2}$.
\end{proof}
\noindent After finding one solution, its finite group of twists supplies
the whole list without further root searches.

\begin{lemma}
% Former title: The field factors of the cyclic torsor

\label{v208:why-whole-algebra}
\noindent Put $e=[H_\chi\cap E:K]$. Then
\begin{equation}\label{v203:kummer-components}
 B_\beta\simeq E\otimes_KH_\chi\simeq(EH_\chi)^e.
\end{equation}
\noindent If $e>1$, the selected descent permutes these $e$ field factors
transitively.
\end{lemma}
\begin{proof}
\noindent The tensor-product decomposition for two Galois extensions gives
\eqref{v203:kummer-components}. A union of factor orbits determines
an idempotent of the fixed algebra. That algebra is the field
$H_\chi$, so there is only one orbit.
\end{proof}
\noindent Thus a single radical field describes only the restricted character;
the cyclic torsor supplies all $n$ points needed for descent.

\begin{example}
% Former title: A quadratic character trivial over $E$

\label{v205:split-torsor-example}
\noindent Let $K=\Q(\sqrt3)$, $E=K(i)$, and let $\chi$ cut out $E/K$.
Its restriction is zero, so $\beta=1$. On
$E[\alpha]/(\alpha^2-1)$ the two descents
\[
 d_\pm|_E=j,\qquad d_\pm(\alpha)=\pm\alpha
\]
\noindent have fixed algebras $K\times K$ and
$K[i\alpha]\simeq E$, respectively.
Here $\alpha+d_-(\alpha)=0$, whereas
$i\alpha+d_-(i\alpha)=2i\alpha$ is primitive; the coefficients
$b_i$ in \eqref{root170:reynolds} are therefore necessary.
\end{example}

\subsection{Selecting the character by its Artin values}
\label{root170:general-selection}
\noindent For a candidate $c$, compute $A_{\beta,c}$ by the linear
equations $d_\delta(x)=x$ in the $K$-basis $b_i\alpha^k$.
Decompose it into fields and choose one factor $L$; all factors
are isomorphic by Proposition~\ref{root170:marked-descent}.
The next two standard calculations determine whether its
character equals $\chi$.

\begin{lemma}
% Former title: Conductor from relative discriminants

\label{exc170:conductor}\label{exc170:directconductor}
\noindent If $[L:K]=p^b>1$ and $L_{b-1}$ is its subfield of degree
$p^{b-1}$, then
\begin{equation}\label{exc170:discconductor}
 f_{\mathfrak q}(L/K)=
 \frac{v_{\mathfrak q}(\mathfrak d_{L/K})-
       v_{\mathfrak q}(\mathfrak d_{L_{b-1}/K})}
      {p^{b-1}(p-1)}.
\end{equation}
\noindent For $L=K$, put $f_{\mathfrak q}=0$. In either case
\begin{equation}\label{exc170:conductoraccept}
 L\subseteq K_{\mathfrak f}^+
 \quad\Longleftrightarrow\quad
 f_{\mathfrak q}(L/K)\le v_{\mathfrak q}(\mathfrak f)
 \quad\text{at every finite }\mathfrak q.
\end{equation}
\end{lemma}
\begin{proof}
\noindent Apply the conductor--discriminant formula
\cite[Chapter~VI, Section~3]{SerreLocal}: subtraction cancels the
characters factoring through $L_{b-1}$, leaving
$p^{b-1}(p-1)$ faithful characters of the same conductor.
Both real places are already allowed in the strict modulus.
\end{proof}
\noindent Use maximal-order discriminants, computable together with the
subfields by \cite[Chapters~3--6]{CohenANT}. Only primes dividing
$\mathfrak d_{L/K}$ need a conductor test.

\begin{lemma}
% Former title: Artin values from the descent coefficients

\label{v216:descended-artin}
\noindent Let $\mathfrak q$ be a prime of $K$ away from $n$, the
discriminant of $E/K$, and the numerators and denominators of
$\beta,c_\delta$. Put $q=\Nm\mathfrak q$,
$\delta=\Frob_{\mathfrak q}(E/K)$, and choose
$\mathfrak p\mid\mathfrak q$ in $E$. Then
\begin{equation}\label{root170:descended-residue}
 \zeta_n^{\chi_c(\mathfrak q)}
 \equiv\left(\frac{\beta^{(q-u_\delta)/n}}{c_\delta}\right)^{
                    u_\delta^{-1}}\pmod{\mathfrak p},
\end{equation}
\noindent where $u_\delta^{-1}$ is an inverse modulo $n$.
\end{lemma}
\begin{proof}
\noindent On a geometric root, $r^a d_\delta$ acts by
$\alpha\mapsto\zeta_n^{au_\delta}c_\delta\alpha^{u_\delta}$,
whereas arithmetic Frobenius reduces to $\alpha^q$.
Since $q\equiv u_\delta\pmod n$, comparison gives the formula.
The expression inside parentheses has $n$th power
$\beta^q/\delta(\beta)\equiv1$, and reduction is injective on
$\mu_n$. This is the usual Kummer Frobenius calculation
\cite[Proposition~4]{FHS}, applied in each geometric factor.
\end{proof}

\begin{proposition}
% Former title: Finite selection of the descent

\label{v216:select-descent}
\noindent A candidate is the required descent precisely when its field
factor satisfies \eqref{exc170:conductoraccept} and
\[
 \chi_c([\mathfrak a_\ell])=\chi([\mathfrak a_\ell])
 \qquad\text{for generators }[\mathfrak a_\ell]
                         \text{ of }G_{\mathfrak f}^+.
\]
\noindent Exactly one candidate passes these tests, and its fixed algebra
is isomorphic to $H_\chi$.
\end{proposition}
\begin{proof}
\noindent The conductor test makes $\chi_c$ a character of $G_{\mathfrak f}^+$;
equality on generators gives $\chi_c=\chi$. Existence and uniqueness
follow from Proposition~\ref{root170:marked-descent}, and the fixed
algebra is a field because $\chi$ has order $n$.
\end{proof}
\noindent Represent the ray class generators by ideals outside the finite
excluded set, apply \eqref{root170:descended-residue} to their
prime factors, and add the resulting values.
The procedure uses only explicitly constructed fields, not an
equation for the unknown $H_\chi$.
\label{v208:exact-conductor}
Exact finite conductor $\mathfrak f$ can additionally be checked
by testing that $\chi$ does not factor through any
$G_{\mathfrak f/\mathfrak q}^+$ with $\mathfrak q\mid\mathfrak f$.

\subsection{A primitive invariant and comparison with the ordinary case}
\begin{proof}[Proof of Theorem~\ref{root170:general-R}]
\noindent Proposition~\ref{v216:select-descent} gives the unique family.
Averaging the $K$-basis $b_i\alpha^k$ projects onto its fixed
algebra, so the sums \eqref{root170:reynolds} span $H_\chi$.
Primitivity is equivalent to having $n$ distinct translates under $r$.
For an explicit finite search, enumerate these sums as $v_\ell$ and
put $\theta(m)=\sum_{\ell=0}^{dn-1}m^\ell v_\ell$. Each polynomial
$r^a\theta(m)-\theta(m)$, $0<a<n$, is nonzero of degree at most
$dn-1$. Some $0\le m\le(n-1)(dn-1)$ therefore gives a primitive
invariant, verified inside $B_\beta$.
\end{proof}

\label{cmp202:chapter}
\noindent In the ordinary anti-invariant case, the same descent has
$H_\chi\cap E=K$ and is selected by the norm condition
$c_u\tau(c_u)=1$. Theorem~\ref{root170:ordinary-descent-theorem}
then gives the primitive unweighted trace, so enumeration and
the conductor--Artin tests are unnecessary. The distinction is
in how the equivariant descent is selected; the cyclic-torsor
interpretation is the same in both cases.

\section{Construction of the strict ray class field}
\label{mx171:chapter}\label{exc170:mixed}\label{root170:algorithm}
\label{v211:galois-action}

\noindent The cyclic generators for prime-power-order characters of the
\emph{full} group $G_{\mathfrak f}^+$ generate the strict ray class field.
Here the prime power is the order of the character, not its modulus;
its conductor may involve several prime ideals.
We first compute these characters, then assemble their cyclic fields
and describe the Artin action. Finally we give an exact criterion
for the compositum attached to separate prime-power moduli to suffice.

\subsection{Computing the characters of the full modulus}
\label{v227:full-modulus-characters}
\label{mx171:fermatselection}\label{mx171:raysequence}\label{mx171:presentation}
\noindent This subsection supplies the separating characters,
including their values on the ideals used in the reciprocity tests.
The computation is the standard ray-group presentation
\cite[Section~4.3.1]{CohenAdvanced}\cite[Section~3]{FHS}.

\noindent Write $\mathfrak f=\prod_i\mathfrak f_i$ with
$\mathfrak f_i=\mathfrak p_i^{a_i}$ for distinct primes $\mathfrak p_i$,
put $R_i=(\OO_K/\mathfrak f_i)^\times$, and choose ideals
$\mathfrak a_1,\ldots,\mathfrak a_h$ representing an invariant-factor
basis $\Cl^+(K)=\bigoplus_{j=1}^h\Z_{m_j}$.
Choose them prime to $\mathfrak f$,
with
\begin{equation}\label{mx171:classrelation}
          \mathfrak a_j^{m_j}=(\gamma_j),\qquad\gamma_j\gg0.
\end{equation}
\noindent Put $\overline R_{\mathfrak f}=(\prod_iR_i)/\langle\bar\epsilon\rangle$.
Let $x_1,\ldots,x_h$ denote the standard generators of a free abelian group
$\Z^h$, mapping to the ray classes $[\mathfrak a_j]_{\mathfrak f}$.
Writing finite abelian groups additively for the following presentation,
the ray sequence~\eqref{v214:ray-exact} gives
\begin{equation}\label{mx171:grouppresentation}
 G_{\mathfrak f}^+\simeq
 \frac{\overline R_{\mathfrak f}\oplus\bigoplus_j\Z x_j}
      {\langle m_jx_j-[\bar\gamma_j]\rangle}.
\end{equation}
\noindent The class of a residue with totally positive lift $b$ maps to
$[(b)]_{\mathfrak f}$. Thus the relation $m_jx_j=[\bar\gamma_j]$
records exactly \eqref{mx171:classrelation}.
\label{mx171:principalmap}\label{mx171:presentationtheorem}

\noindent Equivalently, a character to $\Z_n$ is specified by additive local
characters $\ell_i:R_i\to\Z_n$ and values
$y_j=\chi([\mathfrak a_j])$ satisfying
\begin{equation}\label{mx171:unitcondition}
 \sum_i\ell_i(\epsilon)=0,\qquad
                    m_jy_j=\sum_i\ell_i(\gamma_j)\quad\text{for every }j.
\end{equation}
\noindent Evaluation means reduction modulo $\mathfrak p_i^{a_i}$; denominators
prime to that ideal are allowed. The first relation is the global
unit relation; the others are the chosen ideal-class relations.
If $\mathfrak a\prod_j\mathfrak a_j^{-n_j}=(b)$ with $b\gg0$,
then
\begin{equation}\label{mx171:evaluation}
 \chi([\mathfrak a])=\sum_i\ell_i(b)+\sum_j n_jy_j.
\end{equation}
\noindent The relations make this value independent of the choices.

\noindent Here is the coordinate map supplied by Smith normal form.
Choose an invariant-factor basis $g_1,\ldots,g_t$ of $\prod_iR_i$,
with respective orders $c_1,\ldots,c_t$. Its residue classes and
the ideals $\mathfrak a_j$ define a surjection
\[
 \pi:\Z^{t+h}\longrightarrow G_{\mathfrak f}^+.
\]
\noindent For each $g_l$, the image is the principal ray class of a totally
positive lift; the last $h$ standard vectors map to
$[\mathfrak a_j]_{\mathfrak f}$. Write
\[
 \bar\epsilon=\prod_{l=1}^t g_l^{u_l},\qquad
 \bar\gamma_j=\prod_{l=1}^t g_l^{v_{lj}},\qquad
 \boldsymbol u=(u_l)_l,\quad V=(v_{lj})_{l,j}.
\]
\noindent For $N=t+h$, the relation matrix is
\[
 A=\begin{pmatrix}
       \operatorname{diag}(c_1,\ldots,c_t)&\boldsymbol u&-V\\
       0&0&\operatorname{diag}(m_1,\ldots,m_h)
    \end{pmatrix}\in M_{N,N+1}(\Z),
       \qquad\operatorname{im}A=\ker\pi.
\]
\noindent Smith normal form gives
\[
 P\in\operatorname{GL}_N(\Z),\quad
 Q\in\operatorname{GL}_{N+1}(\Z),\quad
 PAQ=(\operatorname{diag}(d_1,\ldots,d_N)\mid0),
 \qquad d_j>0.
\]
\noindent The group is finite, so there are $N$ nonzero diagonal entries.
The factors $d_j=1$ may be omitted.

\begin{proposition}
% Former title: A separating family of ray characters

\label{v221:computed-characters}
\noindent With this notation the explicit isomorphism is
\[
 \varphi:G_{\mathfrak f}^+\xrightarrow{\sim}\bigoplus_{j=1}^N\Z_{d_j},
 \qquad\varphi(\pi(v))=((Pv)_j\bmod d_j)_j.
\]
\noindent For each $p\mid d_j$, its coordinate character
\[
 \chi_{j,p}:G_{\mathfrak f}^+\longrightarrow\Z_{p^{v_p(d_j)}},
 \qquad g\longmapsto\varphi_j(g)\bmod p^{v_p(d_j)}
\]
\noindent has exact order $p^{v_p(d_j)}$, and
$\bigcap_{j,p}\ker\chi_{j,p}=\{1\}$.
In particular all orders $p^{v_p(d_j)}$ and all values on the
chosen generators are known before constructing any Kummer module.
Equation~\eqref{mx171:evaluation} then evaluates these characters
on the testing ideals.
\end{proposition}
\begin{proof}
\noindent Smith normal form gives the coordinate maps for the presented
group. The Chinese remainder theorem recovers each coordinate
from its primary reductions, proving the common-kernel assertion.
\end{proof}
\noindent Thus prime-power \emph{orders} always suffice when the characters are
computed on $G_{\mathfrak f}^+$. Restricting to order $p$ would not suffice:
all maps $\Z_{p^2}\to\Z_p$ kill the subgroup $p\Z_{p^2}$.
The required orders are the full primary orders of the invariant factors.

\subsection{Assembly and Artin action}
\Needspace{10\baselineskip}
\begin{theorem}
% Former title: Generators of the full strict ray class field

\label{root170:equality}
\noindent For the characters $\chi_{j,p}$ of
Proposition~\ref{v221:computed-characters}, construct
$\theta_{\chi_{j,p}}$ by Theorem~\ref{root170:general-R}.
For ordinary anti-invariant characters one may instead use
Theorem~\ref{root170:ordinary-descent-theorem}. Then
\begin{equation}\label{v203:full-field}
 \boxed{K_{\mathfrak f}^+
        =K(\theta_{\chi_{j,p}}:p\mid d_j).}
\end{equation}
\end{theorem}
\begin{proof}
\noindent The characters have trivial common kernel. Apply
\eqref{v214:character-compositum} to their constructed cyclic fields.
\end{proof}
\noindent Thus the basis of the Kummer module, reciprocity selection
and descent give a field generator for each computed character.
Taking their compositum requires no further compatibility among
the chosen radicals.

\noindent For the Galois action, fix one character
$\chi:G_{\mathfrak f}^+\to\Z_n$ of exact order $n$ and put
\[
 \rho_\chi:G_{\mathfrak f}^+\longrightarrow\mu_n,
             \qquad\rho_\chi(g)=\zeta_n^{\chi(g)}.
\]
\noindent The map $\Z_n\to\mu_n$, $a\mapsto\zeta_n^a$, identifies these
two forms of the same character: $\chi$ gives the exponents in the
linear systems, and $\rho_\chi$ gives the root-of-unity multipliers.
Write the descended generator as
$\theta_\chi=\Theta_\chi(\alpha)$, with $\Theta_\chi(X)\in E[X]$
and $\alpha^n=\beta_\chi$ in $B_\beta$.
Identify the selected fixed algebra $A_{\beta,c}$ with $H_\chi$.
For $\xi\in\mu_n$, the $E$-automorphism
\[
 s_\xi:B_\beta\longrightarrow B_\beta,\qquad
                  \alpha\longmapsto\xi\alpha
\]
\noindent commutes with descent and restricts to a $K$-automorphism of $H_\chi$.

\begin{proposition}
% Former title: Artin action

\label{v211:artin-action}
\noindent For $\mathfrak a\in I_K(\mathfrak f)$,
\[
 \Art_{\mathfrak f}([\mathfrak a])|_{H_\chi}
       =s_{\rho_\chi([\mathfrak a])}|_{H_\chi}.
\]
\noindent Consequently
\begin{equation}\label{v211:artin-generator}
 \boxed{\Art_{\mathfrak f}([\mathfrak a])(\theta_\chi)
       =\Theta_\chi\bigl(\rho_\chi([\mathfrak a])\alpha\bigr).}
\end{equation}
\noindent For the characters in Theorem~\ref{root170:equality}, these formulas
give the action of $\Gal(K_{\mathfrak f}^+/K)$ on all the generators.
\end{proposition}
\begin{proof}
\noindent The selected descent has character $\chi$, so Artin reciprocity
acts through $r^{\chi([\mathfrak a])}=s_{\rho_\chi([\mathfrak a])}$
\cite[Part III, Theorems 7.8 and 9.2]{NeukirchBonn}.
Substitute in $\Theta_\chi$ to obtain the second formula.
\end{proof}
\noindent The polynomials in the ordinary and general constructions are,
respectively,
\begin{equation}\label{v211:ordinary-conjugates}
 \Theta_\chi(X)=\sum_{u\in\Z_n^\times}c_uX^u
\end{equation}
\noindent and
\begin{equation}\label{v211:conjugate-sums}
 \Theta_\chi(X)=
 \sum_{i=1}^{[E:K]}\sum_{k=0}^{n-1}\lambda_{ik}
 \sum_{\delta\in\Gal(E/K)}
 \delta(b_i)c_\delta^k X^{ku_\delta}.
\end{equation}
\noindent The $n$ values $\Theta_\chi(\xi\alpha)$, $\xi\in\mu_n$, are the
distinct conjugates of $\theta_\chi$, and their product polynomial
$\prod_{\xi\in\mu_n}(T-\Theta_\chi(\xi\alpha))$ lies in $K[T]$.

\subsection{What separate prime-power moduli can miss}
\label{mx171:primepower}\label{mx171:kernel}
\noindent Assume $\mathfrak f\ne(1)$ and put
$K_0=\prod_i K_{\mathfrak f_i}^+$, a compositum. The relevant
forgetful map and its kernel are
\[
 \pi_{\mathfrak f}:G_{\mathfrak f}^+\longrightarrow
        \prod_iG_{\mathfrak f_i}^+,\qquad
 [\mathfrak a]_{\mathfrak f}\longmapsto
           ([\mathfrak a]_{\mathfrak f_i})_i,
 \qquad J_{\mathfrak f}=\ker\pi_{\mathfrak f}.
\]
\noindent This kernel measures exactly the part of the full ray class field
missing from $K_0$.

\begin{proposition}
% Former title: The common unit relation

\label{v216:mixed-modulus}
\noindent For $U_i=\langle\bar\epsilon\rangle\subset R_i$ and $o_i=|U_i|$,
one has
\begin{equation}\label{mx171:gluingexact}
 \begin{gathered}
 1\longrightarrow J_{\mathfrak f}\longrightarrow G_{\mathfrak f}^+
 \xrightarrow{\pi_{\mathfrak f}}\operatorname{im}\pi_{\mathfrak f}
 \longrightarrow1,\\
 J_{\mathfrak f}\simeq
      (\prod_i U_i)/\langle(\bar\epsilon)_i\rangle.
 \end{gathered}
\end{equation}
\noindent The image consists of tuples with the same image in $\Cl^+(K)$.
Under Artin reciprocity,
\begin{equation}\label{mx171:fieldkernel}
 \begin{gathered}
 J_{\mathfrak f}=\Gal(K_{\mathfrak f}^+/K_0),\qquad
 [K_{\mathfrak f}^+:K_0]=\frac{\prod_i o_i}{\lcm_i o_i},\\
 K_0=K_{\mathfrak f}^+\quad\Longleftrightarrow\quad
                   o_i\text{ are pairwise coprime}.
 \end{gathered}
\end{equation}
\noindent For a character given by \eqref{mx171:unitcondition},
\begin{equation}\label{v217:local-character-criterion}
 H_\chi\subseteq K_0
       \quad\Longleftrightarrow\quad
          \ell_i(\bar\epsilon)=0\quad\text{for every }i.
\end{equation}
\end{proposition}
\begin{proof}
\noindent In the principal part of \eqref{v214:ray-exact}, the forgetting map is
\[
 (\prod_i R_i)/\langle(\bar\epsilon)_i\rangle
          \longrightarrow\prod_i(R_i/U_i),
\]
\noindent whose kernel is $J_{\mathfrak f}$. The ray maps all induce the identity
on $\Cl^+(K)$. Given ray classes with a common narrow ideal class,
divide by an ideal in that class and use the Chinese remainder theorem
to lift the resulting principal classes. This proves exactness.
The diagonal unit has order $\lcm_i o_i$, giving the degree formula;
the product equals the least common multiple exactly for pairwise
coprime orders. Finally $\chi$ kills $J_{\mathfrak f}$ precisely when
its residue characters kill every $U_i$. These are the usual ray-sequence
and character arguments \cite[Section~4.3.1]{CohenAdvanced}.
\end{proof}

\noindent The full group requires only $\sum_i\ell_i(\bar\epsilon)=0$;
the separate moduli require each summand to be zero.
This is why their compositum can miss characters, although
prime-power-order characters of the full group always suffice.

\noindent The missing group and the characters needed to detect it can both
be computed before the Kummer constructions. If there are $s$
prime-power factors of $\mathfrak f$, then
\[
 J_{\mathfrak f}\simeq
 \Z^s/\operatorname{im}
       \bigl(\operatorname{diag}(o_1,\ldots,o_s)\mid\boldsymbol1\bigr).
\]
\noindent Its Smith normal form gives the invariant factors of the missing
\emph{relative} extension. To locate this group inside
$G_{\mathfrak f}^+$, take totally positive CRT lifts
\[
 w_i\equiv\epsilon\pmod{\mathfrak f_i},\qquad
 w_i\equiv1\pmod{\mathfrak f_j}\quad(j\ne i).
\]
\noindent The classes $[(w_i)]_{\mathfrak f}$ generate $J_{\mathfrak f}$.
Their coordinates under $\varphi$ and the values
$\chi([(w_i)])=\ell_i(\epsilon)$ are obtained from
\eqref{mx171:evaluation}.

\begin{corollary}
% Former title: Adjoining the missing cyclic fields

\label{v218:missing-extension}
\noindent For a subset $\mathcal S$ of the primary coordinate characters in
Proposition~\ref{v221:computed-characters},
\[
 K_0(\theta_\chi:\chi\in\mathcal S)=K_{\mathfrak f}^+
 \quad\Longleftrightarrow\quad
 J_{\mathfrak f}\cap\bigcap_{\chi\in\mathcal S}\ker\chi=\{1\}.
\]
\noindent Such a set can be selected from the computed finite group before
constructing any $\theta_\chi$. Its character orders are the known
numbers $p^{v_p(d_j)}$ from the full ray-group presentation.
\end{corollary}
\begin{proof}
\noindent The subgroup fixing the field on the left is precisely the
intersection on the right. The full family of coordinate characters
has trivial common kernel, so it always contains a suitable subset.
\end{proof}
\noindent One concrete selection is to start with $J_0=J_{\mathfrak f}$,
choose a coordinate character $\chi_a$ nonzero on $J_{a-1}$,
and set $J_a=J_{a-1}\cap\ker\chi_a$. The order decreases at each
step until $J_a=\{1\}$. All these intersections are computations
in the finite presentation above. Afterwards the constructed fields
$K_a=K_{a-1}(\theta_{\chi_a})$ satisfy
\[
 \Gal(K_{\mathfrak f}^+/K_a)=J_a,\qquad
 [K_a:K_{a-1}]=|\chi_a(J_{a-1})|.
\]
\noindent The order of $\chi_a$ on the full group can exceed its order on
$J_{a-1}$. Thus the invariant factors of $J_{\mathfrak f}$ alone
do not specify the character orders needed for this construction.

\begin{example}
% Former title: A missing quadratic extension from an order-eight character

\label{mx171:example13}
\noindent Let $K=\Q(\sqrt{13})$ and $\epsilon=(11+3\sqrt{13})/2$.
Here $h_K^+=1$, and the unit has orders $6$ and $8$ at the inert
primes $5$ and $7$. Residue generators can be chosen so that
$\epsilon=g_5^4=g_7^6$. Writing $e_5,e_7$ for their ray classes gives
\[
 G_{(35)}^+\simeq
 \frac{\Z e_5\oplus\Z e_7}
      {\langle24e_5,48e_7,4e_5+6e_7\rangle}
 \simeq\Z_2\oplus\Z_{24},\qquad
 G_{(5)}^+\simeq\Z_4,\quad G_{(7)}^+\simeq\Z_6.
\]
\noindent The character $\chi(e_5)=1$, $\chi(e_7)=2$ in $\Z_8$ kills the
common unit relation since $4+6\cdot2=0$ in $\Z_8$.
Its value on $(\epsilon,1)$ is $4$, so it is nontrivial on the
order-two group $J_{(35)}$. Hence
\[
 [K_0:K]=24,\qquad[K_0H_\chi:K_0]=2,\qquad
                K_0H_\chi=K_{(35)}^+.
\]
\noindent Thus an order-eight character of the full modulus recovers the missing
quadratic extension over $K_0$. Every character to $\Z_2$ or
$\Z_4$ kills $J_{(35)}$: its generator $4e_5$ is already a
multiple of four. The order of the missing relative extension
therefore cannot replace the order of the required full-group character.
\end{example}

\subsection{The strict Hilbert class field}
\label{exc170:hilbert}\label{hi171:hilbert}\label{hi171:algorithm}
\noindent For $\mathfrak f=(1)$ the result is the classical strict Hilbert
class field $H_K^+$, the maximal abelian extension unramified at finite
primes \cite[Chapter VI, Proposition 6.8]{NeukirchANT}.
The following simplification concerns the norm-one representatives,
not the definition of that field. Compare the computational
constructions in \cite[Sections 3, 4.2--4.3]{FHS}
and \cite{CohenRoblotHilbert}.

\begin{proposition}
% Former title: Norm-one representatives for strict Hilbert characters

\label{hi171:hilbertcorrectness}\label{exc170:hilbertalgorithm}
\label{exc170:fullcyclo}
\noindent Let $n=p^a$, $E=\Q(\zeta_{\lcm(\Delta,n)})$, and let $\Sigma$
be finite and $\Gal(E/\Q)$-stable, containing the primes above $p$
and all primes of norm at most the Minkowski bound $B_E$.
For every character $\chi:\Cl^+(K)\to\Z_n$, one has
\[
                   \nu^{-1}(\chi|_{G_E})\in W^{(1)}_\Sigma.
\]
\noindent In particular no signature correction or extra dyadic coefficient
enlargement is needed.
\end{proposition}
\begin{proof}
\label{v203:hilbertsignature}
\noindent The relation
$\mathfrak a\tau_K\mathfrak a=(\Nm\mathfrak a)$ implies that
$\tau_K$ acts on $\Cl^+(K)$ by inversion. Reciprocity gives the
same action on $H_K^+$ \cite[Part III, Section 6]{NeukirchBonn}.
Complex conjugation lies in the $2$-torsion, which inversion fixes;
it is therefore central. Hence each cyclic field
$H\subset H_K^+$ is totally real or CM.

\noindent Put $L=HE$ and let $J$ be complex conjugation. For a Kummer
eigenvector $\alpha$ of the restricted character,
$A=\alpha J(\alpha)\in(E^+)^{\times,+}$.
The extension $L/E$ is unramified at finite primes. At an inert
unramified prime of $E/E^+$, the valuations of $\alpha$ at all
primes above it agree, since its conjugates differ by roots of
unity. Thus $A$ has even valuation and is a local norm.
It is also a norm at split primes and real places.
The cyclotomic inertia description
\cite[Chapter I, Proposition 10.3]{NeukirchANT} shows that $E/E^+$
has at most one finite ramified prime: $j$ can lie in inertia at
$\ell$ only if the prime-to-$\ell$ part of the conductor divides two.
The product
formula for local Artin symbols handles that possible last prime.
The Hasse norm theorem gives $b\in E^\times$ with $bj(b)=A$
\cite[Part III, Theorem 6.15 and Corollary 4.8]{NeukirchBonn}.
Consequently $(\alpha/b)^n$ is norm one. Since
$\Cl(R_\Sigma)=0$, Lemma~\ref{v203:supported-normone} puts its class in
$W^{(1)}_\Sigma$.
\end{proof}

\begin{corollary}
% Former title: Selection by discriminant one

\label{v216:hilbert-enumeration}
\noindent Enumerate the coordinate products $\prod_{j=0}^r z_j^{e_j}$,
$0\le e_j<n$, and their descents over the coefficient field above.
The compositum of the field factors $L$ satisfying
$\mathfrak d_{L/K}=\OO_K$ is
\begin{equation}\label{hi171:hilbertdiscriminant}
                         (H_K^+)^{n\Cl^+(K)}.
\end{equation}
\noindent If $\Cl^+(K)\simeq\bigoplus_j\Z_{m_j}$, take
$n=p^{\max_jv_p(m_j)}$ for each $p\mid\prod_jm_j$.
These exponents are known from its Smith normal form, and the
compositum of the resulting fields is $H_K^+$.
\end{corollary}
\begin{proof}
\noindent The proposition includes every strict Hilbert character of exponent
dividing $n$. Conversely the selected cyclic fields are unramified
at finite primes by the discriminant criterion
\cite[Chapter I, Section 9]{NeukirchANT}. Character recovery proves
the equality.
\end{proof}
\noindent This is the conductor test of Lemma~\ref{exc170:directconductor}
with $\mathfrak f=(1)$: no separate comparisons on ray generators
are needed when all unramified characters are included.
\section{Kummer modules and explicit generators: examples}
\label{root170:examples}
\noindent In each example we give the Kummer module, the class
selected by reciprocity and a generator obtained by descent.
The displayed bases come from the indicated RM curves. All roots
of unity use $\zeta_m=e^{2\pi i/m}$.

\noindent The symbols $\log_g(x)$ below are finite discrete logarithms:
$g^{\log_g(x)}=x$ in the stated finite field, with the exponent
taken modulo $|g|$. They specify the \emph{input} ray character
$\chi$, using the residue-group presentation of
Section~\ref{v227:full-modulus-characters}; they are not the Kummer
coordinates $e_j$. Selection computes those coordinates from
$C\boldsymbol e=\boldsymbol b$. We give the Jacobi expressions
producing each good-prime basis, followed by the selected class
and descended generator; the arithmetic evaluations are omitted.

\subsection{An ordinary quintic}
\label{v221:ordinary-quintic}
\noindent Let $K=\Q(h)$, $h=\sqrt2$, $E=K(\zeta_5)$,
$n=5$ and $\mathfrak f=(1025)$.
For a totally positive $u\in\mathcal O_K$ prime to $\mathfrak f$,
define $r,s,\ell_\pm\in\Z_5$ and $\chi$ by
\[
 \begin{gathered}
 u^{24}\equiv1+5(r+sh)\pmod{25},\qquad
 \ell_\pm=\log_6(u|_{h=\pm17}\bmod41)\pmod5,\\
 \chi([(u)])=4s+2\ell_++3\ell_-.
 \end{gathered}
\]
\noindent For $\Sigma$ the primes above $41$, use $C_{5,8}$.
On $\F_{41}^\times=\langle6\rangle$, let
$\vartheta(6)=\zeta_5$, $\xi(6)=\zeta_8$ and put
\[
 J_{i,c}=J(\vartheta^i,\xi^c),\qquad
 A_i=\frac{J_{i,1}J_{i,7}}{J_{i,3}J_{i,5}}\quad(i=1,2),
\]
\noindent with $J$ as in \eqref{geo170:jacobi}. These determinant quotients
already give a Kummer basis. Put
$s_0=\zeta_5+\zeta_5^{-1}$, $t_0=\zeta_5-\zeta_5^{-1}$. Then
\begin{equation}\label{v221:quintic-inputs}
 \begin{aligned}
 z_1=A_1&=(-11+16s_0+14ht_0+2hs_0t_0)/41,\\
 z_2=A_2&=(-27-16s_0-2ht_0+14hs_0t_0)/41,\\
 W_\Sigma&=\Z_5[A_1]\oplus\Z_5[A_2].
 \end{aligned}
\end{equation}
\noindent Selection by reciprocity gives
\begin{equation}\label{v221:quintic-beta}
                  (e_1,e_2)=(2,1),\qquad\beta_\chi=A_1^2A_2.
\end{equation}
\noindent For $\alpha^5=\beta_\chi$, descent gives
\begin{equation}\label{v221:quintic-theta}
 \theta_\chi=\alpha+\alpha^2/A_1+
                 \alpha^3/(A_1A_2)+\alpha^{-1}.
\end{equation}
\noindent Its minimal polynomial over $K$ is
\[
 X^5-10X^3+\frac{380}{41}X^2+\frac{185}{41}X-\frac{404}{1681}.
\]
\noindent The cyclic extension $H_\chi=K(\theta_\chi)$ is totally real,
with conductor $(1025)$ and relative discriminant $(1025)^4$.

\subsection{An ordinary cubic without saturation}
\label{v224:ordinary-cubic}
\noindent Let $K=\Q(h)$, $h=\sqrt2$, $E=K(\zeta_3)$,
$n=3$, $s=h(2\zeta_3+1)=\sqrt{-6}$ and $\mathfrak f=(63)$.
For a totally positive $u=a+bh$ prime to $\mathfrak f$,
define $r,t,\ell_\pm\in\Z_3$ and $\chi$ by
\[
 \begin{gathered}
 u^8\equiv1+3(r+th)\pmod9,\qquad
 \ell_\pm=\log_3(a\pm3b\bmod7)\pmod3,\\
 \chi([(u)])=t+\ell_+-\ell_-.
 \end{gathered}
\]
\noindent For $\Sigma$ the primes above $7$, use $C_{3,8}$ over
$k=\F_7(\varrho)$ with $\varrho^2=-1$ and
$k^\times=\langle2+\varrho\rangle$.
Let $\vartheta(2+\varrho)=\zeta_3$, $\xi(2+\varrho)=\zeta_8$
and put $J_c=J(\vartheta,\xi^c)$. The determinant quotient gives
\[
 A=\frac{J_1J_7}{J_3J_5}=\frac{1-20s}{49},\qquad
 W_\Sigma=\Z_3[A],\qquad z_1=A,\qquad
 e_1=1,\qquad\beta_\chi=A.
\]
\noindent Here $[U_\Sigma:\langle-1,A\rangle]=B=2$, so no internal
cube root is needed. For $\alpha^3=A$,
\[
 \theta_\chi=\alpha+\alpha^{-1},\qquad
 \theta_\chi^3-3\theta_\chi-\frac2{49}=0.
\]
\noindent The field $H_\chi=K(\theta_\chi)$ is totally real and cyclic cubic,
with conductor $(63)$ and relative discriminant $(63)^2$.

\subsection{A dyadic character on the RM curve \texorpdfstring{$C_{8,8}$}{C8,8}}
\label{v221:dyadic-example}
\noindent Let $K=\Q(h)$, $h=\sqrt2$, $E=\Q(\zeta_8)=K(i)$,
$n=4$ and $\mathfrak f=(17)$. For a totally positive generator
$a+bh$ of an ideal prime to $17$, define
\[
 \chi([(a+bh)])=\log_3(a+6b)-\log_3(a-6b)\pmod4,
\]
\noindent where the logarithms are in $\F_{17}^\times$.
Put $w=ih$, and let $\Sigma$ be the primes above $17$.
For $C_{8,8}$, choose $\lambda(3)=\zeta_8$ on $\F_{17}^\times$.
The Jacobi sums and positive RM quotients are
\[
 \begin{gathered}
 J_{2,1}=J(\lambda^2,\lambda)=1-4i,\qquad
 J_{7,7}=J(\lambda^7,\lambda^7)=-3+2w,\\
 A_1=J_{7,7}/j(J_{7,7}),\qquad
 A_2=J_{2,1}J_{7,7}/j(J_{2,1}J_{7,7}).
 \end{gathered}
\]
\noindent An internal square root gives
$z_1=-j(J_{2,1})J_{7,7}/17$ and $z_2=A_1/A_2$;
thus $z_1^2=A_1^2/A_2$. Explicitly,
\[
 \begin{gathered}
 z_0=\zeta_8,\qquad
 z_1=-\frac{(1+4i)(-3+2w)}{17},\qquad
 z_2=\frac{-15+8i}{17},\\
 W_\Sigma^{(1)}=\Z_4[z_0]\oplus\Z_4[z_1]\oplus\Z_4[z_2].
 \end{gathered}
\]
\noindent Selection gives $(e_0,e_1,e_2)=(0,2,-1)$ and
\[
 \beta_\chi=z_1^2z_2^{-1}=\frac{1-12w}{17}.
\]
\noindent For $\alpha^4=\beta_\chi$, the character selects the lift
$d_j(\alpha)=-\alpha^{-1}$. Hence
\[
 \theta_\chi=\alpha-\alpha^{-1},\qquad
 \theta_\chi^4+4\theta_\chi^2+\frac{32}{17}=0.
\]
\noindent The extension $H_\chi/K$ is cyclic quartic, with finite conductor
$(17)$ and relative discriminant $(17)^3$; both real places ramify.
The lift $d_j(\alpha)=\alpha^{-1}$ instead gives conductor
divisible by the prime above two and is excluded. This is why
descent uses the full character, not just its restriction to $E$.

\subsection{Shintani's example over \texorpdfstring{$\Q(\sqrt5)$}{Q(sqrt5)}}
\label{v206:five-example}\label{v206:five-bridge}
\noindent Let $K=\Q(d)$, $d=\sqrt5$, $\phi=(1+d)/2$,
$\epsilon=\phi^2$ and $\mathfrak p=(4-d)$, of norm $11$.
The group $G_{\mathfrak p}^+$ has order two; let $\chi$ be its
nontrivial character. Take $n=2$, $E=\Q(\zeta_5)$ and
$\Sigma$ the primes above $11$.
Use the quotient $C_{10,5}\to C_{5,5}$,
$(Y,Z)\mapsto(Y^2,Z)$; the coefficient field is unchanged since
$\Q(\zeta_{10})=E$. On $C_{5,5}$ take $\lambda(2)=\zeta_5$ on
$\F_{11}^\times$. The Jacobi values
$J_1=J(\lambda,\lambda)$ and $J_2=J(\lambda^2,\lambda)$ are
\[
 J_1=-\zeta_5^3-2\zeta_5^2+2\zeta_5,\qquad
 J_2=2\zeta_5^3+4\zeta_5^2+\zeta_5+2
\]
\noindent and give the basis
\[
 z_0=-\zeta_5,\quad z_1=J_1J_2/11,\quad z_2=J_1/j(J_1),
 \qquad W_\Sigma^{(1)}=\bigoplus_{j=0}^2\Z_2[z_j].
\]
\noindent The corrected character $\psi$ has coordinates $(1,1,0)$.
After restoring its signature factor,
\begin{equation}\label{v218:five-radicand}
 \beta_\chi=d z_0z_1\equiv3\phi-2\pmod{E^{\times2}}.
\end{equation}
\noindent Descent gives
\begin{equation}\label{v218:five-ray-field}
 \theta_\chi=\sqrt{3\phi-2},\qquad
                  K(\theta_\chi)=H_\chi=K_{\mathfrak p}^+.
\end{equation}
\noindent Its relative discriminant is $\mathfrak p$.
With the square root positive at $\sigma_1$, Shintani's
evaluated double-sine invariant is
\begin{equation}\label{v218:shintani-five}
 X_{\mathrm{Sh}}=\frac{\epsilon+\sqrt{3\phi-2}}2,
 \qquad K(X_{\mathrm{Sh}})=K_{\mathfrak p}^+.
\end{equation}
\noindent See \cite[Section~3, Example~2, pp.~191--192]{Shintani1977}
and \cite[Section~2, (2.1)--(2.2)]{Yamamoto2007}.
The minimal polynomial of $X_{\mathrm{Sh}}$ over $\Q$ is
$X^4-3X^3+3X^2-3X+1$.

\subsection{An RM input at a bad reduction prime}
\label{v222:bad-prime-example}
\noindent Let $K=\Q(h)$, $h=\sqrt6$, $n=3$ and $\mathfrak f=(9)$.
The ray class group is $\Z_6\oplus\Z_3$, with generators
$[(-4+10h)]$, $[(17-2h)]$; let $\chi$ take value one on both.
Put
\[
 E=K(\zeta_3)=\Q(s,t),\quad s^2=-2,\quad t^2=-3,
 \quad h=-st,\quad\zeta_3=(-1+t)/2,\quad b=\frac{1+s}{1-s},
\]
\noindent with $s,t$ of positive imaginary part, and let $\Sigma$ be the
primes above three. Here the RM curve $C_{24,24}$ is evaluated
through its CM ideal, as in Section~\ref{v222:principalized-input}.
The principalization can be written out completely. In
$F=\Q(w)$, $w=\zeta_{24}$, identify $s=w^3+w^9$, and denote
the two primes above three by $\mathfrak q_+,\mathfrak q_-$,
where $v_{\mathfrak q_+}(1+s)=2$, and choose
$\mathfrak q_3=\mathfrak q_+$; then $\mathfrak q_-=j(\mathfrak q_+)$.
For $v=3$, the positive
RM indices are $I_3^+=\{1,5,19,23\}$. The CM-type formula gives
\[
 H_{3,3}=12[\mathfrak q_+]+4[\mathfrak q_-],\qquad
 \mathfrak a_{3,3}=\mathfrak q_+^{12}\mathfrak q_-^4
                          =(9(1+s)^4).
\]
\noindent Indeed $v_{\mathfrak q_\pm}(3)=2$, while $1+s$ has valuations
$2,0$. Thus its ideal class order is one. Choosing $x=9(1+s)^4$ gives
\begin{equation}\label{v228:bad-prime-calculation}
 \frac{x}{j(x)}=b^4,\qquad
 A_1=N_{F/E}\left(\frac{x}{j(x)}\right)=b^8,
                         \qquad [F:E]=2.
\end{equation}
\noindent This is an explicit algebraic input; no good-reduction extension
or bad-prime Frobenius value has been computed. It gives
\[
 z_0=\zeta_3,\qquad z_1=A_1,\qquad
 W_\Sigma^{(1)}=\Z_3[z_0]\oplus\Z_3[z_1].
\]
\noindent Here $B=16$, so this input requires no internal cube root.
Its relation to the geometric Frobenius quotient is
\eqref{v222:principalized-RM}.
Selection gives $(e_0,e_1)=(0,2)$ and $\beta_\chi=b^{16}$.
For $\alpha^3=\beta_\chi$, put $\eta=\alpha/b^5$. Then
\[
 \eta^3=b,\qquad\theta_\chi=\eta+\eta^{-1},\qquad
 \theta_\chi^3-3\theta_\chi+\frac23=0.
\]
\noindent The field $H_\chi=K(\theta_\chi)$ is totally real and cyclic
cubic, with conductor $(9)$ and relative discriminant $(81)$.

\noindent For comparison, the alternative of
Section~\ref{v222:diagonal-alternative} uses $C_8$ at the same
prime three, where this curve has good reduction. In
$\F_9=\F_3(\varrho)$, $\varrho^2=-1$, choose
$\lambda(2+\varrho)=\zeta_8$. Since
$E\cap\Q(\zeta_8)=\Q(s)$, direct Jacobi evaluation gives
\[
 J(\lambda,\lambda)=(1-s)^2=-1-2s,\qquad
 A^{\rm diag}=N_{\Q(\zeta_8)/\Q(s)}
       \left(\frac{J(\lambda,\lambda)}{j(J(\lambda,\lambda))}\right)
                         =b^{-4}.
\]
\noindent Thus $A_1=(A^{\rm diag})^{-2}$ and
$[A_1]=[A^{\rm diag}]$ in $E^\times/E^{\times3}$.
Both inputs give the displayed Kummer basis and the same selected
extension. The RM route computes a principal ideal generator at
level $24$; the diagonal route evaluates a Jacobi sum at level $8$.

\appendix

\section{Computing ray characters from ideals and units}
\label{v207:ray-character-appendix}

\noindent This standard algorithm uses the ordinary class group and both
real signs, giving the same ray group as the narrow-class
presentation \eqref{v214:ray-exact}
\cite[Section 4.3.1]{CohenAdvanced}\cite[Section 3]{FHS}.
It fixes the residue and sign conventions used in the construction.

\begin{proposition}
% Former title: Computing the ray characters

\label{v204:arithmetic-input-proposition}
\noindent Arithmetic in $K$ gives a finite presentation of $G_{\mathfrak f}^+$
and all homomorphisms $G_{\mathfrak f}^+\to \Z_n$, with their
values on ideals, without constructing a class field. More explicitly:
\begin{enumerate}
\item Choose an invariant-factor presentation of the ordinary ideal
class group, with representatives $\mathfrak a_1,\ldots,\mathfrak a_s$
prime to $\mathfrak f$, of orders $d_i$, and find
$\mathfrak a_i^{d_i}=(\alpha_i)$. Compute a fundamental unit
$\epsilon_0$, so $\OO_K^\times=\langle-1,\epsilon_0\rangle$.
\item Compute the finite abelian group and the homomorphism
\[
 A_{\mathfrak f}=(\OO_K/\mathfrak f)^\times\times\{\pm1\}^2,
\]
\[
\begin{gathered}
 a:\{x\in K^\times:v_{\mathfrak p}(x)=0
                 \ (\mathfrak p\mid\mathfrak f)\}\longrightarrow A_{\mathfrak f},\\
 x\longmapsto(x\bmod\mathfrak f,\operatorname{sgn}\sigma_1(x),
                             \operatorname{sgn}\sigma_2(x)).
\end{gathered}
\]
\noindent For $\mathfrak f=\OO_K$ the first factor is the trivial group.
Take generators and defining relations of $A_{\mathfrak f}$, add
symbols $e_i$ for $\mathfrak a_i$, and impose the further relations
\[
                  a(-1)=a(\epsilon_0)=0,\qquad
                         d_i e_i=a(\alpha_i).
\]
\noindent Here the finite abelian group is written additively. The resulting
group is $G_{\mathfrak f}^+$.
\item Let $R\in\operatorname{Mat}_{t\times s}(\Z)$ have these
relations as columns, in a list of $t$ generators. Thus
$G_{\mathfrak f}^+\simeq\Z^t/R\Z^s$.
A character is exactly a row vector
$\boldsymbol c\in(\Z_n)^t$ satisfying
$\boldsymbol cR=0$, and, for a column vector $\boldsymbol x$,
\begin{equation}\label{v204:character-from-relations}
 \chi_{\boldsymbol c}([\boldsymbol x])
                =\boldsymbol c\boldsymbol x\pmod n.
\end{equation}
\noindent Smith normal form lists these characters as in
Proposition~\ref{v221:computed-characters}.
\item To evaluate it on $\mathfrak a\in I_K(\mathfrak f)$, find
$m_i$ and $\alpha\in K^\times$ such that
\[
             \mathfrak a=(\alpha)\prod_i\mathfrak a_i^{m_i}.
\]
\noindent The coordinate vector is that of $a(\alpha)+\sum_i m_i e_i$ in
the displayed presentation. Apply the preceding formula.
\end{enumerate}
\end{proposition}
\begin{proof}
\noindent The standard ray exact sequence is
\[
 A_{\mathfrak f}/a(\OO_K^\times)\lhook\joinrel\longrightarrow
 G_{\mathfrak f}^+\longrightarrow\Cl(\OO_K)\longrightarrow0;
\]
\noindent see \cite[Section 4.3.1]{CohenAdvanced}. Lifting ideal-class generators gives
$d_i e_i=a(\alpha_i)$ with our principal-ideal convention.
The class-group, unit, and principal-ideal
computations are standard \cite[Chapters 5--6]{CohenANT}.
\end{proof}

\section{The CM interpretation of the Frobenius ideal}
\label{geo203:CM-appendix}
\noindent The construction uses the local CM ideal formula in
Proposition~\ref{v222:RM-divisor}; the proof there is independent
of this appendix. We record its relation with the classical
reflex norm, which is different from the coefficient norm
$N_{F/E}$ in \eqref{v222:RM-input}.

\noindent For the ordinary factor primitive in both deck indices, put
$F=\Q(\zeta_{\Delta n})$. Its CM type is
\begin{equation}\label{geo203:CM-Phi}
 \Phi=\{\sigma_u:u\in\Z_{\Delta n}^\times,
                     \{u/n\}+\{u/\Delta\}<1\}.
\end{equation}
\noindent Here an automorphism of the curve acts on the Jacobian by
pushing forward degree-zero divisors. Its action on the tangent
space has the eigenvalues of its pullback on holomorphic
differentials. Those differentials give \eqref{geo203:CM-Phi}
\cite[pp.~1183--1184]{KoblitzRohrlich}; the factor need not be
simple. The uniform factor has the CM types
\eqref{v222:RM-CM-type}, with the same action convention.

\noindent For this abelian CM field, define
\[
 S_\Phi=\{s\in\Gal(F/\Q):s\Phi=\Phi\},\qquad
 F^*=F^{S_\Phi},\qquad
 \Phi^*=\{\sigma^{-1}|_{F^*}:\sigma\in\Phi\},
\]
\noindent using distinct restrictions in the last set. The reflex norm is
\[
 N_{\Phi^*}:(F^*)^\times\longrightarrow F^\times,
                 \qquad a\longmapsto\prod_{\psi\in\Phi^*}\psi(a),
\]
\noindent with the analogous definition on ideals. At a prime
$\mathfrak q\nmid\Delta n$ of $F$, the Jacobi factorization is
\begin{equation}\label{v213:reflex-ideal}
 N_{\Phi^*}(N_{F/F^*}\mathfrak q)\mathcal O_F
       =\prod_{\sigma\in\Phi}\sigma^{-1}\mathfrak q
       =(J_{\mathfrak q,1}(1)).
\end{equation}
\noindent Indeed $N_{F/F^*}\mathfrak q$ extends to
$\prod_{h\in S_\Phi}h\mathfrak q$; the reflex type chooses
exactly the $S_\Phi$-cosets in $\Phi^{-1}$. This proves the
middle equality with all prime multiplicities, without assuming
complete splitting or a primitive CM type. It is the abelian
form of the Shimura--Taniyama ideal formula
\cite[Chapter~II, Section~8; Chapter~III, Section~13]{ShimuraCM}.
At ramified primes we use its local form
\eqref{v223:local-CM-slope}, as explained in
Section~\ref{v223:general-CM-divisors}.

\noindent Thus the reflex norm describes the Frobenius \emph{ideal} from
the CM type. The coefficient norm $N_{F/E}$ instead puts a
normalized determinant value in the chosen Kummer coefficient
field. At good primes the normalized determinant is a quotient of Jacobi products;
at bad primes Section~\ref{v222:principalized-input} computes
substitute values by principalizing the CM ideals. No such
substitute is asserted to be a Frobenius eigenvalue.

\section{Classical growth of ray towers}
\label{mx201:towers}
\noindent Fix a nonzero integral ideal $\mathfrak f$, and write
\[
 S=\{\mathfrak p:\mathfrak p\mid\mathfrak f\},\qquad
 K_a=K_{\mathfrak f^a}^+,\qquad
 K_\infty=\bigcup_{a\ge1}K_a.
\]
\noindent Thus $K_\infty$ is the maximal abelian extension of $K$ unramified
at every finite prime outside $S$. We ask which part of this tower
can keep growing when the support $S$ is fixed.

\noindent Say that $S$ contains the \emph{full support above $p$} if every
prime of $K$ above the rational prime $p$ belongs to $S$. Put
\begin{equation}\label{v227:cyclotomic-support}
 \begin{gathered}
 T=\{p:\text{$S$ contains the full support above $p$}\},\\
 K_T^{\mathrm{cyc}}=K\bigl(\Q_\infty^{(p)}:p\in T\bigr),
 \end{gathered}
\end{equation}
\noindent where $\Q_\infty^{(p)}$ is the cyclotomic
$\widehat\Z_p$-extension of $\Q$; put $K_T^{\mathrm{cyc}}=K$
when $T$ is empty. An inert or ramified prime in $S$ supplies full
support above its rational prime. If $p$ splits, both primes above
$p$ must lie in $S$.

\noindent The following is a classical consequence of the structure of abelian
$S$-ramification
\cite[Section~1.3, equations~(1.1) and~(1.3), pp.~252--253]{GrasIncomplete2019}
and of cyclotomic $\widehat\Z_p$-extensions
\cite[Section~13.1]{Washington}. It is independent of the Fermat construction.

\begin{theorem}
% Former title: Eventual cyclotomic growth
\label{mx201:stabilization}
\noindent The extension $K_\infty/K_T^{\mathrm{cyc}}$ is finite. There is an
integer $a_0\ge1$ such that
\begin{equation}\label{mx201:cyclotomic-tail}
 \boxed{\quad K_\infty=K_{a_0}K_T^{\mathrm{cyc}},\qquad
 K_a=K_{a_0}\bigl(K_a\cap K_T^{\mathrm{cyc}}\bigr)
                       \quad(a\ge a_0).\quad}
\end{equation}
\noindent The tower is eventually constant exactly when $T$ is empty.
\end{theorem}
\noindent In particular, after one finite strict ray class field $K_{a_0}$
is known, every further layer requires only cyclotomic adjunctions.
For $\mathfrak f=m\OO_K$, all rational primes dividing $m$ have
full support, so $T=\{p:p\mid m\}$. In contrast, a tower supported
at just one prime above a split rational prime eventually stops.

\begin{proof}
\noindent Put $G=\Gal(K_\infty/K)$, and let $G_p$ be its maximal pro-$p$
quotient. The cited structure and rank formulas give a finitely
generated $\widehat\Z_p$-module with
\begin{equation}\label{mx201:rank}
 \operatorname{rank}_{\widehat\Z_p}G_p
 =\begin{cases}1,&p\in T,\\0,&p\notin T.\end{cases}
\end{equation}
\noindent Indeed, the rank is the local degree
$d_p=\sum_{\mathfrak p\in S,\,\mathfrak p\mid p}[K_{\mathfrak p}:\Q_p]$
minus the dimension of the span of the unit logarithms. This dimension
is one when $d_p>0$: the fundamental unit cannot be a local root of
unity. Since $K$ is real quadratic, $d_p$ is two exactly when $p\in T$.
Gras uses unramified real places; allowing their ramification adds only
a finite $2$-group, so the same conclusion holds here, including $p=2$.

\noindent The ray sequence~\eqref{v214:ray-exact} also shows that only finitely many
$G_p$ are nontrivial: outside the rational primes below $S$, their orders
divide the product of the narrow class number and the orders of the
residue-unit groups at $S$. For $p\in T$, the cyclotomic
$\widehat\Z_p$-extension of $K$ lies in $K_\infty$. Restriction to it is
a surjection from $G_p$ onto $\widehat\Z_p$, with finite kernel
by~\eqref{mx201:rank}. The other $G_p$ are finite. Consequently
$H=\Gal(K_\infty/K_T^{\mathrm{cyc}})$ is finite.

\noindent Let $N_a=\Gal(K_\infty/K_a)$. Since $\bigcap_aN_a=1$, choose
$a_0$ with $N_{a_0}\cap H=1$. Galois correspondence gives
$K_\infty=K_{a_0}K_T^{\mathrm{cyc}}$. For $a\ge a_0$, the elementary
identity $N_{a_0}\cap HN_a=N_a$ gives the second formula in
\eqref{mx201:cyclotomic-tail}. Finally, $K_T^{\mathrm{cyc}}/K$ is
infinite exactly when $T\ne\varnothing$.
\end{proof}
\noindent The support is fixed throughout. The initial layer $K_{a_0}$ depends
on that support, and the theorem does not specify a minimal $a_0$.
\begingroup
\raggedbottom
\linespread{1}\selectfont
\let\savedthebibliography\thebibliography
\def\thebibliography#1{\savedthebibliography{#1}%
  \setlength{\itemsep}{0pt}\setlength{\parsep}{0pt}}
\bibliographystyle{plain}
\bibliography{FERMAT_RAY_RECIPROCITY_V232}
\enlargethispage{4\baselineskip}
\endgroup
\end{document}